\documentclass[preprint]{siamart}

\usepackage{booktabs}
\usepackage{fec}
\usepackage{algorithm}
\usepackage[noend]{algpseudocode}
\usepackage{xcolor}
\usepackage{footnote,longtable}
\hypersetup{colorlinks}

\usepackage{amsmath,amssymb,latexsym}%,amsthm}
\usepackage{graphicx}
\usepackage{epsf}
\usepackage{subfigure}
\usepackage{amsfonts}
\usepackage{graphics}

\usepackage{enumerate}

 \newtheorem{assu}[theorem]{Assumption}

 \numberwithin{equation}{section}

\newcommand{\beq}{\begin{equation}}
\newcommand{\eeq}{\end{equation}}

\newcommand{\bB}{\mathbb{B}}

\newcommand{\bR}{\mathbb{R}}

\newcommand{\R}{\bR}
\newcommand{\Rn}{\bR^n}
\newcommand{\Rm}{\bR^m}

\newcommand{\eR}{{\bar\bR}}

\newcommand{\argmax}{\mathop{\mathrm{argmax}}}

\newcommand{\grad}{\nabla}
\newcommand{\half}{\frac{1}{2}}
\newcommand{\thalf}{\tfrac{1}{2}}

\newcommand{\norm}[1]{\left\Vert #1\right\Vert}
\newcommand{\tnorm}[1]{\norm{#1}_2}

\newcommand{\dist}[2]{\mathrm{dist}\left(#1\,\left|\, #2\right.\right)}

\newcommand{\tdist}[2]{\mathrm{dist}_2\left(#1\,\left|\, #2\right.\right)}

\newcommand{\diag}{\mathrm{diag}}
\newcommand{\interior}[1]{\mathrm{int}\left(#1\right)}

\newcommand{\gam}{\gamma}

\newcommand{\sig}{\sigma}

\usepackage[figuresright]{rotating}

\allowdisplaybreaks

\title    {Inexact Sequential Quadratic Optimization with Penalty Parameter Updates Within the QP Solve: Extended Version}
\author   {James V.~Burke\thanks{Dept.~of Mathematics, University of Washington; \email{jvburke@uw.edu}}
      \and Frank E.~Curtis\thanks{Dept.~of Ind.~and Sys.~Engr., Lehigh University; \email{frank.e.curtis@gmail.com}}
      \and Hao Wang\thanks{Sch.~of Inf.~Sci.~and Tech., ShanghaiTech University; \email{wanghao1@shanghaitech.edu.cn}}
      \and Jiashan Wang\thanks{Dept.~of Mathematics, University of Washington; \email{jsw1119@math.washington.edu}.}
          }
\date     {\today}

\begin{document}

% Title
\maketitle

%**********
% Abstract
%**********
\begin{abstract}
This paper focuses on the design of sequential quadratic optimization (commonly known as SQP) methods for solving large-scale nonlinear optimization problems.  The most computationally demanding aspect of such an approach is the computation of the search direction during each iteration, for which we consider the use of matrix-free methods.  
%In particular, to reduce overall computational costs, we consider the use of matrix-free methods in such a way that only inexact subproblem solutions are required 
In particular, we develop a method that requires an inexact solve of a single QP subproblem
to establish the convergence of the overall SQP method.  It is known that SQP methods can be plagued by poor behavior of the global convergence mechanism.  To confront this issue, we propose the use of an exact penalty function with a dynamic penalty parameter updating strategy to be employed \emph{within} the subproblem solver in such a way that the resulting search direction predicts progress toward both feasibility and optimality.  We present our parameter updating strategy and prove that, under reasonable assumptions, the strategy does not modify the penalty parameter unnecessarily.  We also discuss a matrix-free subproblem solver in which our updating strategy can be incorporated.  We close the paper with a discussion of the results of numerical experiments that illustrate the benefits of our proposed techniques.
\end{abstract}

%**********
% Keywords
%**********

%*****
% AMS
%*****

%*********
% Section
%*********
\section{Introduction}\label{sec.introduction}

In this paper, we consider the use of sequential quadratic optimization (commonly known as SQP) methods for solving large-scale nonlinear optimization problems (NLPs) \cite{Burk89,Burk92,BurkCurtWang14, BurkHan89, Han77,Powe78,Wils63}.  While they have proved to be effective for solving small- to medium-scale problems, SQP methods have traditionally faltered in large-scale settings due to the expense of (accurately) solving large-scale quadratic subproblems (QPs) during each iteration.  However, with the use of matrix-free methods for solving these subproblems, one may consider the acceptance of inexact subproblem solutions.  Such a feature offers the possibility of terminating the subproblem solver early, perhaps well before an accurate solution has been computed.  This characterizes the type of strategy that we propose in this paper.

Some work has been done to provide global convergence guarantees for SQP methods that allow inexact subproblem solves \cite{CurtJohnRobiWaec14}.  However, the practical efficiency of such an approach remains an open question.  A critical aspect of their implementation is the choice of subproblem solver since it must be able to provide good inexact solutions quickly, as well as have the ability to compute highly accurate solutions---say, by exploiting well-chosen starting points---in the neighborhood of a solution of the NLP.  In addition, while a global convergence mechanism such as a merit function or filter is necessary to guarantee convergence from remote starting points, any NLP algorithm can suffer when such a mechanism does not immediately guide the algorithm toward promising regions of the search space.  To confront this issue when an exact penalty function is used as a merit function, we propose a dynamic penalty parameter updating strategy to be incorporated \emph{within} the subproblem solver so that each computed search direction predicts progress toward both feasibility and optimality.  This strategy represents a stark contrast to previously proposed techniques that only update the penalty parameter after a sequence of iterations, in hindsight at the end of an iteration \cite{Burk89,Han77,HanMang79}, or at the expense of numerous subproblem solves within a single iteration \cite{BurkCurtWang14,ByrdLopeNoce12,ByrdNoceWalt08}. 

 To provide some context about how the algorithm proposed in this paper compares to other recently proposed SQP-type methods in the literature, let us contrast our approach with those proposed in \cite{BurkCurtWang14} and \cite{CurtJohnRobiWaec14}.  The penalty SQP method proposed in \cite{BurkCurtWang14} was motivated by the desire to formulate an SQP approach that attains strong convergence guarantees when solving problems regardless of whether they involve constraints that are feasible or infeasible.  Toward this end, the approach involved a novel dynamic updating scheme for the penalty parameter that, e.g., quickly drives the algorithm toward constraint violation minimization when infeasibility is detected.  The approach relies on exact solves of two QP subproblems per iteration; the first determines the reduction that can be obtained in a local model of an infeasibility measure while the second minimizes a local model of the objective while ensuring that the reduction in a local model of the infeasibility measure is proportional to that attained by the solution of the first QP.  In this manner, rapid convergence can be attained when solving either a feasible or infeasible problem, although a high price is paid by needing exact subproblem solutions.  The method in \cite{CurtJohnRobiWaec14} overcomes this obstacle by allowing inexact subproblem solves.  However, it also potentially requires (approximate) solutions of two QPs per iteration, one aimed at minimizing constraint violation and one aimed at reducing the objective subject to an appropriate bound on constraint violation.  The approach proposed in this paper also allows inexactness in the QP solves, but only requires solving a \emph{single} QP in each iteration.  This is made possible by a new strategy for dynamically updating the penalty parameter \emph{within} the QP solver.  
This dynamic penalty parameter updating strategy is the focus of our investigation.
We prove that our algorithm does not reduce the penalty parameter unnecessarily and that one can ensure convergence to an optimal solution (when a given problem is feasible) or to an infeasible stationary point (when a given problem is infeasible).

Overall, the contributions in this paper can be summarized as the following.
\bitemize
  \item Our proposed SQP technique is specifically designed to be effective in large-scale settings.  In particular, it allows for the use of iterative methods for solving the QP subproblems, allowing inexactness in the subproblem solves.
  \item Our technique involves a dynamic penalty parameter updating strategy to be employed \emph{within} the subproblem solve.  This makes the approach efficient while not having to accurately solve multiple QPs in a single iteration.
  \item By ensuring that each computed step predicts progress toward minimizing constraint violation, our technique allows for automatic infeasibility detection.
\eitemize

%************
% Subsection
%************
\subsection{Organization}

In the remainder of this section, we outline our notation and introduce various concepts that will be employed throughout the paper.  In \S\ref{sec.SQP}, we introduce a basic penalty-SQP algorithm. Our penalty parameter updating strategy is detailed in \S\ref{sec.penalty}.  A complete algorithm is presented and analyzed in \S\ref{sec.complete_algorithm}. The results of numerical experiments are presented in \S\ref{sec.numerical}. Concluding remarks are provided in \S\ref{sec.conclusion}.

%************
% Subsection
%************
\subsection{Notation}

Let $\Rn$ be the space of real $n$-vectors, $\Rn_+$ be the nonnegative orthant of $\Rn$ (i.e., $\Rn_+ := \{x \in \Rn : x \geq 0\}$), and $\Rn_{++}$ be the interior of $\Rn_{+}$ (i.e., $\Rn_{++} := \{x \in \Rn : x > 0\}$).  The set of $m\times n$ real matrices is denoted $\R^{m\times n}$.  On $\R^{n}$, the $\ell_2$ (i.e., Euclidean) norm is indicated as $\tnorm{\cdot}$, with the unit $\ell_2$-norm ball defined as $\bB_2 := \{x \in \Rn : \|x\|_2 \leq 1\}$.  For a pair of vectors $(u,v) \in \Rn\times\Rn$, their inner product is written as $\langle u, v \rangle := u^Tv$ and the line segment between them is written as $[u,v]$.  The middle value operator applied to $(a,b,c) \in \R{} \times \R{} \times \R{}$, denoted by $\text{mid} \{a,b,c\}$, returns the median of $\{a,b,c\}$.  For a scalar $a$, let $(a)_+ := \max\{a,0\}$ and $(a)_- := \min\{a,0\}$.  The set of nonnegative integers is denoted by $\N{}$. The extended real number line is defined as  $\bar{\mathbb{R}}=\mathbb{R}\cup\{-\infty,+\infty\}$. 

For a set of scalars $b_i \in \R$ for $i \in \{1,\dots,m\}$, we denote the vector $\bmbf = [b_1,b_2,\dots,b_m]^T \in \Rm$.  For convenience, we use $\pmb{1}_n$ to denote the $n$-vector of all ones and $\pmb{0}_n$ to denote the $n$-vector of all zeros.  Given vectors $y^i \in \R^{d_i}$ for $i \in \{1,\dots,m\}$, we use boldface to denote the element $\ymbf = (y^1,\dots,y^m)$ on the product space $\R^{d_1} \times \cdots \times \R^{d_m}$.  Conversely, given $\ymbf \in \R^{d_1} \times \cdots \times \R^{d_m}$, the $i$-th component of $\ymbf$ (an element of $\R^{d_i}$) is denoted $y^i$ while the $j$-th element of $y^i$ is written as $y_j^i$. 
% In this product space, we define the norm % and its dual norm
%\bequationn
%  \norm{\ymbf} = \norm{(y^1, \dots, y^m)} := \sum_{i=1}^m \tnorm{y^i}%\ \ \text{and}\ \ \dnorm{\ymbf} = \sup_{i\in\{1,\dots,m\}} \tnorm{y^i}.
%\eequationn
For convex sets $C_i \in \R^{d_i}$ for $i \in \{1,\dots,m\}$, % we define the product set
%\bequationn
%  \Cmbf := C_1 \times \cdots \times C_m \subset \R^{d_1} \times \cdots \times \R^{d_m}.
%\eequationn
%T
the distance functions are defined as 
\bequationn
  \tdist{y^i}{C_i} := \inf_{z^i \in C_i} \|y^i - z^i\|_2.%\ \ \text{and}\ \ \dist{\ymbf}{\Cmbf} := \sum_{i=1}^m \tdist{y^i}{C_i},
\eequationn
The interior of a set $C$ is denoted by $\text{int}(C)$. 

For an extended-real-valued function $f : \Rn \to \eR$, the Legendre-Fenchel conjugate of $f$ is denoted as $f^\star$.  For a convex set $X \subseteq \Rn$, we define the characteristic function $\delta(d|X)$ which evaluates to 0 if $d \in X$ and evaluates to $\infty$ otherwise.  The conjugate of $\delta(\cdot|X)$ is the support function of $X$, which we denote by
%\bequationn
  $\delta^*( y |X) = \sup_{d \in X}\ \langle y, d \rangle$.
%\eequationn
For example, for a hyperplane $C := \{ d : \langle a, d \rangle + b = 0\}$ (respectively, half space $C=\{ d : \langle a, d\rangle + b\le 0\}$), one finds that $\delta^*(y|C) < \infty$ if and only if $\langle y, a\rangle= \pm \|y\|_2\|a\|_2$ (respectively, $\langle y, a\rangle=  \|y\|_2\|a\|_2$).  In this case, 
\begin{equation}\label{delta*}
  y = \zeta a\ \ \text{with}\ \ \zeta = \frac{1}{\|a\|_2^2} \langle y, a \rangle,\ \ \text{meaning that}\ \ \delta^*(y|C) = -\zeta b.
\end{equation}

For an iterative algorithm,  we use superscript $k$ to indicate the iteration number for vectors and subscript $k$ for scalars to avoid confusion with the $k$th power of the scalar, e.g., $x^k$ and $\rho_k$.  For an algorithm for solving the subproblem, we use superscript $(j)$ to indicate the iteration number for vectors and subscript $(j)$ for scalars.  

% If $f$ is convex, then the subdifferential of $f$ at $x$ is the set
%\bequationn
%  \partial f(x) := \{ y \in \Rn : f(x) + \langle y, z-x \rangle \leq f(z)\ \text{for all}\ z \in \bB_2\}.
%\eequationn
%For example, the subdifferentials of our distance functions are given by (see \cite{RW98})
%\bequalin
%  \partial \tdist{y^i}{C_i} &:= \bcases \frac{(I-P_{C_i})y^i}{\|(I-P_{C_i})y^i\|_2} & \text{if $y^i \notin C_i$} \\ \bB_2 \cap N(y^i|C_i) & \text{if $y^i \in C_i$,} \ecases \\
%  \text{and}\ \ \partial \dist{\ymbf}{C} &:= \partial \tdist{y^1}{C_1} \times \cdots \times \partial \tdist{y^m}{C_m}.
%\eequalin
%\bequalin
%  \partial \tdist{y_i}{C_i} &:= \bcases \frac{(I-P_{C_i})y_i}{\|(I-P_{C_i})y_i\|_2} & \text{if $i \notin \cA(\ymbf)$} \\ \bB_2 \cap N(y_i|C_i) & \text{if $i \in \cA(\ymbf)$,} \ecases \\
%  \text{and}\ \ \partial \dist{\ymbf}{C} &:= \partial \tdist{y_1}{C_1} \times \cdots \times \partial \tdist{y_m}{C_m},
%\eequalin
%where the index set $\cA(\ymbf)$ is defined as
%\bequationn
%  \cA(\ymbf) := \{i \in \{1,\ldots,m\} : \tdist{y_i}{C_i} = 0\}
%\eequationn
%The normal cone to $C_i$ at $y^i \in C_i$ is defined by
%\bequationn
%  N(y^i | C_i) := \{z^i \in \R^{d_i} : \langle z^i, p^i-y^i\rangle \leq 0\ \text{for all}\ p^i \in C_i\}.
%\eequationn

%*********
% Section
%*********
\section{A Penalty-SQP Framework}\label{sec.SQP}

%In this section, we formulate our problem of interest and outline the basic components of a penalty-SQP algorithm \cite{Flet87}.  This method represents the framework in which we will define our dynamic penalty parameter updating strategy.  With this strategy in hand, we present the details of a complete algorithm---an instance of this framework---in \S\ref{sec.complete_algorithm}.

%We formulate our problem of interest as 
Consider the following nonlinear optimization problem with equality and inequality constraints where we assume that the functions $f : \Rn \to \R$ and $c : \Rn \to \Rm$ are continuously differentiable:
\bequation\tag{NLP}\label{prob.nlp}
  \baligned
    \min_{x\in\Rn} &\ f(x) & & \\
    \st &\ c_i(x) =    0\ \ \text{for all}\ \ i \in \{1,\dots,\mbar\}; \\
        &\ c_i(x) \leq 0\ \ \text{for all}\ \ i \in \{\mbar+1,\dots,m\}.
  \ealigned
\eequation
Our penalty-SQP framework uses two functions for use in the algorithm and for characterizing first-order stationary solutions.  First, with a penalty parameter $\rho \in \R_+$, we define the measure of infeasibility and exact penalty function
\bequationn
  v(x) = \sum_{i=1}^{\mbar} \left|c_i(x)\right| + \sum_{i=\mbar+1}^m (c_i(x))_+\ \ \text{and}\ \ \phi(x,\rho) = \rho f(x) + v(x).
\eequationn
Generally speaking, our penalty-SQP framework aims to solve \eqref{prob.nlp} through systematic minimization of $\phi(\cdot,\rho)$ for appropriately chosen values of $\rho \in \R{}_{++}$.  However, if the constraints of \eqref{prob.nlp} are infeasible, then the algorithm is designed to return an infeasibility certificate in the form of a stationary point for the \emph{feasibility problem}
\bequation\label{feas prob}
  \min_{x\in\mathbb{R}^n} \phi(x,0),\ \ \text{where}\ \ \phi(x,0) = v(x).
\eequation
Given $\rho \in \R_+$ and $\eta \in \Rm$,
we define the Fritz John function for \eqref{prob.nlp} by
\bequationn
  F(x,\rho,\eta) = \rho f(x) + \langle \eta , c(x) \rangle.
\eequationn
Note that $\rho \in \R_+$ plays a double role as penalty parameter in~$\phi$ and objective multiplier in $F$.  This makes sense from both theoretical and practical perspectives.  First-order stationarity conditions for \eqref{prob.nlp} can be written in terms of~$\nabla F$, the constraint function $c$, and bounds on the dual variables \cite{CurtJohnRobiWaec14}.

In the $k$th iteration of our penalty-SQP framework, the search direction computation is based on a local model of the penalty function about a primal iterate $x^k \in \Rn$ that 
can make use of a dual iterate $\eta^k \in \Rm$.  We define this model over a convex set $X \subseteq \Rn$ containing $\{0\}$ by
\bequationn
  J(d,\rho;x^k,\eta^k) := l(d, \rho; x^k) + \tfrac{1}{2}\langle d, H(\rho;x^k,\eta^k) d\rangle + \delta(d|X),
\eequationn
where $l$ is a linearized model of the penalty function (ignoring $\rho f(x^k)$) defined by
\bequationn
  l(d,\rho;x^k) = \rho \langle \nabla f(x^k),d \rangle + \sum_{i=1}^{\mbar} \left| c_i(x^k) +\langle \nabla c_i(x^k), d\rangle \right| + \sum_{i=\mbar+1}^m (c_i(x^k) + \langle \nabla c_i(x^k), d\rangle)_+
\eequationn
and $H$ represents an approximation of $\nabla_{xx}^2 F$ with
\bequationn
  H(\rho;x^k,\eta^k) \approx \nabla_{xx}^2 F(\rho;x^k,\eta^k) = \rho \nabla_{xx}^2 f(x^k) + \sum_{i=1}^m \eta_i^k \nabla_{xx}^2 c_i(x^k).
\eequationn 
In particular, the search direction $d^k$ is computed as an approximate minimizer of $J(\cdot,\rho_k;x^k,\eta^k)$ for some $\rho_k \in (0,\rho_{k-1}]$, i.e.,
\bequation\tag{QP}\label{prob.qp}
  d^k \approx \arg\min_{d \in \mathbb{R}^n}\ J(d,\rho_k;x^k,\eta^k)\ \ \text{for some}\ \ \rho_k \in (0,\rho_{k-1}].
\eequation
We introduce the set $X$ to allow for the possibility of employing, e.g., a trust region constraint; e.g., for some $\Delta \in \R{}_+$, one may define $X$ such that $X \subset \{d : \|d\|_2 \leq \Delta\}$.

The value $\rho_k \in (0,\rho_{k-1}]$ is computed \emph{during} the iterative solve of~\eqref{prob.qp}.  Roughly speaking, we aim to adjust this value so that the (inexact) solution $d^k$ to~\eqref{prob.qp} predicts progress toward both feasibility and optimality.  In particular, this occurs if the reduction in a linearized model of the feasibility measure, 
\begin{align}
%  \Delta J(d^k,0;x^k,\eta^k) = J(0,0;x^k,\eta^k) - J(d^k,0;x^k,\eta^k),
  \Delta l(d^k, 0;x^k) &:= l(0, 0;x^k) - l(d^k, 0;x^k), \label{eq.l_model_reduction}\\
   \text{where generally}\ \ \Delta l(d^k, \rho_k;x^k) &:= l(0, \rho_k;x^k) - l(d^k, \rho_k;x^k), \label{eq.l_rho_model_reduction}
\end{align}
and the reduction in the local model of the penalty function,
\bequation\label{eq.J_model_reduction}
  \Delta J(d^k,\rho_k;x^k,\eta^k) := J(0,\rho_k;x^k,\eta^k) - J(d^k,\rho_k;x^k,\eta^k),
\eequation
are sufficiently positive, in which case $d^k$ represents a direction of sufficient descent for both $v$ and $\phi(\cdot,\rho_k)$ from $x^k$.  However, if $x^k$ is (nearly) stationary for $v$ and/or for~$\phi(\cdot,\rho_k)$, then requiring both of these reductions to be positive can force the algorithm to compute a highly accurate solution of \eqref{prob.qp} when one is not entirely needed.  Therefore, the precise conditions that $(d^k,\rho_k)$ must satisfy---introduced in the next section---involve margins that allow one or both of these reductions to be small or even negative for an acceptable step.

%The improvement toward reducing $v$ by $d^k$ from $x^k$ can also be reflected by the model reduction
%\bequationn
%  \Delta l(d^k;x^k):=l(0;x^k)-l(d^k;x^k).
%\eequationn
%By Lemma 4.2 and Lemma 4.3 in \cite{BurkCurtWang14}, the quantity $-\Delta l(d^k;x^k)$ plays the role of surrogate for the directional derivative of $v$ at $x^k$ along $d^k$.

Overall, the $k$th iteration of our penalty-SQP strategy proceeds as in Algorithm~\ref{alg.sqo}.  First, a search direction and penalty parameter pair $(d^k,\rho_k)$ is computed by a subproblem solver such that $d^k$ yields reductions in the local 
models of the penalty function and measure of infeasibility that satisfy our conditions in~\S\ref{sec.penalty}.  Then, a line search is performed with respect to the merit function $\phi(\cdot,\rho_k)$ from $x^k$ along the search direction~$d^k$, yielding a stepsize $\alpha_k \in \R_{++}$.  Finally, the new iterate is set as $x^{k+1} \gets x^k + \alpha_kd^k$ and the algorithm proceeds to the $(k+1)$st iteration.  We discuss choices for the new dual iterate $\eta^{k+1}$ with the complete algorithm in \S\ref{sec.complete_algorithm}.

\begin{algorithm}
  \caption{Penalty-SQP Algorithm (Preliminary)}
  \label{alg.sqo}
  \balgorithmic[1]
    \Require $(\gamma, \theta) \in (0,1)$ and $\rho_{-1} \in (0,\infty)$.
    \State Choose $(x^0,\eta^0) \in \R^n \times \R^m$.
    \For{\textbf{all} $k \in \N{}$}
      \State Solve (approximately) \eqref{prob.qp} to obtain $(d^k, \rho_k) \in \R^n \times (0,\rho_{k-1}]$.
      %\State \qquad or stop if a stationarity certificate is satisfied.
      \label{step.qp}
      \State Let $\alpha^k$ be the largest value in $\{\gamma^0,\gamma^1,\gamma^2,\dots\}$ such that 
      \begin{equation*}
%        \phi(x^k+\alpha_kd^k,\rho_k) - \phi(x^k,\rho_k) \leq - \theta \alpha_k   \Delta J(d^k,\rho_k;x^k,\eta^k) .
         \phi(x^k+\alpha_k d^k, \rho_k) - \phi(x^k, \rho_k) \leq -\theta\alpha_k \Delta l(d^k, \rho_k; x^k).
      \end{equation*}
      \State Set $x^{k+1} \gets x^k + \alpha_kd^k$ and choose $\eta^{k+1} \in \R^m$.
    \EndFor
  \ealgorithmic
\end{algorithm}

Before proceeding, it is worthwhile to emphasize the benefit of ignoring the term $\rho f(x^k)$ in our definitions of the models $J$ and $l$ above.  It is valid to do this since this term has no effect on the solution of \eqref{prob.qp}, and since its presence would not affect the model reduction values in \eqref{eq.l_model_reduction} and \eqref{eq.J_model_reduction}.  On the other hand, ignoring this term simplifies our presentation and analysis significantly since it allows us to avoid the fact that, if this term were not ignored, then the optimal value of \eqref{prob.qp} for a given~$x^k$ would shift with changes in the penalty parameter.

%********* 
% Section
%********* 
\section{A Dynamic Penalty Parameter Updating Strategy}\label{sec.penalty}
 
In this section, we present a dynamic penalty parameter updating strategy.  
As mentioned, the method is novel since the update is employed 
\emph{within} a solver for the subproblem arising in our penalty-SQP framework.  A potential pitfall of such an approach is that, since the penalty parameter dictates the weight between the objective terms in \eqref{prob.qp}, one may disrupt typical convergence guarantees of the subproblem solver by manipulating this weight during the solution process.  However, under reasonable assumptions, we prove that for sufficiently small values of the penalty parameter, our updating strategy will no longer be triggered.  Consequently, once the penalty parameter reaches a sufficiently small value, it will remain fixed and the subproblem solver will effectively be applied to solve \eqref{prob.qp} for a fixed value $\rho_k$.  
%The updating strategy is described in a manner that allows it to be incorporated into various subproblem solvers; see \S\ref{sec.solvers}.

%************
% Subsection
%************
\subsection{Preliminaries}

%As our penalty parameter updating strategy is to be employed in each iteration of our penalty-SQP framework, we can present our strategy generically by focusing on the $k$th iteration of the framework.  Thus, 
For ease of exposition in this section, %we utilize the following shorthand notation to 
we drop the dependence of certain quantities on the iteration number:
\bequation\label{eq.shorthand}
  \baligned
    & g = \nabla f(x^k),\ a^i = \nabla c_i(x^k),\ b_i = c_i(x^k),\ A = [a^1, \cdots, a^m]^T, \\
    & H_f \approx \nabla_{xx}^2 f(x^k),\ H_0 \approx \sum_{i=1}^m \eta_i^k \nabla_{xx}^2 c_i(x^k),\ \text{and}\ H_\rho = \rho H_f + H_0.
  \ealigned
\eequation
We also temporarily drop the dependence of the functions $J$, $l$, etc.~on the $k$th iterate.

We make the following assumption about the subproblem data.

\bassumption\label{ass.qp}
  The subproblem data matrices $A$, $H_f$, and $H_0$ are such that
  \benumerate
    \item[(i)]  $H_\rho$ is positive definite for any $\rho \in [0,\rho_{k-1}]$; and
    \item[(ii)] $\|a^i\|_2 > 0$ for all $i \in \{1,\dots,m\}$.
  \eenumerate
\eassumption

\noindent
We claim that this assumption is reasonable due to the following considerations.  First, in large-scale contexts, it is typically impractical to construct complete second-derivative matrices.  Hence, as indicated in \eqref{eq.shorthand}, one can assume that $H_f$ and $H_0$ represent (limited memory) Hessian approximations with at least~$H_0$ being positive definite.  %(See \S\ref{sec.low_rank} for further discussion.)  
Second, if $a^i = 0$ for any $i \in \{1,\dots,m\}$, then the model of the $i$th constraint is constant with respect to~$d$, meaning that the $i$th constraint can be removed from the subproblem.  Such a phenomenon can be detected during a preprocessing phase before solving the subproblem, so for simplicity, we assume that each constraint gradient is nonzero.  Under Assumption~\ref{ass.qp}, we define the scaled quantities $\abar^i := a^i/\|a^i\|_2$ and $\bbar_i := b_i/\|a^i\|_2$ for all $i \in \{1,\dots,m\}$.

Of central importance in the subproblems are the convex sets
\bequalin
  C_i &:= \{d \in \Rn : \langle \bar a^i,d \rangle + \bar b_i = 0\}\ \ \text{for all}\ \ i \in \{1,\dots,\mbar\} \\
  \text{and}\ \ C_i &:= \{d \in \Rn : \langle\bar a^i,d\rangle + \bar b_i \leq 0\}\ \ \text{for all}\ \ i \in \{\mbar+1,\dots,m\}.
\eequalin
%At a feasible iterate $x^k$, the penalty term in the model $J$ can be written as $\sum\limits_{i=1}^m \delta(d| C_i)$. 
%, and it is suggested that $\rho$ should be fixed while solving the subproblem.
% We focus on designing the updating strategy for $\rho$ when \emph{$x^k$ is neither feasible nor infeasible stationary}. However, we claim in Lemma~\ref{lem.vstat} that our proposed strategy is in fact not triggered for a feasible $x^k$. 
The quadratic and penalty terms in $J$ can %thus 
be written, respectively, as
\bequationn
  \psi(d,\rho) = \rho \langle g, d\rangle + \thalf \langle d, H_\rho d\rangle \ \ \text{and}\ \ l(d,0)=\sum_{i=1}^m \|a^i\|_2 \tdist{d}{C_i},
\eequationn
meaning that %, without loss of generality (i.e., assuming $\|a_i\|_2=1$ for all $i \in \{1,\dots,m\}$) 
we may rewrite the penalty-SQP subproblem \eqref{prob.qp} as
\bequation\tag{QPrho}\label{prob.J}
  \min_{d\in\Rn}\ J(d,\rho),\ \ \text{where}\ \ J(d,\rho) = \psi(d,\rho)+ l(d,0) + \delta(d| X).
\eequation
We refer to \eqref{prob.J} with $\rho>0$ as a \emph{penalty subproblem} and we refer to \eqref{prob.J} with $\rho=0$ as the \emph{feasibility subproblem}.  The Fenchel--Rockafellar dual of~\eqref{prob.J} is
\bequation\tag{DQPrho}\label{prob.D}
  \baligned
    \max_{\umbf\in\Rn\times\cdots\times\Rn}\ D(\umbf,\rho)\ \ \st &\ u^0+\sum_{i=1}^{m}\|a^i\|_2 u^i + u^{m+1}= 0 \\ \text{and} &\ u^i \in \bB_2\ \ \text{for all}\ \ i \in \{1,\dots,m\},
  \ealigned
\eequation
where the dual objective function is given by
\bequationn
  D(\umbf,\rho) = -\thalf \langle u^0 - \rho g,  H_\rho^{-1} (u^0 - \rho g)\rangle - \sum_{i=1}^m\|a^i\|_2 \delta^*(u^i | C_i) - \delta^*(u^{m+1} | X).
\eequationn
Letting $\zeta_i(\umbf) := \langle u^i, \bar a^i\rangle$ for a dual feasible $\umbf$, one finds from \eqref{delta*} and the constraint in \eqref{prob.D} that 
$D(\umbf, \rho)$ is finite if and only if
\bequation\label{eq.zeta}
  \baligned
    &\ u^i = \zeta_i(\umbf) \bar a^i, \\
    \text{which means}\ \ &\ \zeta_i(\umbf) \in \begin{cases} [-1,1] & \text{for all $i \in \{1,...,\bar m\}$} \\ [0,1] & \text{for all $i \in \{\bar m+1, ..., m\}$}, \end{cases} \\
    \text{and}\ \ &\ \delta^*(u^i|C_i) = - \zeta_i(\umbf) \bar b_i.
  \ealigned
\eequation

An interesting aspect of the dual subproblem \eqref{prob.D} is that the penalty parameter appears only in the objective.  Thus, if $\umbf$ satisfies the constraints of \eqref{prob.D}, then it is dual-feasible regardless of the value of $\rho$ appearing in the subproblem.  As a result, by weak duality, we have for any primal-dual feasible pair $(d,\umbf)$ that both
\bequation\label{weak.duality}
  D(\umbf,0) \leq J(d,0)\ \ \text{and}\ \ D(\umbf,\rho) \leq J(d,\rho).
\eequation

We close this subsection by noting that  
 the projection onto the set $C_i$
 \bequationn
  P_{C_i}(y^i) := \arg\min_{z^i \in C_i} \tnorm{z^i -y^i}%\ \ \text{and}\ \ P_\Cmbf(\ymbf) := \arg\min_{\pmb{z} \in \Cmbf} \|\pmb{z} - \ymbf\|.
\eequationn
 is easy to compute for any $i \in \{1,\dots,m\}$; in particular,
\bequalin
  P_{C_i}(d)  = &\begin{cases} d-(\langle \abar^i,d \rangle +\bbar_i)   \abar^i & \text{for all}\ \ i \in \{1,\dots,\mbar\} \\
  d-(\langle \abar^i,d \rangle +\bbar_i)_+ \abar^i & \text{for all}\ \ i \in \{\mbar+1,\dots,m\}.
  \end{cases}
\eequalin
%We can see that if $X = \Rn$, then the Fenchel-Rockafellar dual of \eqref{prob.J} reduces to
%\bequationn
%  \baligned
%    \max_{\umbf\in\Rn\times\cdots\times\Rn} &\ -\thalf (u^0 - \rho g)^T H_\rho^{-1} (u^0 - \rho g) - \sum_{i=1}^m\|a^i\|_2 \delta^*(u^i | C_i) \\ 
%    \st &\ u^0+\sum_{i=1}^{m}\|a^i\|_2 u^i + u^{m+1}= 0 \\ 
%    \text{and} &\ u^i \in \bB_2\ \ \text{for all}\ \ i \in \{1,\dots,m\}.
%  \ealigned
%\eequationn

%************
% Subsection
%************
\subsection{Updating the penalty parameter}\label{sec.update}

%We now describe our dynamic penalty parameter updating strategy.  
Given $\rho \ge 0$, 
let $(d_{\rho}^*,\umbf_{\rho}^*)$ represent an optimal primal-dual pair for the penalty subproblem \eqref{prob.J} corresponding to~$\rho$; in particular, $(d_0^*,\umbf_0^*)$ represents an optimal primal-dual pair for the feasibility subproblem.  
%We describe 
The algorithm is presented 
in the context of a subproblem solver that generates two sequences of iterates:  
the first sequence, call it $\{(d^{(j)},\umbf^{(j)})\}$, is a sequence of primal-dual feasible solution estimates for a penalty subproblem, while the second sequence, 
call it $\{\wmbf^{(j)}\}$, is a sequence of dual feasible solution estimates for the feasibility subproblem.  (In our strategy, we do not make separate use of a sequence of primal solution estimates for the feasibility subproblem; rather, the sequence $\{d^{(j)}\}$ plays this role as well.)  Without loss of generality, we assume that the $j$th primal solution estimate $d^{(j)}$ represents a better (or no worse) primal solution estimate for the penalty subproblem than a zero step in the sense that
\bequation\label{eq.primal_no_worse_than_zero}
  J(d^{(j)},\rho_{(j)}) \leq J(0,\rho_{(j)}).
\eequation
Similarly, we assume that the dual solution estimate $\wmbf^{(j)}$ represents a better (or no worse) dual solution estimate for the feasibility subproblem than $\umbf^{(j)}$, and that each dual solution estimate $\umbf^{(j)}$ is no worse than the feasible $\umbf^{(0)}$, in that
\bequation\label{eq.dual_no_worse}
  D(\wmbf^{(j)},0) \geq D(\umbf^{(j)},0) \geq D(\umbf^{(0)},0) > - \infty.
\eequation
These are both reasonable assumptions since if \eqref{eq.primal_no_worse_than_zero} (resp.~\eqref{eq.dual_no_worse}) were not to hold, then one could consider $d^{(j)} = 0$ (resp.~$\wmbf^{(j)} = \umbf^{(j)} = \umbf^{(0)}$) for the $j$th iterate (even if the subproblem solver works with a different estimate in its internal operations). %{\red In practice, we know $D(\wmbf^{(j)},0) \ge 0$, so that we can further assume the dual solution estimate achieves dual value no less than 0.}

Observe that, by the definition of the model $J$, we have for any $\rho \in (0,\infty)$ that
\bequationn
  J^{(0)} := J(0,\rho) = J(0,0) = l(0, 0) = \sum_{i=1}^{\bar m} |b_i|+\sum_{i=\bar m+1}^m (b_i)_+ \geq 0.
\eequationn
Let $J^{(0)}_\omega := J^{(0)} + \omega$ for any scalar $\omega\in(0,\infty)$.  (As discussed later, $\omega$ is held fixed during a given subproblem solve, but will sequentially be reduced to zero over the course of the overall penalty-SQP framework.)  We then define the following ratios corresponding to the $j$th subproblem solver iterate:
\bequation
  r_v^{(j)} := \frac{J^{(0)}_\omega - l(d^{(j)}, 0) }{J^{(0)}_\omega - (D(\wmbf^{(j)},0))_+}\ \ \text{and}\ \ r_\phi^{(j)} := \frac{J^{(0)}_\omega - J(d^{(j)},\rho_{(j)})}{J^{(0)}_\omega - D(\umbf^{(j)},\rho_{(j)})}. 
\eequation
(Referring back to our discussion surrounding \eqref{eq.l_model_reduction} and \eqref{eq.J_model_reduction}, note that the numerators of these ratios are $\Delta l(d^{(j)}, 0) + \omega$ and $\Delta J(d^{(j)},\rho_{(j)}) + \omega$, respectively.)  The critical property of these ratios is that, if they are sufficiently large, then the corresponding subproblem solver iterates must yield reductions in the feasibility and penalty function models that are proportional to those obtained by corresponding exact subproblem solutions.  In particular, suppose that for some prescribed $\beta_v \in (0,1)$ we have
\begin{equation}\tag{Rv}\label{red.fea}
   r_v^{(j)} \ge \beta_v.
\end{equation}
Then the reduction in the linearized constraint violation model obtained by the subproblem solver iterate $d^{(j)}$ relative to a zero step satisfies
\bequation\label{alter.red.constr}
  \begin{aligned}
J_\omega^{(0)} - l(d^{(j)},0) & \geq  \beta_v\(J^{(0)}_\omega - ( D(\wmbf^{(j)},0))_+\)\\
    &\geq \beta_v\(J^{(0)}_\omega - D(\umbf_0^*,0)\) = \beta_v\(J^{(0)}_\omega - J(d_0^*,0)\),
  \end{aligned}
\eequation
where the first inequality follows by \eqref{red.fea}, the second follows by the optimality of~$\umbf_0^*$ with respect to the feasibility subproblem (for which it is known that $D(\umbf_0^*,0) \geq 0$), and the last follows by strong duality.  Similarly, if for $\beta_\phi \in (0,1)$ we have
\begin{equation}\tag{Rphi}\label{red.penalty}
  r_\phi^{(j)} \geq \beta_\phi,
\end{equation}
then it follows that
\bequation\label{alter.red.penalty}
  \begin{aligned}
    J^{(0)}_\omega - J(d^{(j)},\rho_{(j)})
      &\geq \beta_\phi(J^{(0)}_\omega - D(\umbf^{(j)},\rho_{(j)})) \\
      &\geq \beta_\phi(J^{(0)}_\omega - D(\umbf_{\rho_{(j)}}^*,\rho_{(j)})) = \beta_\phi(J^{(0)}_\omega - J(d_{\rho_{(j)}}^*,\rho_{(j)})).
  \end{aligned}
\eequation

The last component of our updating strategy involves an estimate of the complementarity of a primal-dual solution estimate.  This is needed since we only reduce the penalty parameter if a primal-dual solution estimate is approximately complementary.  We do this in the following manner.  First, defining the index sets
\[ \begin{aligned}
 \Ecal_+(d) & := \{ i \in \{1,\dots,\bar m\} : \langle \bar a^i, d\rangle + \bar b_i > 0 \},\\
 \Ecal_- (d) & := \{ i\in\{1,\dots,\bar m\} : \langle \bar a^i, d\rangle + \bar b_i < 0 \},\\
 \text{and}\ \ \Ical_+(d)  & := \{ i\in\{\bar m+1,\dots, m\} : \langle \bar a^i, d\rangle + \bar b_i > 0\},
 \end{aligned}
\]
we define the complementarity measure
\[\chi(d,\umbf):=
 \sum_{i\in\Ecal_+\cup\Ical_+} (1- \zeta_i(\umbf)) \|a^i\|_2  \dist{d}{C_i}  + \sum_{i\in\Ecal_-}( 1+\zeta_i(\umbf) ) \|a^i\|_2 \dist{d}{C_i}. 
 \]
To reduce the penalty parameter, we require that $(d^{(j)}, \umbf^{(j)})$ satisfies
\[ \chi^{(j)} := \chi(d^{(j)},\umbf^{(j)})  \le (1 - \beta_v)^2J^{(0)}_\omega,\]
or, equivalently, 
\begin{equation}\tag{Rc}\label{red.comp}
r_c^{(j)} := 1-  \sqrt{\frac{ \chi^{(j)}}{J^{(0)}_\omega}} \ge \beta_v.
\end{equation}

In our strategy, if the optimality QP subproblem is solved sufficiently accurately, then we turn to verify whether feasibility has also been improved to a satisfactory extent. Therefore, the key idea here is to determine a criterion reflecting that the optimality QP has been solved sufficiently accurately.  Making this determination requires us to check a measure of complementarity.  In particular, if the initial objective $J^{(0)}$ is far from optimal, then $r_\phi^{(j)} \approx 1$ might not indicate that the subproblem solution 
is nearly primal-dual optimal since 
a large $J^{(0)}$ can cause  
 the numerator  of $r_\phi^{(j)}$ to be very close to the denominator, even though 
 the dual value is far from dual optimality. As a result, the updating strategy 
 may be triggered too early, so that  $\rho$ is inappropriately driven to zero.  Therefore, we need a 
 certification showing the progress achieved by the dual estimates, which can be reflected 
 by the complementary condition \eqref{red.comp}.

Overall, our penalty parameter strategy is motivated by the desire to ensure that if the $j$th iterate of the subproblem solver offers a sufficiently accurate solution of the penalty subproblem for $\rho_{(j)} > 0$, then it should also offer a sufficiently accurate solution of the feasibility subproblem; otherwise, the penalty parameter should be reduced.  Specifically, choosing parameters
\bequation\label{eq.parameters}
  0 < \beta_v < \beta_\phi < 1,
\eequation
we initialize $\rho_{(0)} \gets \rho_{k-1}$ 
(from the preceding iteration of the penalty-SQP framework) and apply the subproblem solver to \eqref{prob.J} to initialize $\{(d^{(j)},\umbf^{(j)},\wmbf^{(j)})\}$.  If, at the end of the $j$th subproblem solver iteration we have that \eqref{red.penalty} or \eqref{red.comp} is not satisfied, then we continue to iterate toward solving \eqref{prob.J} with $\rho = \rho_{(j)}$.  Otherwise, if \eqref{red.penalty} and \eqref{red.comp} hold but \eqref{red.fea} does not, then we reduce the penalty parameter by setting
\begin{equation}\label{update.rho}
  \rho_{(j+1)} \gets \theta_\rho \rho_{(j)}
\end{equation}
for some prescribed $\theta_\rho\in (0,1)$.  (A special case that one should consider occurs when \eqref{red.penalty}, \eqref{red.comp}, and \eqref{red.fea} all hold with $d^{(j)} = 0$.  For simplicity in our presentation, in such a case, we have the subproblem solver terminate with $d^{(j)} = 0$, causing the penalty-SQP framework to take a null step in the primal space.  As previously mentioned, this would be followed by a decrease in $\omega$, prompting the penalty-SQP framework to eventually make further progress or terminate with a stationarity certificate.  In practice, this decrease in $\omega$ in this scenario need not occur over a sequence of iterations.  It can occur immediately within a subproblem solve.  We merely state the occurrence of a null step for simplicity in our discussions.)

We state our \emph{d}ynamic \emph{u}pdating \emph{st}rategy (DUST) as:
\begin{equation}\tag{DUST}\label{dust}
  \boxed{ 
  \baligned
   & \text{\small Given $\rho_{(j)}$ and the $j$th iterate $(d^{(j)},\umbf^{(j)},\wmbf^{(j)})$, perform the following:} \\
   & \text{\small \qquad $\bullet$ if %$d^{(j)} = 0$ and 
   \eqref{red.penalty}, \eqref{red.comp}, and \eqref{red.fea} hold, then terminate;} \\
   & \text{\small \qquad $\bullet$ else if \eqref{red.penalty} and \eqref{red.comp} hold, but \eqref{red.fea} does not, then apply \eqref{update.rho};} \\
   & \text{\small \qquad $\bullet$ else set $\rho_{(j+1)} \gets \rho_{(j)}$.}
  \ealigned
  }
\end{equation}

We formally analyze \eqref{dust} in the following subsections.  
%For now, to preview our analysis, let us present 
We begin with the following intuitive arguments to motivate
the strategy for adjusting the penalty parameter in a few cases of interest.  
These cases depend on properties of the $k$th iterate of the penalty-SQP framework, namely, $x^k$, with respect to the constraint violation measure and the penalty function.

\bitemize
  \item First, observe that with an optimal primal-dual solution $(d_\rho^*, \umbf_\rho^*)$ for a penalty subproblem, one has $\zeta_i(\umbf_\rho^*) = 1$ for $i \in \Ecal_+(d_\rho^*)$, $\zeta_i(\umbf_\rho^*) = -1$ for $i \in \Ecal_-(d_\rho^*)$, and $\zeta_i(\umbf_\rho^*) = 1$ for $i \in \Ical_+(d_\rho^*)$, from which it follows that $\chi(d_\rho^*,\umbf_\rho^*)=0$.  Therefore, for a given $\omega \in (0,\infty)$, the condition \eqref{red.comp} will hold for sufficiently accurate primal-dual solutions of the penalty subproblem. 
  \item If $x^k$ is not stationary with respect to $\phi(\cdot,\rho)$ for any $\rho \in (0,\rho_{k-1}]$, then, with $(d^{(j)},\umbf^{(j)},\rho_{(j)}) = (d_{\rho}^*,\umbf_{\rho}^*,\rho)$ for any such $\rho$, one finds that $r_\phi^{(j)} = 1 > \beta_\phi$.  In turn, this means that \eqref{red.penalty} holds for any $(d^{(j)},\umbf^{(j)})$ in a neighborhood of $(d_{\rho}^*,\umbf_{\rho}^*)$.  If, in addition, $x^k$ is not stationary with respect to $v$, then one should expect that for a sufficiently small $\rho_{(j)}$ the condition \eqref{red.fea} would also be satisfied for such a~$d^{(j)}$.  This should be expected since for $(d_0^*,\umbf_0^*)$ one has
  \bequationn
    \frac{J^{(0)}_\omega - l(d_0^*, 0)}{J^{(0)}_\omega - (D(\umbf_0^*,0))_+} \geq \frac{J^{(0)}_\omega - J(d_0^*,0)}{J^{(0)}_\omega - D(\umbf_0^*,0)} = 1,
  \eequationn
  meaning that $r_v^{(j)} > \beta_v$ for $(d^{(j)},\wmbf^{(j)})$ in a neighborhood of $(d_0^*,\umbf_0^*)$.  Overall, in this case, one should expect that \eqref{dust} would only reduce the penalty parameter a finite number of times, if at all.
  \item If $x^k$ is not stationary with respect to $\phi(\cdot,\rho)$ for any $\rho \in (0,\rho_{k-1}]$, but is stationary with respect to $v$, then for $(d_0^*,\umbf_0^*)$ one has
  \bequationn
    \frac{J^{(0)}_\omega - l(d_0^*,0)}{J^{(0)}_\omega - (D(\umbf_0^*,0))_+} = \frac{\omega}{\omega} = 1,
  \eequationn
  meaning that $r_v^{(j)} > \beta_v$ for $(d^{(j)},\wmbf^{(j)})$ in a neighborhood of $(d_0^*,\umbf_0^*)$.  Hence, as in the previous bullet, one should expect that \eqref{dust} would only reduce the penalty parameter a finite number of times.
  \item If $x^k$ is stationary with respect to $\phi(\cdot,\rho_{(j)})$ for $\rho_{(j)} > 0$ encountered during the subproblem solve, then, under Assumption~\ref{ass.qp}, the only primal iterate satisfying \eqref{red.penalty} is $d^{(j)} = 0$.  For this value, one finds that
  \bequationn
    r_v^{(j)} = \frac{\omega}{\omega + J^{(0)} - (D(\wmbf^{(j)},0))_+}.
  \eequationn
  There are now two cases to consider.  If $r_v^{(j)} < \beta_v$, then \eqref{dust} decreases the penalty parameter, as is appropriate.  Otherwise, if $r_v^{(j)} \geq \beta_v$, then---with a sufficiently accurate dual solution---\eqref{dust} returns a null step to the penalty-SQP framework.  (In a later subproblem solve with a smaller $\omega$, one would either find that \eqref{red.penalty} holds for $d^{(j)} = 0$---and a sufficiently accurate dual solution---but \eqref{red.fea} does not, prompting a decrease of the penalty parameter, or---again with a sufficiently accurate dual solution---one would terminate the overall algorithm with certificate of stationarity for $x^k$.)
\eitemize
 
We close this subsection by making a few practical remarks regarding the use of \eqref{dust} within a subproblem solver for \eqref{prob.J}.  In particular, while we have defined the sequence $\{(d^{(j)},\umbf^{(j)},\wmbf^{(j)})\}$ as being generated by the solver, it may be reasonable to reinitialize the solver---or at least perform some auxiliary computations---after any iteration in which \eqref{update.rho} is invoked.  (Such auxiliary computations may involve scaling vectors and/or matrices due to the change in the penalty parameter.)  %; e.g., see the discussion of the Hessian approximation strategy in \S\ref{sec.low_rank}.)  
That being said, it is reasonable to assume that, during any sequence of iterations in which the penalty parameter does not change, the subproblem solver would be applied as if it were being applied to a static instance of \eqref{prob.J}.  In such a manner, any convergence guarantees for the subproblem solver would hold if/when the penalty parameter stabilizes at a fixed value, as is guaranteed to occur under common conditions described next.

%************

% Subsection
%************
\subsection{Finite Updates for a Single Subproblem}\label{sec.bound.rhoj}

The purpose of this subsection is to show that if \eqref{dust} is employed within an algorithm for solving \eqref{prob.J}, then, under reasonable assumptions on the subproblem data, for any $\rho_{(j)} \in (0,\tilde\rho]$ for some sufficiently small $\tilde\rho > 0$ whose value depends only on the subproblem data,  if \eqref{red.penalty} and \eqref{red.comp} are satisfied, then \eqref{red.fea} is also satisfied.  In other words, after a finite number of iterations, the update \eqref{update.rho} will never be triggered.  Let $\underline{\lambda}_0$ and $ \overline{\lambda}_0$ be the smallest and largest eigenvalues of $H_0$, and similarly for $\underline{\lambda}_\rho$ and $\overline\lambda_\rho$ with respect to the matrix $H_\rho$.  Notice that, since $\rho_{(j)} \in (0, \rho_{(0)}]$, it follows that
\beq\label{eig.bd} 
\underline{\lambda}_{\rho_{(j)} } \ge \underline\lambda:=  \min\{\underline{\lambda}_{\rho_{(0)}}, \underline\lambda_0\}\quad \text{and}\quad 
 \overline{\lambda}_{\rho_{(j)} }  \le \overline\lambda: = \max\{ \overline{\lambda}_{\rho_{(0)}}, \overline\lambda_0\}.
 \eeq
%  Let $\{\rho_{(j)}\}$ be a sequence generated along with $\{(d^{(j)},\umbf^{(j)},\wmbf^{(j)})\}$ by applying \eqref{dust}. 
  % Then, defining
 % \bequationn
  %    \Ucal = \{j \in \N{} : (d^{(j)},\umbf^{(j)})\ \text{satisfies}\ \eqref{red.penalty}\},
 % \eequationn
We formalize our assumption for this analysis as the following.

%\begin{assu}\label{ass.algorithm}
%  Let $\{(d^{(j)},\umbf^{(j)},\wmbf^{(j)})\}$ be a sequence such that, for all $j \in \N{}$, $d^{(j)} \in X$ while $\umbf^{(j)}$ and $\wmbf^{(j)}$ are feasible for \eqref{prob.D}.  Moreover, let $\{\rho_{(j)}\}$ be a sequence generated along with $\{(d^{(j)},\umbf^{(j)},\wmbf^{(j)})\}$ by applying \eqref{dust}.  Then, defining
%  \bequationn
%      \Ucal = \{j \in \N{} : (d^{(j)},\umbf^{(j)})\ \text{satisfies}\ \eqref{red.penalty}\},
%  \eequationn
%  the subsequences $\{\|\umbf^{(j)}\|_2\}_{k\in\Ucal}$ and $\{\|\wmbf^{(j)}\|_2\}_{k\in\Ucal}$ are bounded by a constant $\kappa_0 > 0$ independent of $\{\rho_{(j)}\}$.
%\end{assu}

%This boundedness assumption on the dual estimates is reasonable since our subproblems are assumed to be strongly convex. \footnote{\red I wonder whether this Assumption is not needed here if we assume $w^{(j)}$ and $u^{(j)}$ are not getting ``worse and worse". (In fact, we haven't assumed $u^{(j)}$ is getting better and better while solving the subproblem)  This is easy to see, since 
%in \eqref{prob.D}, $u^i, i=1,...,m$ are bounded, $\delta^*(u^{m+1}|X) \ge 0$ for any $u^{(m+1)}$ since $0\in X$. Then if $u^{m+1}\to \infty$, we know $u^0\to\infty$ and consequently 
%$D(\umbf, \rho) \to - \infty$.  Thus, we actually have that all dual iterates are bounded, not just the subsequence $\Ucal$.}

\begin{assu}\label{ass.algorithm}
  For all $j \in \N{}$, the sequence $\{(d^{(j)},\umbf^{(j)},\wmbf^{(j)})\}$ has $d^{(j)} \in X$, \eqref{eq.primal_no_worse_than_zero} and \eqref{eq.dual_no_worse} hold, and $\umbf^{(j)}$ and $\wmbf^{(j)}$ are feasible for \eqref{prob.D}.
\end{assu}

We first show that the dual sequences $\{\umbf^{(j)}\}$ and $\{\wmbf^{(j)}\}$ are bounded in norm. 
\begin{lemma}\label{lem.bounds.dual} 
  Under Assumption~\ref{ass.qp}, there exists $\kappa_0>0$ such that, for all $j \in \N{}$,
\[\|\umbf^{(j)}\|_2 \le \kappa_0\quad \text{and}\quad \|\wmbf^{(j)}\|_2 \le \kappa_0.\]
\end{lemma} 

\begin{proof}

Since $\umbf^{(j)}$ is feasible for \eqref{prob.D},  the elements $ \{ (u^i)^{(j)} \}$ for all $i\in \{1,\dots,m\}$ are bounded in norm by 1.  Therefore, by the first constraint of \eqref{prob.D}, it suffices to show that $\{ (u^0)^{(j)} \}$ is bounded.  We show this by contradiction. Suppose there exists an infinite index set $\Jcal$ such that $\{ \| (u^0)^{(j)} \|_2\}_{j\in\Jcal} \nearrow \infty$.  Notice that for $(u^{m+1})^{(j)}$ it holds  that  $\delta^*( (u^{m+1})^{(j)} |X) = \sup\limits_{x\in X} \langle (u^{m+1})^{(j)}, x\rangle \ge 0$  since it is assumed that $0\in X$.  All together, with these facts and Assumption~\ref{ass.qp}, we may conclude that $\{ D(\umbf^{(j)}, 0) \}_{j\in\Jcal} \to -\infty$, which contradicts \eqref{eq.dual_no_worse}. Therefore, $\{ (u^0)^{(j)} \}$ must be bounded, so overall the sequence $\{ \umbf^{(j)} \}$ is bounded.

Following the same argument for $\wmbf^{(j)}$, it follows that $\{\wmbf^{(j)}\}$ is bounded. 
\end{proof}

We now show that the primal variables $\{d^{(j)}\}$ are also bounded in norm.

\begin{lemma}\label{lem.bounds}
  Under Assumptions~\ref{ass.qp} and \ref{ass.algorithm}, it follows that, for all $j \in \N{}$,  
  \beq\label{d.bd}
    \|d^{(j)} \|_2 \le \kappa_1 := \(\rho_{(0)}\|g\|_2 + \sqrt{\rho_{(0)}^2\|g\|^2_2 + 2 \overline\lambda J^{(0)}}\)/\underline\lambda.
  \eeq
\end{lemma}
\begin{proof}
  By Assumption~\ref{ass.algorithm}, it follows that $d^{(j)} \in X$ for all $j \in \N{}$, which implies that $\delta(d^{(j)}| X)=0$ for all $j\in \N{}$.  By  \eqref{eq.primal_no_worse_than_zero}, every $(d^{(j)}, \umbf^{(j)}, \rho_{(j)})$ for $j\in\N{}$ must satisfy
  \bequationn
    \rho_{(j)} \langle g, d^{(j)}\rangle + \thalf \langle d^{(j)}, H_{\rho_{(j)}} d^{(j)} \rangle \le  J(d^{(j)}, \rho_{(j)}) \leq J(0,\rho_{(j)}) = J^{(0)}.
  \eequationn
  It follows that 
  \bequationn
    \thalf \underline\lambda_{\rho_{(j)}} \|d^{(j)}\|_2^2 \le  J^{(0)} + |\rho_{(j)} \langle g, d^{(j)} \rangle| \le J^{(0)} + \rho_{(0)} \| g\|_2 \|d^{(j)}\|_2,
  \eequationn
  which, using the quadratic formula, implies that
  \bequationn 
    \|d^{(j)}\|_2   \le \(\rho_{(0)}\|g\|_2 + \sqrt{\rho_{(0)}^2\|g\|^2_2 + 2 \underline\lambda_{\rho_{(j)}} J^{(0)}  } \)  /\underline\lambda_{\rho_{(j)}}.
  \eequationn
  Together with \eqref{eig.bd}, this proves \eqref{d.bd}, as desired.
\end{proof}

The next lemma shows that the differences between the primal and dual values of the penalty and feasibility subproblems are bounded with respect to $\rho$. 

\begin{lemma}\label{lem primal dual value}
  Under Assumptions~\ref{ass.qp} and \ref{ass.algorithm}, it follows that, for any $j\in \N{}$,
  \begin{subequations}
    \begin{align}
      |J (d^{(j)}, \rho_{(j)}) - J (d^{(j)},0) | & \le \kappa_2 \rho_{(j)} \label{J.krho}\\ \text{and}\ \ 
      |D(\umbf^{(j)}, \rho_{(j)}) - D(\umbf^{(j)},0) | & \le \kappa_3 \rho_{(j)},\label{D.krho}
    \end{align}
  \end{subequations}
  where, with $\kappa_1 > 0$ defined in Lemma~\ref{lem.bounds},
  \begin{align*}
  \kappa_2 & := \|g\|_2\kappa_1 + \tfrac{1}{2}\|H_f\|_2 \kappa_1^2 \\
  \text{and}\ \ \kappa_3 & :=\frac{ \kappa_0+\rho_{(0)}\|g\|_2}{2\underline\lambda} (\kappa_0 \|H_0 ^{-1}\|_2\|H_f \|_2 +  \|g\|_2) + \thalf \kappa_0 \|H_0 ^{-1}\|_2 \|g\|_2. 
  \end{align*}
\end{lemma}
\begin{proof}
  For the primal values,  it holds true that   
  \begin{equation*}
    \begin{aligned}
      |J(d^{(j)},\rho_{(j)}) - J(d^{(j)}, 0)| & = |\rho_{(j)} \langle g, d^{(j)}\rangle + \thalf \langle d^{(j)}, H_{\rho_{(j)}} d^{(j)}\rangle - \thalf \langle d^{(j)},  H_0d^{(j)}\rangle |\\
      & = |\rho_{(j)} \langle g, d^{(j)}\rangle  + \tfrac{1}{2} \rho_{(j)} \langle d^{(j)},  H_f d^{(j)}\rangle |  \\
      &   \le \rho_{(j)} (\|g\|_2\| d^{(j)} \|_2 + \tfrac{1}{2} \| H_f \|_2 \|d^{(j)} \|_2^2),
    \end{aligned}
  \end{equation*}
  which combined with Lemma~\ref{lem.bounds}   proves \eqref{J.krho}.
   
   We now aim to prove \eqref{D.krho}.  Toward this goal, let $\hat y^{(j)} := H_{\rho_{(j)}}^{-1}  (u_0^{(j)} - \rho_{(j)} g)$ and $\bar y^{(j)} := H_0 ^{-1}  u_0^{(j)}$.  Then, by Assumption~\ref{ass.algorithm}, it follows that
  \[ \|\hat y^{(j)}\|_2 \le  (\kappa_0+\rho_{(j)}\|g\|_2)/\underline{\lambda}_{\rho_{(j)}} \le 
   (\kappa_0+\rho_{(0)}\|g\|_2)/\underline\lambda.
  \]
  In addition, it follows that
  \[
    \rho_{(j)} g = u_0^{(j)} - (u_0^{(j)} - \rho_{(j)} g)= H_0  \bar y^{(j)} - H_{\rho_{(j)}} \hat y^{(j)} = H_0 (\bar y^{(j)} - \hat y^{(j)} )-  \rho_{(j)} H_f \hat y^{(j)},
  \]
  which implies that, for all $j\in \N{}$,
\begin{equation}\label{auxiliary y}
 \begin{aligned}
 \| \bar y^{(j)} - \hat y^{(j)} \|_2 
     &   = \|\rho_{(j)} H_0 ^{-1}   (H_f \hat y^{(j)} + g) \|_2   \\
    &  \le \rho_{(j)} \|H_0 ^{-1}  \|_2 \|H_f \hat y^{(j)} + g\|_2  \\
    &   \le \rho_{(j)} \|H_0 ^{-1}\|_2 \left(\|H_f\|_2\frac{ \kappa_0+\rho_{(0)}\|g\|_2}{\underline\lambda} + \norm{g}_2\right). 
    \end{aligned}    
  \end{equation}
  The difference between the dual values   is then given by 
  \begin{equation*} 
    \begin{aligned}
      &| D (\umbf^{(j)},\rho_{(j)}) - D (\umbf^{(j)},0) | \\
       =\ & | - \thalf \langle u_0^{(j)}- \rho_{(j)} g, H_{\rho_{(j)}}^{-1}(u_0^{(j)}-\rho_{(j)} g)\rangle  + \thalf \langle u_0^{(j)}, H_0 ^{-1} u_0^{(j)}\rangle| \\
       =\ & | \thalf \langle \bar  y^{(j)} - \hat y^{(j)}, u_0^{(j)} \rangle + \thalf \rho_{(j)} \langle g, \hat y^{(j)}\rangle | \\
       \le\  & \thalf  \|\bar  y^{(j)} - \hat y^{(j)}\|_2 \|u_0^{(j)} \|_2 + \thalf \rho_{(j)}\|g\|_2\| \hat y^{(j)} \|_2\\
       \le\ &  \rho_{(j)} \left( \thalf \|H_0 ^{-1} \|_2 \left(\|H_f \|_2\frac{ \kappa_0+\rho_{(0)}\|g\|_2}{\underline\lambda} + \norm{g}_2\right)\kappa_0
       + \thalf \|g\|_2
   \frac{ \kappa_0+\rho_{(0)}\|g\|_2}{\underline\lambda}    \right) \\
   =\ &  \rho_{(j)} \left( 
    \frac{ \kappa_0+\rho_{(0)}\|g\|_2}{2\underline\lambda} \left( \kappa_0\|H_0 ^{-1} \|_2 \|H_f \|_2 + \|g\|_2\right) + \thalf \kappa_0\|H_0 ^{-1}\|_2  \|g\|_2
   \right),
    \end{aligned}
  \end{equation*}
  where the last inequality follows by \eqref{auxiliary y} and  Assumption~\ref{ass.algorithm}.
\end{proof}

Let us now define 
\bequationn
  \Ucal = \{j : (d^{(j)},\umbf^{(j)})\ \text{satisfies}\ \eqref{red.penalty}\text{ and }\eqref{red.comp}\ \text{but not } \eqref{red.fea}\},
\eequationn
meaning that $\Ucal$ is the set of subproblem iterations in which \eqref{update.rho} is triggered.  Now we are ready to prove our main result in this section. %{\red The proof of this theorem only considers the case when only condition \eqref{red.penalty} are satisfied by $(d^{(j)}, \umbf^{(j)})$, which is a less strict condition than that \eqref{red.penalty} and \eqref{red.comp} are both satisfied. } 

\begin{theorem}\label{thm.local}
  Suppose Assumptions~\ref{ass.qp} and \ref{ass.algorithm} hold and let 
  \bequationn
    \kappa_4 := \inf_{j\in\Ucal}\{ J^{(0)} -J(d^{(j)},\rho_{(j)})\} \ge 0\ \  \text{and}\ \ \kappa_5 := \inf_{j\in\Ucal} \{J^{(0)}-D(\umbf^{(j)},0)\} \ge 0.
  \eequationn
  Then, for $\rho_{(j)}\in(0,\tilde\rho]$, where 
  \begin{equation}\label{tilde.rho}
    \tilde \rho : = \frac{\omega + \min\{\kappa_4, \kappa_5\}}{\max\{\kappa_2,\kappa_3\}}\left( 1-\sqrt{\beta_v/\beta_\phi}\right),
  \end{equation}
  if $(d^{(j)},\umbf^{(j)})$ satisfies \eqref{red.penalty} and \eqref{red.comp}, then $(d^{(j)},\wmbf^{(j)})$ satisfies \eqref{red.fea}.  In other words, for any $\rho_{(j)}\in(0,\tilde\rho]$, the update \eqref{update.rho} is never triggered by \eqref{dust}. 
\end{theorem}
\begin{proof}
  In order to derive a contradiction, suppose that $\Ucal$ is infinite, meaning that the subproblem solver is never terminated and $\rho_{(j)} \to 0$.  We have from \eqref{J.krho} that 
  \bequationn
    -\kappa_2\rho_{(j)} \le  J(d^{(j)},\rho_{(j)})  -  J(d^{(j)},0)  \le \kappa_2\rho_{(j)}\ \ \text{for any}\ \ j \in \Ucal,
  \eequationn
  which, after adding and dividing through by $J_\omega^{(0)} - J(d^{(j)},\rho_{(j)})$, yields for $j \in \Ucal$ that
\begin{equation}\label{lele2}
 1 -\frac{\kappa_2\rho_{(j)}}{J^{(0)}_\omega-J(d^{(j)},\rho_{(j)})}  \le \frac{J^{(0)}_\omega - J (d^{(j)},0)}{J ^{(0)}_\omega-J(d^{(j)},\rho_{(j)})} \le 1+\frac{\kappa_2\rho_{(j)}}{J^{(0)}_\omega-J(d^{(j)},\rho_{(j)})}.
 \end{equation}
  Thus, for any
  \bequationn
    \rho_{(j)} \le \frac{ \omega + \kappa_4}{\kappa_2}\left(1-\sqrt{\frac{\beta_v}{\beta_\phi}}  \right) \leq \frac{ J^{(0)}_\omega - J(d^{(j)}, \rho_{(j)})}{\kappa_2}\left(1-\sqrt{\frac{\beta_v}{\beta_\phi}}\right),
  \eequationn
  it follows from the first inequality of \eqref{lele2} that
  \begin{equation}\label{satisfy1}
    \frac{J ^0_\omega-J (d^{(j)},0)}{J ^0_\omega-J (d^{(j)},\rho_{(j)})} \geq \sqrt{\frac{\beta_v}{\beta_\phi}}.
  \end{equation}
  Following an argument similar to that for \eqref{D.krho}, we have that for any
  \[\rho_{(j)} \le  \frac{\omega + \kappa_5}{\kappa_3}\left( 1-\sqrt{\frac{\beta_v}{\beta_\phi}}\right) \le  \frac{J^{(0)}_\omega - D(\umbf^{(j)}, 0)}{\kappa_3}\left( 1-\sqrt{\frac{\beta_v}{\beta_\phi}}\right),
  \]
  one finds that
  \begin{equation}\label{satisfy2} 
\frac{J ^0_\omega-D (\umbf^{(j)},\rho_{(j)})}{J ^0_\omega-D (\umbf^{(j)},0)} \geq \sqrt{\frac{\beta_v}{\beta_\phi}}. 
  \end{equation}
  Overall, we have shown that for any $\rho_{(j)} \le \tilde\rho$ with $\tilde\rho$ defined in \eqref{tilde.rho}, it follows that 
\eqref{satisfy1} and \eqref{satisfy2} both hold true and, since $D (\wmbf^{(j)},0) \ge D (\umbf^{(j)},0)$, that 
  \begin{equation}\label{satisfy3}
    \frac{J ^0_\omega-D (\umbf^{(j)},\rho_{(j)})}{J ^0_\omega- D (\wmbf^{(j)},0)} \ge \frac{J ^0_\omega-D (\umbf^{(j)},\rho_{(j)})}{J ^0_\omega-D (\umbf^{(j)},0)} > \sqrt{\frac{\beta_v}{\beta_\phi}}.
  \end{equation}
  
  Since our supposition that $\Ucal$ is infinite implies that $\rho_{(j)} \to 0$, we may now proceed under the assumption that $j \in \Ucal$ with $\rho_{(j)} \in (0,\tilde\rho]$.  Let us now define the ratios
  \bequationn
    \hat r_v^{(j)} := \frac{J^{(0)}_\omega - J(d^{(j)},0) }{J^{(0)}_\omega -  (D(\wmbf^{(j)},0))_+   }\ \ \text{and}\ \ \bar r_v^{(j)} := \frac{J^{(0)}_\omega - J(d^{(j)},0)}{J^{(0)}_\omega -   D(\wmbf^{(j)},0)},
  \eequationn
  where, since $J(d^{(j)},0) = l(d^{(j)},0) + \thalf \langle d^{(j)}, H_0 d^{(j)} \rangle \geq l(d^{(j)},0)$ and by the definition of the operator $(\cdot)_+$, it follows that $r_v^{(j)} \ge \hat r_v^{(j)} \ge \bar r_v^{(j)}$.  From \eqref{satisfy1} and \eqref{satisfy3}, 
  \bequationn
    \frac{\bar r^{(j)}_v}{r^{(j)}_\phi} 
=   \frac{J ^0_\omega-J (d^{(j)},0)}{J ^0_\omega-J (d^{(j)},\rho_{(j)})} \frac{J ^0_\omega-D ( \umbf^{(j)},\rho_{(j)})}{J ^0_\omega-D ( \wmbf^{(j)},0)}  
\ge \frac{\beta_v}{\beta_\phi},
  \eequationn
 yielding 
\[ r^{(j)}_v \ge \bar r^{(j)}_v \ge \frac{\beta_v}{\beta_\phi} r^{(j)}_\phi \ge \beta_v. \]
  However, this contradicts the fact that $j\in\Ucal$.  Overall, since we have reached a contradiction, we may conclude that $\Ucal$ is finite. 
\end{proof}

%************
% Subsection
%************
\section{A Complete Penalty-SQP Algorithm}\label{sec.complete_algorithm}

In the previous section, a dynamic penalty parameter updating strategy was proposed to guarantee that the computed search direction simultaneously offers progress toward reducing the penalty function and reducing infeasibility.  In this section, a complete algorithm for solving \eqref{prob.nlp} that employs this strategy is proposed and analyzed.  It follows the general strategy in Algorithm~\ref{alg.sqo}, but includes additional details.

Our complete algorithm involves an additional check of the penalty parameter after the search direction has been computed as is similarly done in various algorithms that employ a penalty function as a merit function.  Let $\tilde\rho_k$ be the value of the penalty parameter obtained by applying \eqref{dust} within the $k$th subproblem solve.  Then, given a constant $\beta_l \in (0,\beta_\phi(1-\beta_v)]$, we require $\rho_k \in (0,\tilde\rho_k]$ so that
\begin{equation}\label{dust.after}
\Delta l(d^k, \rho_k; x^k) + \omega_k \geq \beta_l (\Delta l(d^k, 0; x^k)+\omega_k),
\end{equation}
where the right-hand side of this inequality is guaranteed to be positive due to \eqref{red.fea}.  More precisely, we employ the following \emph{P}osterior \emph{S}ubproblem \emph{ST}rategy:
\begin{equation}\tag{PSST}\label{psst}
  \boxed{
  \rho_k \gets \begin{cases} \tilde\rho_k  & \text{if this yields \eqref{dust.after}}  \\
   \cfrac{(1-\beta_l) ( \Delta l(d^k, 0; x^k) + \omega_k) }{\langle \nabla f(x^k), d^k\rangle + \tfrac{1}{2} \langle d^k, H(\rho_k; x^k, \eta^k) d^k\rangle } & \text{otherwise.}
 \end{cases}
 }
 \end{equation}
Observe that if the choice $\rho_k = \tilde\rho_k$ does not yield \eqref{dust.after}, then, by setting $\rho_k$ according to the latter formula in \eqref{psst}, it follows (since $H(\rho_k; x^k, \eta^k) \succeq 0$) that
 \[ \rho_k\langle \nabla f(x^k), d^k\rangle \le (1-\beta_l)(\Delta l(d^k, 0; x^k) + \omega_k),\]
 which means that
 \[ \Delta l(d^k, \rho_k; x^k) + \omega_k =  \Delta l(d^k, 0; x^k) - \rho_k\langle \nabla f(x^k), d^k\rangle + \omega_k \ge \beta_l(\Delta l(d^k, 0; x^k) + \omega_k),\]
implying that \eqref{dust.after} holds.
%Therefore, \eqref{dust.after} is satisfied by $\rho=\rho_k$.   In fact, it is obvious that \eqref{dust.after} is true for any $\rho \le \rho_k$, which implies that $\rho_k \in(0, \tilde\rho_k)$.

The intuition  of this posterior updating strategy is to detect whether the iterate may be near an infeasible stationary point.  If a step has  achieved improvement on optimality but not very much on feasibility, then the algorithm should decrease $\rho$ to reduce the effect of the objective in the penalty function.  This is the typical approach used by penalty methods that update the penalty parameter in hindsight at the end of an iteration.  This idea is similar to the updating strategy in  \cite{BurkCurtWang14}.  A novel aspect of \eqref{psst}, however, is that this model reduction condition is imposed inexactly (due to the presence of $\omega_k > 0$). In fact, for a relatively large $\omega_k$, the model reduction in $l(\cdot, \rho_k; x^k)$ is not necessarily at least a fraction of that in $l(\cdot, 0; x^k)$. This difference makes \eqref{psst} more suitable for an inexact penalty-SQP framework. 
 
Our complete algorithm employing \eqref{dust} and \eqref{psst} is given as Algorithm~\ref{alg.sqo.full}.  While we do not complicate the notation by making the dependence explicit on $k \in \N{}$, it should be clear that in the inner loop (over $j$) one is solving a subproblem with quantities dependent on the $k$th iterate; see \eqref{eq.shorthand}.  Also, while our analysis does not depend on this choice, we remark that a reasonable choice for $\eta^{k+1}$ for all $k \in \N{}$ are the \emph{QP multipliers}, i.e., $\eta^{k+1} = \zeta(\umbf^{(j)})$, where $\zeta(\umbf)$ is defined prior to \eqref{eq.zeta}.  We do not specify this choice since one might also consider using, e.g., \emph{least squares multipliers}~\cite{NoceWrig06}.  Our analysis, which focuses on primal convergence, works with any such choice as long as the sequence of dual estimates remains bounded (see below).

In the remainder of this section,  we show that if \eqref{dust} and \eqref{psst} are employed within a penalty-SQP algorithm for solving \eqref{prob.nlp}, then, under reasonable assumptions, the algorithm converges from any starting point.  Specifically, if \eqref{dust} and \eqref{psst} are only triggered a finite number of times, then every limit point of the iterates is either infeasible stationary or first-order stationary for \eqref{prob.nlp}.  Otherwise, if \eqref{dust} and \eqref{psst} are triggered an infinite number of times, driving the penalty parameter to zero, then every limit point of the iterates is either an infeasible stationary point or a feasible point at which a constraint qualification fails to hold.

For our analysis in this section, we extend our use of the sub/superscript $k$ to denote the value of quantities associated with iteration $k \in \N{}$.  For example, $\Ucal^k$ denotes the set $\Ucal$ defined in \S\ref{sec.bound.rhoj} while solving the $k$th subproblem and $\kappa_{0, k}$ is the constant $\kappa_0$ in Assumption~\ref{ass.algorithm} for the $k$th subproblem.

\begin{algorithm}[t]
  \caption{Penalty-SQP with a Dynamic Penalty Parameter Updating Strategy}
  \label{alg.sqo.full}
  \balgorithmic[1]
    \Require $(\gamma, \theta_\rho, \theta_\alpha, \theta_\omega,\beta_v,\beta_\phi) \in (0,1)$, $\beta_l \in (0,\beta_\phi(1-\beta_v))$, and $(\rho_{-1},\omega_0) \in (0,\infty)$
    \State Choose $(x^0,\eta^0) \in \R^n \times \R^m$.
    \For{$k \in \N{}$}
     \State Set $\rho_{(0)} \gets \rho_{k-1}$% and $d^{(0)} \gets 0$
      \For{$j \in \N{}$} 
      
      \State Generate a primal-dual feasible  solution estimate $(d^{(j)}, \umbf^{(j)}, \wmbf^{(j)})$
      %\If{ $r_\phi^{(j)} \ge \beta_\phi$, $r_c^{(j)} \ge \beta_v$ and $r_v^{(j)} < \beta_v$ }{ set $\rho_{({j+1})} \gets \theta_\rho \rho_{(j)}$.} \EndIf
      %\If{ $r_\phi^{(j)} \ge \beta_\phi$, $r_c^{(j)} \ge \beta_v$ and $r_v^{(j)} \ge \beta_v$ }{ terminate} \EndIf
      \State Set $\rho_{(j+1)}$ by applying \eqref{dust}
      \EndFor
      \State Set $d^k \gets d^{(j)}$ and $\tilde \rho_k \gets \rho_{(j)}$.
      %\If{ \[\Delta p(d^k, \tilde \rho_k; x^k) + \omega_k \ge \beta_l (\Delta p(d^k, 0; x^k)+\omega_k)\]}
      %{Set $\rho_k = \tilde\rho_k$.}
      %\Else{ Set  
 %\[\rho_k = \frac{(1-\beta_l) ( \Delta p(d^k, 0; x^k) + \omega_k) }{\langle f(x^k), d^k\rangle + \tfrac{1}{2} \langle d^k, H(\rho_k; x^k, \eta^k) d^k\rangle }.\]}
 %\EndIf
      \State Set $\rho_k$ by applying \eqref{psst}
      \State Let $\alpha^k$ be the largest value in $\{\gamma^0,\gamma^1,\gamma^2,\dots\}$ such that 
      \bequation\label{line_search}
        \phi(x^k+\alpha_k d^k, \rho_k) - \phi(x^k, \rho_k) \le -\theta_\alpha\alpha_k \Delta l(d^k, \rho_k; x^k).
      \eequation
     \State  Choose $\omega_{k+1} \in (0, \theta_\omega\omega_k]$. 
      \State Set $x^{k+1} \gets x^k + \alpha_kd^k$ and choose $\eta \in \R^m$.
    \EndFor
  \ealgorithmic
\end{algorithm} 

We make the following assumption throughout this analysis.
  
\bassumption\label{global}
  The   compact convex set $X \subset \R^n$ with $0 \in \interior{X}$   is used in defining all subproblems, and there exist positive scalar constants 
    $\underline\Lambda, \overline\Lambda$ 
  %$\underline\Lambda_0$, $\overline\Lambda_0$, $\underline\Lambda_\rho$, $\overline\Lambda_\rho$, 
  and $K_0$ with $\underline\Lambda \le \overline\Lambda$  such that the following hold true.
  \benumerate
    \item[(i)] $f$ and $c_i$ for all $i \in \{1,\dots,m\}$, and their first- and second-order derivatives, are all  bounded in an open convex set containing $\{x^k\}$ and $\{x^k+d^k\}$.
    \item[(ii)]% There exist constants $\overline\Lambda \ge \underline\Lambda > 0$ such that 
    For all $k\in\N{}$ and  any $\rho\in[0, \rho_0]$,  
    \[ 0< \underline\Lambda \le \underline\lambda_{0,k} \leq \overline\lambda_{0,k} \leq \overline\Lambda\quad \text{and}\quad 0 < \underline\Lambda \leq \underline\lambda_{\rho,k} \leq \overline\lambda_{\rho,k} \le \overline\Lambda.\]
    \item[(iii)] $\kappa_{0,k} \leq K_0$ for all $k \in \N{}$. 
    \item[(iv)] $\|\nabla c_i(x^k)\|_2 > 0$ for all $k \in \N{}$ and $i \in \{1,\dots,m\}$.
    \item[(v)] $\{\eta^k\}$ is bounded. 
  \eenumerate
\eassumption

Recalling Lemmas~\ref{lem.bounds} and \ref{lem primal dual value}, it follows under Assumption~\ref{ass.qp}, \ref{ass.algorithm}, and \ref{global} that there exist positive scalar constants $K_1$, $K_2$, and $K_3$ such that
\begin{equation}\label{k3bd}
  0 < \kappa_{1,k} \leq K_1,\ \ 0<\kappa_{2,k} \le K_2,\ \ \text{and}\ \ 0 < \kappa_{3,k} \le K_3\ \ \text{for all}\ \ k \in \N{}.
\end{equation}

Let us define the index set 
\bequationn
  \Dcal := \{k \in \N{} : \Ucal^k\ne\emptyset \}.
\eequationn
Moreover, for every $k\in\Dcal$, let $j_k$ be the subproblem iteration number corresponding to the value of the smallest ratio $r_v$, i.e., such that 
\[   r_v^{(j_k)} \le r_v^{(i_k)} \quad \text{for any }\  i_k\in\Ucal^k.\]
Let us also define the index set 
\[ 
\Tcal := \{ k \in\N{}: \rho_k \  \text{is reduced by}\ \eqref{psst} \}.
\]
It follows from these definitions that $\rho_k < \rho_{k-1}$ if and only if $k\in\Dcal\cup\Tcal$.

Before analyzing the behavior of the iterates of our algorithm, we first provide a couple results related to our subproblem and its solutions.  For this result and the remainder of this section, let $d^*(\rho; x, \eta)$ denote a minimizer of $J(d, \rho; x, \eta)$.  From  \cite[Lemma 4.2, 4.3, and 4.4]{BurkCurtWang14}, we have the properties stated in the following lemma.
\begin{lemma}\label{lem.subproblem}
  Under Assumption~\ref{global}, the following hold at any $(x^k,\eta^k)$.
\begin{enumerate}
\item[(i)] The minimizer of $J(\cdot,\rho; x^k, \eta^k)$ is unique for any $\rho\ge0$. 
\item[(ii)]  $\Delta l ( d^*(0,x^k,\eta^k); x^k) \ge 0$ where equality holds if and only if $d^*(0; x^k,\eta^k) = 0$. 
\item[(iii)]  $d^*(0; x^k,\eta^k) = 0$ if and only if $x^k$ is stationary for $v$. 
\item[(iv)] If $d^*(\rho;x^k,\eta^k)=0$ for $\rho > 0$ and $v(x^k)=0$, then $x^k$ is stationary for \eqref{prob.nlp}. 
\end{enumerate}
\end{lemma}

We also have the following fact about the subproblem solutions.
 
\begin{lemma}\label{lem.dk_is_bounded}
  Under Assumption~\ref{global}, $\{ d^*(0; x^k, \eta^k) \}$ and $\{ d^*(\rho_k; x^k, \eta^k)\}$ are bounded.
\end{lemma}
\begin{proof}
  %Assumption~\ref{global}$(ii)$ ensures that there exists $\overline v$ such that, for all $k \in \N{}$,
%\begin{equation}\label{f_v_bounded}  \|\nabla f(x^k)\| \leq \bar v \quad \text{and} \quad v(x^k) \leq \bar v.\end{equation}
%Assumption~\ref{global}$(iii)$ ensures $ \langle d, H(\rho_k; x^k ,\eta^k) d \rangle  \ge \underline\Lambda\|d\|^2$,
%which, together with \eqref{f_v_bounded}, \eqref{eq.primal_no_worse_than_zero}, and the fact that $\{\rho_k\}$ is monotonically decreasing, yield that, for all $k \in \N{}$, 
%\begin{align*}
%  -\rho_0 \bar v \|d^k\| + \tfrac{1}{2}\underline\Lambda\|d^k\|^2  &\leq \rho_k \langle \nabla f(x^k), d^k\rangle + \tfrac{1}{2} \langle d^k, H(\rho_k; x^k ,\eta^k) d^k \rangle + l(d^k, 0; x^k)\\
% &= J(d^k,\rho_k;x^k,\eta^k) \leq J(0, \rho_k; x^k, \eta^k) =  v(x^k) \leq  \bar v. 
%\end{align*}
%This indicates that $\{\|d_k\|\}$ is bounded since, for all $k \in \N{}$,
%\[ \|d^k\| \leq \(\rho_0\bar v + \sqrt{ \rho_0^2 \bar v^2 + 2\underline\Lambda\bar v }\)/\underline\Lambda.\]
  The proof follows the same line of argument for bounding each primal step in norm as is used in the proof of Lemma~\ref{lem.bounds}, where the facts that
\[ \begin{aligned}
 J(d^*(0; x^k, \eta^k), 0; x^k, \eta^k) & \le J(0, 0;  x^k, \eta^k) \\
\text{and}\quad  J(d^*(\rho_k; x^k, \eta^k), \rho_k; x^k, \eta^k) & \le J(0, 0;  x^k, \eta^k)
\end{aligned} \]
  follow from the definitions of $d^*(0; x^k, \eta^k)$ and $d^*(\rho_k; x^k, \eta^k)$.
\end{proof}

We now prove a useful lower bound for the stepsize in each iteration.

\begin{lemma}\label{lem.stepsize} Under Assumption~\ref{global}, it follows that, for all $k \in \N{}$, the stepsize satisifies 
$\alpha_k \ge C \Delta l(d^k, \rho_k; x^k)$ for some constant $C>0$ independent of $k$. 
\end{lemma}
\begin{proof}
  If $d^k = 0$, then \eqref{line_search} holds with $\alpha^k = \gamma^0 = 1$.  Hence, for the remainder of the proof, let us assume that $d^k \neq 0$.  Under Assumption~\ref{global}, applying Taylor's theorem and \cite[Lemma 4.2]{BurkCurtWang14}, we have that for all positive $\alpha$ that are sufficiently small, there exists $\tau > 0$ such that 
\[ \phi(x^k+\alpha d^k, \rho_k ) - \phi(x^k ,\rho_k)  \le - \alpha \Delta l(d^k,\rho_k; x^k ) + \tau \alpha^2\|d^k\|^2_2.\]
Thus, for any $\alpha \in[ 0, (1-\theta_\alpha)\Delta l(d^k,\rho_k; x^k)/(\tau\|d^k\|^2_2)]$, it follows that
\[ -\alpha \Delta l(d^k, \rho_k; x^k) + \tau\alpha^2\|d^k\|^2_2 \le - \alpha \theta_\alpha \Delta l(d^k, \rho_k; x^k),\]
meaning that the sufficient decrease condition \eqref{line_search} holds. During the line search, the stepsize is multiplied by $\gamma$ until \eqref{line_search} holds, so we know by the above inequality that the backtracking procedure terminates with 
\[ \alpha_k \ge \gamma(1-\theta_\alpha) \Delta l(d^k, \rho_k; x^k) / (\tau\|d^k\|^2_2).\]
The result follows from this inequality since $\{\|d^k\|_2\}$ is bounded above by $K_1$.
\end{proof}

%We now prove that, in the limit, 
Next we show that the reductions in the models of the constraint violation and the penalty function both vanish in the limit. 
For this purpose, it will be convenient to work with the shifted penalty function
\[ \varphi(x,\rho) := \rho ( f(x) - \underline f)  + v(x) \ge 0, \]
where $\underline f$ is the infimum of $f$ over the smallest convex set containing $\{x^k\}$. The existance of $\underline f$ follows from 
Assumption~\ref{global}$(i)$. The function $\varphi$ possesses a useful monotonicity property proved in the following lemma. 
 
\begin{lemma}\label{lem.varphi}
  Under Assumption~\ref{global}, it holds that, for all $k \in \N{}$,
 \[ \varphi(x^{k+1}, \rho_{k+1}) \le \varphi(x^k,\rho_k)  - \theta_\alpha \alpha_k\Delta l(d^k, \rho_k; x^k).\]
 \end{lemma}

\begin{proof}
 By the line search condition \eqref{line_search}, it follows that 
 \[ \varphi(x^{k+1}, \rho_k) \le \varphi(x^k, \rho_k) - \theta_\alpha \alpha_k \Delta l(d^k,\rho_k; x^k),\]
 which implies 
 \[ \varphi(x^{k+1}, \rho_{k+1}) \le \varphi(x^k,\rho_k) - (\rho_k-\rho_{k+1})(f(x^{k+1}) - \underline f) - \theta_\alpha \alpha_k\Delta l(d^k, \rho_k; x^k).\]
 The result then follows from this inequality, the fact that $\{\rho_k\}$ is monotonically decreasing, and since $f(x^{k+1}) \ge \underline f$ for 
 all $k\in\N{}$. 
\end{proof}

We now show that the model reductions and duality gap all vanish asymptotically.

\begin{lemma}\label{lem.delta_2_zero}
  Under Assumption~\ref{global}, the following limits hold.
\begin{enumerate}
\item[(i)] $0 =  \lim\limits_{k\to \infty} \Delta l(d^k, \rho_k; x^k) = \lim\limits_{k\to\infty} \Delta J(d^k,\rho_k; x^k, \eta^k)$,
\item[(ii)] $0 =  \lim\limits_{k\to\infty} \Delta l(d^k, 0; x^k) = \lim\limits_{k\to\infty} \Delta J(d^k, 0; x^k, \eta^k)$,
\item[(iii)] $0 =  \lim\limits_{k\to\infty} \Delta J(d^*(0;x^k,\eta^k), 0; x^k, \eta^k)=  \lim\limits_{k\to\infty} \Delta J(d^*(\rho_k;x^k,\eta^k),\rho_k; x^k, \eta^k)$,
\item[(iv)] $0 = \lim\limits_{k\to\infty} [J(0, \rho_k; x^k, \eta^k) - D(\umbf^k, \rho_k; x^k, \eta^k)],$ 
\item[(v)] $0 = \lim\limits_{k\to\infty} [J(0, 0; x^k, \eta^k) - D(\wmbf^k, 0; x^k, \eta^k)].$
 \end{enumerate}
\end{lemma}

\begin{proof}
  Let us first prove $(i)$ by contradiction.  Suppose that $\Delta l(d^k, \rho_k; x^k)$ does not converge to 0. Then, there exists a constant 
$\epsilon>0$ and an infinite $\Kcal \subseteq \N{}$ such that $\Delta l(d^k, \rho_k; x^k) \geq \epsilon$ for all $k\in \Kcal$. It then follows from Lemma~\ref{lem.stepsize} and \ref{lem.varphi} that $\varphi(x^k; \rho_k) \to -\infty$, which contradicts the fact that $\{\varphi(x^k,\rho_k)\}$ is bounded below by zero. Therefore, $\Delta l(d^k, \rho_k; x^k)\to 0$.  The second limit in $(i)$ then follows from the first limit, the fact that $H(\rho_k;x^k,\eta^k) \succeq 0$ for all $k \in \N{}$, and the fact that
\begin{equation}\label{model_red_J}  
\begin{aligned}
\Delta l(d^k, \rho_k; x^k)  & =  \Delta J(d^k, \rho_k ; x^k, \eta^k) + \tfrac{1}{2} \langle d^k, H(\rho_k; x^k, \eta^k) d^k\rangle \\
 &\ge \Delta J(d^k, \rho_k ; x^k, \eta^k).
\end{aligned}
\end{equation}
Next, from \eqref{dust.after} and \eqref{model_red_J}, it follows that 
\begin{equation*}
\Delta l(d^k, \rho_k; x^k) + \omega_k   \ge \beta_l (\Delta l(d^k, 0; x^k)+\omega_k) \ge \beta_l( \Delta J(d^k, 0; x^k, \eta^k) + \omega_k).   
\end{equation*}
The limits in $(ii)$ follow from these inequalities, the first limit in $(i)$, and the fact that $\{\omega_k\} \to 0$.  Finally, the limits in $(iii)$, $(iv)$, and $(v)$ follow from the limits in parts $(i)$ and $(ii)$ along with the inequalities in \eqref{alter.red.constr} and \eqref{alter.red.penalty}.
\end{proof}

%
%
%\begin{lemma}\label{lem.limits}
%Suppose Assumption~\ref{global} holds. Then, the following limits hold:
%\begin{align}
% J(0,0; x^k, \eta^k) - D(\wmbf^k, \rho_k; x^k, \eta^k)  \to 0\quad\text{and}\quad  G(d^k, \wmbf^k, 0; x^k, \eta^k) \to 0.   \\
%  G(0, \umbf^k, \rho_k; x^k, \eta^k)  \to 0\quad\text{and}\quad  G(0, \wmbf^k, 0; x^k, \eta^k) \to 0.   
%\end{align}
%\end{lemma}
%
%\begin{proof}
%By Lemma~\ref{lem.delta_2_zero} (a), Assumption~\ref{global} (v) and \eqref{alter.red.penalty}, it holds  that  
%\[  J(0, \rho_k; x^k, \eta^k) - D(\umbf^k, \rho_k; x^k, \eta^k)  \to 0. \]
%
%By Lemma~\ref{lem.delta_2_zero} (b), Assumption~\ref{global} (v) and \eqref{alter.red.constr}, it holds  that  
%\[  J(0, \rho_k; x^k, \eta^k) - D(\umbf^k, \rho_k; x^k, \eta^k)  \to 0. \]
%
%
%By Lemma~\ref{lem.delta_2_zero}, $\Delta p(d^k, \rho_k; x^k) = J(0,\rho_k; x^k, \eta^k) - p(d^k, \rho_k; x^k) \to 0$, which together with Assumption~\ref{global} (v) and \eqref{alter.red.constr}, implies 
%\[  J(d^k, 0; x^k, \eta^k) - D(\umbf^k, 0; x^k, \eta^k)  \to 0 \quad\text{and}\quad   J(d^k, 0; x^k, \eta^k) - J(d_\rho_k, 0; x^k, \eta^k)\to0 \]
%\end{proof}
%
%

We now show that the primal steps and the exact subproblem solutions vanish.

\begin{lemma}\label{lem.d_to_zero}
  Suppose Assumption~\ref{global} holds and $\{\rho_k\} \to \rho_*$.  Then, $\{d^k\} \to 0$ and for any limit point $x^*$ of $\{x^k\}$ it follows that $d^*(0; x^*, \cdot) = 0$ and $d^*(\rho_*; x^*, \cdot) = 0$.
\end{lemma}

\begin{proof}
  From Lemma~\ref{lem.delta_2_zero}$(ii)$, it follows that
  \[\begin{aligned} 0 = & \lim_{k\to\infty}- \Delta J(d^k, 0; x^k, \eta^k) = \lim_{k\to\infty} - \Delta l(d^k, 0; x^k) + \thalf \langle d^k,  H(0; x^k, \eta^k) d^k \rangle \\
 = & \lim_{k\to\infty}   \thalf \langle d^k,  H(0; x^k, \eta^k) d^k \rangle \ge \lim_{k\to\infty}   \thalf \underline{\Lambda} \|d^k\|^2_2.
\end{aligned} \]
This implies that $\{d^k\} \to 0$, as desired.  Next, from Lemma~\ref{lem.delta_2_zero}$(iii)$ and continuity, it follows that $\Delta J(d^*(0; x^*, \cdot),0;x^*, \cdot)=0$, from which it follows that
\begin{equation*}
  J(d^*(0; x^*, \cdot),0;x^*, \cdot) = J(0, 0;x^*, \cdot).
\end{equation*}  From the strong convexity of 
$J(\cdot,0;x^*, \cdot)$ and the fact that $d^*(0; x^*, \cdot)$ is its minimizer, it follows that $d^*(0; x^*, \cdot)=0$.  Using a similar argument and Lemma~\ref{lem.delta_2_zero}$(iii)$ again, it follows that $d^*(\rho_*; x^*, \cdot)=0$, completing the proof.
%Similarly, from Lemma~\ref{lem.delta_2_zero}$(iii)$, it holds that $\Delta J(d^*(\rho_*; x^*, \eta^*),\rho_*;x^*, \eta^*)=0$. 
%This implies that $J(d^*(\rho_*; x^*, \eta^*),0;x^*, \eta^*)=J(0, \rho_*;x^*, \eta^*)$. From the strong convexity of 
%$J(\cdot, \rho_*;x^*, \eta^*)$ and that $d^*(\rho_*; x^*, \eta^*)$ is the minimizer of $J(\cdot,\rho_*;x^*, \eta^*)$, it holds true that $d^*(\rho_*; x^*, \eta^*)=0$.
 %\[\begin{aligned} 0 = & \lim_{k\to\infty}- \Delta J(d^*(\rho_k; x^k, \eta^k), \rho_k; x^k, \eta^k)\\
%= & \lim_{k\to\infty} - \Delta l(d^*(\rho_k; x^k, \eta^k), \rho_k; x^k) + \thalf \langle d^*(\rho_k; x^k, \eta^k),  H(\rho_k; x^k, \eta^k) d^*(\rho_k; x^k, \eta^k) \rangle \\
%  & \textbf{(FEC: How does this follow? Lemma~\ref{lem.delta_2_zero} does not talk about $\Delta l(d^*)$.)} \\
%= & \lim_{k\to\infty}   \thalf \langle d^*(\rho_k; x^k, \eta^k),  H(\rho_k; x^k, \eta^k) d^*(\rho_k; x^k, \eta^k) \rangle \\
%\ge & \lim_{k\to\infty}   \thalf \underline{\Lambda} \|d^*(\rho_k; x^k, \eta^k)\|^2,
%\end{aligned} \]
%where the last inequality follows by Assumption~\ref{global}$(iii)$.  This implies the third desired limit.  The proofs for the other limits follow the same line of argument.
\end{proof}

%We now provide 
Our first global convergence theorem follows.
 
\begin{theorem}\label{thm.global1}
 Under Assumption~\ref{global}, the following statements hold.  
\begin{enumerate}
\item[(i)] Any limit point of $\{x^k\}$ is first-order stationary for $v$, i.e., it is feasible or an infeasible stationary point for \eqref{prob.nlp}. 
\item[(ii)] If $\rho_k \to \rho_*$ for some $\rho_*>0$ and $v(x^k)\to0$, then any limit point $x^*$ of $\{x^k\}$ with $v(x^*)=0$ is a KKT point for \eqref{prob.nlp}. 
\item[(iii)] If $\rho_k \to 0$, then either all limit points of $\{x^k\}$ are feasible for \eqref{prob.nlp} or all are infeasible. 
\end{enumerate}
\end{theorem}
\begin{proof}
  Part $(i)$ follows by combining Lemma~\ref{lem.d_to_zero} with Lemma~\ref{lem.subproblem}$(iii)$.  Similarly, part $(ii)$ follows by combining Lemma~\ref{lem.d_to_zero} with Lemma~\ref{lem.subproblem}$(iv)$.

We prove $(iii)$ by contradiction.  Suppose there exist infinite $\Kcal^* \subseteq \N{}$ and $\Kcal^\times  \subseteq \N{}$ such that $\{x^k\}_{k\in\Kcal^*}\to x^*$ with $v(x^*)=0$ and 
$\{x^k\}_{k\in\Kcal^\times} \to x^\times$ with $v(x^\times) = \epsilon > 0$.  Since $\rho_k \to 0$, there exists $k^* \ge 0$ such that for all $k\in\Kcal^*$ and $k\ge k^*$ one has that 
$\rho_{k+1}(f(x^k)-\underline f) < \epsilon/4$ and $v(x^k) < \epsilon/4$, meaning that $\varphi(x^k, \rho_{k+1}) < \epsilon/2$.  On the other hand, it follows that $\rho_{k+1}(f(x^k)-\underline f) \ge 0$ for all $k \in \N{}$ and there exists $k^\times \in \N{}$ such that $v(x^k)\ge \epsilon/2$ for all $k\ge k^\times$ with $k \in \Kcal^\times$, meaning that $\varphi(x^k, \rho_{k+1}) \ge \epsilon/2$.  This contradicts Lemma~\ref{lem.varphi}, which shows that $\varphi(x^k, \rho_{k+1})$ is monotonically decreasing.  Thus, the set of limit points of $\{x^k\}$ must be all feasible or all infeasible. 
\end{proof}

Theorem~\ref{thm.global1} is satisfactory in the case when $\rho_k \to \rho_* > 0$, 
%in which case it is proved 
since it shows that any limit point of the primal sequence is a KKT point for \eqref{prob.nlp}.  
%However, more should 
But more needs to be said  %in the case that 
when $\rho_k \to 0$. 
%in particular, in the following results, we look further at 
We now address this case, showing that it only occurs if a limit point of the algorithm is either an infeasible stationary point or a feasible point at which a constraint qualification fails to hold.  We begin with the following lemma.

\begin{lemma}\label{lem.dual}
Suppose Assumption~\ref{global} holds and $\rho_k \to 0$. Let $x^*$ be a limit point of $\{x^k\}_{k\in\Dcal\cup\Tcal}$ that is feasible for \eqref{prob.nlp} with infinite $\Scal\subseteq\Dcal\cup\Tcal$ such that $\{x^k\}_{k\in\Scal} \to x^*$.  Then, the following hold true.
\begin{enumerate}
\item[(i)] $|\Scal\cap\Dcal|$ is finite or $\{ \Delta J(d^{(j_k)},\rho_{(j_k)}; x^k, \eta^k)\}_{k\in\Scal\cap\Dcal} \to 0$;
\item[(ii)] $|\Scal\cap\Dcal|$ is finite or $\{ d^{(j_k)} \}_{k\in\Scal\cap\Dcal} \to 0$;
\item[(iii)] any limit point of $\{\umbf^{(j_k)}\}_{k\in\Scal \cap\Dcal} \cup \{\umbf^k\}_{k\in\Scal\cap\Tcal}$ is optimal for $D(\cdot, 0; x^*, \cdot)$;
%\item[(d)]  $\Ecal \cup \Acal(x^*) \ne \emptyset$. 
\item[(iv)] $\{\umbf^{(j_k)}\}_{k\in\Scal \cap\Dcal} \cup \{\umbf^k\}_{k\in\Scal\cap\Tcal}$ has a nonzero limit point.
\end{enumerate}
\end{lemma}
\begin{proof}
  For part $(i)$, if $|\Scal\cap\Dcal|$ is finite, then there is nothing left to prove.  Hence, let us assume that $|\Scal\cap\Dcal| = \infty$.  Observe that, for all $k \in \N{}$, it holds that
    \begin{equation*}
    \begin{aligned}
    0 \le &\     \Delta J(d^{(j_k)},\rho_{(j_k)}; x^k, \eta^k) \\
    = &\ v(x^k) -   \rho_{(j_k)} \langle  \nabla f(x^k), d^{(j_k)}\rangle - \tfrac{\rho_{(j_k)}}{2} \langle
    d^{(j_k)}, H_f(x^k)d^{(j_k)}\rangle - J(d^{(j_k)}, 0; x^k, \eta^k)\\
      \le &\ v(x^k) -   \rho_{(j_k)} \langle  \nabla f(x^k), d^{(j_k)}\rangle - \tfrac{\rho_{(j_k)}}{2} \langle
    d^{(j_k)}, H_f(x^k)d^{(j_k)}\rangle,
     \end{aligned}
  \end{equation*}
  where the first inequality follows from \eqref{eq.primal_no_worse_than_zero} and the second inequality follows from the definition of $J$, which ensures that $J(d^{(j_k)}, 0; x^k, \eta^k)\ge 0$.  In addition, $\{d^{(j_k)}\}$ is bounded due to Lemma~\ref{lem.bounds} and Assumption~\ref{global}$(ii)$-$(iii)$.  Consequently, since $|\Scal\cap\Dcal| = \infty$ and $\{v(x^k)\}_{k\in\Scal\cap\Dcal} \to 0$ with $\rho_{(j_k)}\to0$, the limit in part $(i)$ holds. 
  
  For part $(ii)$, again, if $|\Scal\cap\Dcal|$ is finite, then there is nothing left to prove.  Otherwise, since $\{J(0,0; x^k, \eta^k)\}_{k\in\Scal\cap\Dcal}=\{v(x^k)\}_{k\in \Scal\cap\Dcal}\to0$ and $\rho_{(j_k)}\to0$, the limit in part~$(ii)$ holds due to Lemma~\ref{lem.bounds} and Assumption~\ref{global}$(ii)$-$(iii)$.
  
  Now consider part $(iii)$.  If $|\Scal\cap\Dcal|$ is infinite, then for a limit point $\umbf^*$ there must exist an infinite $\Scal_\Dcal \subseteq \Scal \cap \Dcal$ such that $\{ \umbf^{(j_k)}\}_{k\in\Scal_D} \to \umbf^*$.  Then, % from Lemma~\ref{lem.delta_2_zero}$(iv)$,
  it follows that 
  \bequation\label{limlimlim}
    \baligned
      0 &\leq J(0, 0; x^*, \cdot) - D(\umbf^*, 0 ;x^*, \cdot) \\ &= \lim_{  k\in\Scal_\Dcal \atop k\to\infty} J(0,\rho_{(j_k)}; x^k, \cdot) -D(\umbf^{(j_k)}, \rho_{(j_k)}; x^k, \cdot)\\
     & \le \lim_{  k\in\Scal_\Dcal \atop k\to\infty}  \beta_\phi  [J(0,\rho_{(j_k)}; x^k, \cdot) - J(d^{(j_k)},\rho_{(j_k)}; x^k, \cdot)]\\
     & = \lim_{  k\in\Scal_\Dcal \atop k\to\infty}  \beta_\phi  [J(0,0; x^k, \cdot) - J(d^{(j_k)},0; x^k, \cdot)] \le \lim_{  k\in\Scal_\Dcal \atop k\to\infty}  \beta_\phi  J(0,0; x^k, \cdot) = 0,
    \ealigned
  \eequation
  where the second inequality is by \eqref{red.penalty} and the third inequality is by the fact that $J(d^{(j_k)},0; x^k, \cdot) \ge 0$. This 
means that $\umbf^*$ is optimal for $D(\cdot, 0; x^*, \cdot)$.  On the other hand, if $|\Scal\cap\Dcal|$ is finite, then $|\Scal\cap\Tcal|$ must be infinite, in which case for a limit point $\umbf^*$ there must exist an infinite $\Scal_\Tcal \subseteq \Scal \cap \Tcal$ such that $\{ \umbf^k\}_{k\in\Scal_\Tcal} \to \umbf^*$.  Then, again from Lemma~\ref{lem.delta_2_zero} and \eqref{limlimlim}, it follows that $\umbf^*$ is optimal for $D(\cdot, 0; x^*, \cdot)$.

  For part $(iv)$, first observe that 
 \[ l(d, 0; x^k) = \sum_{i\in\Ecal_+(d)\cup\Ecal_-(d)\cup\Ical_+(d)} \|\nabla c_i(x^k)\|_2 \dist{d}{C_i^k},\]
 and that $\chi(d,\umbf; x^k)$ can be viewed as a weighted variant of this sum with weights
 \[1-\zeta_i(\umbf)\ \ \text{for all}\ \ i \in \Ecal_+(d)\cup\Ical_+(d)\ \ \text{and}\ \ 1+\zeta_i(\umbf)\ \ \text{for all}\ \ i\in \Ecal_-(d).\]
  Also observe that \eqref{red.comp} holds at any primal-dual point
  \[ (d,\umbf) \in \{(d^{(j_k)}, \umbf^{(j_k)})\}_{k\in\Scal\cap\Dcal} \cup \{(d^k, \umbf^k)\}_{k\in\Scal\cap\Tcal}\]
  due to the facts that
\begin{align}
  \chi(d^{(j_k)}, \umbf^{(j_k)}; x^k) &\le (1-\beta_v)^2 (v(x^k)+\omega_k)\ \ \text{for all}\ \  k\in\Scal\cap\Dcal\ \ \text{and}\label{comp.1}\\
  \chi(d^k, \umbf^k; x^k) & \le (1-\beta_v)^2 (v(x^k)+\omega_k)\ \ \text{for all}\ \  k\in\Scal\cap\Tcal.\label{comp.2}
\end{align}
We now consider three cases.

\bitemize
  \item[Case $(a)$:] Assume there exists  an infinite $\Scal_\Dcal \subseteq \Scal\cap\Dcal$ such that    
 \begin{equation}
 \label{p_too_large_d}
 l(d^{(j_k)}, 0; x^k) > (1-\beta_v) (v(x^k)+\omega_k)\ \ \text{for all}\ \ k\in\Scal_D.
  \end{equation}
  Then, $\|\zeta(\umbf^{(j_k)}) \|_\infty \ge \beta_v$ for all $k\in\Scal_\Dcal$; indeed, if this were not the case, then for some $k \in \Scal_\Dcal$ one would find from the definition of $\chi$ and \eqref{p_too_large_d} that
 \[\chi(d^{(j_k)}, \umbf^{(j_k)}; x^k) \geq (1-\beta_v)  l(d^{(j_k)}, 0; x^k) > (1-\beta_v)^2 (v(x^k)+\omega_k), \]
  contradicting \eqref{comp.1}.  In this case, combining Lemma~\ref{lem.bounds.dual}, Assumption~\ref{global}$(iv)$, and the fact that $\|\zeta(\umbf^{(j_k)}) \|_\infty \ge \beta_v$ for all $k\in\Scal_\Dcal$ shows that $\{\umbf^{(j_k)}\}_{k\in\Scal\cap\Dcal}$ has a nonzero limit point, proving part $(iv)$, as desired.  
  
  \item[Case $(b)$:] 
Assume there exists  an infinite $\Scal_\Tcal \subseteq \Scal\cap\Tcal$ such that 
 \begin{equation}
 \label{p_too_large_t}
 l(d^{k}, 0; x^k)  > (1-\beta_v) (v(x^k)+\omega_k)\ \ \text{for all}\ \ k\in\Scal_\Tcal.
  \end{equation}
  Then, $\|\zeta(\umbf^k)\|_\infty \ge \beta_v$ for all $k\in\Scal_\Tcal$; indeed, if this were not the case, then for some $k\in\Scal_\Tcal$ one
would find from the definition of $\chi$ and \eqref{p_too_large_d} that
 \[\chi(d^k, \umbf^k; x^k) \geq (1-\beta_v)  l(d^k, 0; x^k) > (1-\beta_v)^2 (v(x^k)+\omega_k), \]
  contradicting \eqref{comp.2}.  In this case, combining Lemma~\ref{lem.bounds.dual}, Assumption~\ref{global}$(iv)$, and the fact that $\|\zeta(\umbf^k) \|_\infty \ge \beta_v$ for all $k\in\Scal_\Tcal$ shows that $\{\umbf^k\}_{k\in\Scal\cap\Tcal}$ has a nonzero limit point, proving part $(iv)$, as desired.  
 % \textbf{Similarly, if \eqref{p_too_large_t} holds for infinite~$\Scal_\Tcal$, then part $(iv)$ follows using a similar argument and \eqref{comp.2}.}
  \item[Case $(c)$:] Suppose that \eqref{p_too_large_d} and \eqref{p_too_large_t} only hold for finite subsets of $\Scal\cap\Dcal$ and $\Scal\cap\Tcal$.  In this case, there exists a sufficiently large $\bar k \in \N{}$ such that 
 \begin{align}
  l(d^{(j_k)}, 0; x^k) \leq (1-\beta_v) (v(x^k)+\omega_k) &\ \text{for all}\ k\in\Scal\cap\Dcal\ \text{with}\ k \geq \bar k; \label{chi_less_d}\\ 
  l(d^{k}, 0; x^k) \leq (1-\beta_v) (v(x^k)+\omega_k) &\ \text{for all}\ k\in\Scal\cap\Tcal\ \text{with}\ k \geq \bar k. \label{chi_less_t}
 \end{align}
 We can further assume that 
 \bequationn
   \baligned
     \|\zeta(\umbf^{(j_k)})\|_\infty &< \beta_v\ \text{for all}\ k\in\Scal\cap\Dcal\ \text{with}\ k \geq \bar k\ \text{and} \\
     \|\zeta(\umbf^k)\|_\infty &< \beta_v\ \text{for all}\ k\in\Scal\cap\Tcal\ \text{with}\ k \geq \bar k;
   \ealigned
 \eequationn
 since otherwise, as in Cases $(a)$ and $(b)$, respectively, part $(iv)$ would hold.  Now, for $k \geq \bar k$ with $k\in\Scal\cap\Dcal$, it follows from \eqref{chi_less_d} that
 % similar to \eqref{J_D_p}, we have 
 % \begin{equation}\label{J_D_p_d}
  % J(0,0; x^k, \eta^k)  - D(\umbf^{(j_k)}, \rho_{(j_k)}; x^k, \eta^k) \ge (1-\beta_v) v(x^k).\end{equation} 
\[
\begin{aligned}
      &\ J(0,0;x^k,\eta^k) +\omega_k - l(d^{(j_k)},0; x^k) \\ 
  \geq&\ v(x^k)+\omega_k - (1-\beta_v) (v(x^k)+\omega_k) \\
     =&\ \beta_v( v(x^k)+\omega_k)\\
  \geq&\ \beta_v [v(x^k)+\omega_k - (D(\wmbf^{(j_k)}, 0; x^k, \eta^k))_+],
\end{aligned}
\]
from which it follows that
\[
 r_v^{(j_k)} = \frac{J(0,0;x^k,\eta^k)+\omega_k - l(d^{(j_k)},0; x^k)}{v(x^k)+\omega_k - (D(\wmbf^{(j_k)}, 0; x^k, \eta^k))_+} \geq \beta_v.
\]
  This indicates that \eqref{dust} is not triggered at any iteration $k \geq \bar k$ with $k\in\Scal\cap\Dcal$.  By the definition of $\Dcal$, this implies that $\Scal \cap \Dcal$ is finite.  On the other hand, for $k\in\Scal\cap\Tcal$ with $k \geq \bar k$,  it holds that 
\begin{equation}\label{J_D_p_t}
\begin{aligned} 
      & J(0,0; x^k, \eta^k)  - D(\umbf^k, \rho_k; x^k, \eta^k) \\ 
\ge &\ v(x^k)   + \sum_{i=1}^m \|\nabla c_i(x^k)\|_2 \delta^*(u_i^k|C_i^k) \\
  = & \ \sum_{i=1}^{\bar m}|c_i(x^k)| + \sum_{i=\bar m+1}^m (c_i(x^k))_+   - \sum_{i=1}^m  \|\nabla c_i(x^k)\|_2 \zeta^i(\umbf^k)  \frac{c_i(x^k)}{\|\nabla c_i(x^k)\|_2} \\
= & \   \sum_{i=1}^{\bar m}|c_i(x^k)| + \sum_{i=\bar m+1}^m (c_i(x^k))_+   - \sum_{i=1}^m  \zeta^i(\umbf^k)  c_i(x^k)  \\
= & \   \sum_{i=1}^{\bar m}[|c_i(x^k)| - \zeta^i(\umbf^k) c_i(x^k)] + \sum_{i=\bar m+1}^m [(c_i(x^k))_+   -    \zeta^i(\umbf^k) c_i(x^k)]  \\
\ge & \   \sum_{i=1}^{\bar m}(1-|\zeta^i(\umbf^k)|) |c_i(x^k)|   + \sum_{i=\bar m+1}^m (1-|\zeta^i(\umbf^k)|)(c_i(x^k))_+    \\
\ge & \   (1-\beta_v) \sum_{i=1}^{\bar m} |c_i(x^k)|   +  (1-\beta_v) \sum_{i=\bar m+1}^m(c_i(x^k))_+ = (1- \beta_v) v(x^k),
 \end{aligned}
 \end{equation}
 where the first inequality is  from the positive definiteness of $H(0,x^k,\eta^k)$ and $\delta^*(u_{m+1}^k|X)=\sup_{d\in X}\langle u_{m+1}^k, d\rangle \ge 0$, and the first equality is from \eqref{eq.zeta}.  Since \eqref{red.penalty} is satisfied, the first inequality in \eqref{alter.red.penalty} and \eqref{J_D_p_t} imply 
\[ \begin{aligned}
&\Delta J(d^k, \rho_k; x^k, \eta^k) + \omega_k = J(0,0; x^k, \eta^k)  - J(d^k, \rho_k; x^k, \eta^k) + \omega_k \\
\ge &\ \beta_\phi [J(0,0; x^k, \eta^k)   - D(\umbf^k, \rho_k; x^k, \eta^k)  + \omega_k]  \\ 
\ge & \  \beta_\phi[(1-\beta_v)v(x^k)+ \omega_k] \ge \beta_\phi(1-\beta_v) (v(x^k)+ \omega_k)\\
\ge  & \  \beta_l (v(x^k)+ \omega_k) \ge \beta_l (\Delta l(d^k, 0; x^k)+\omega_k) .
 \end{aligned}
 \]
which, together with \eqref{model_red_J}, yields
 \[ \Delta l(d^k, \rho_k; x^k) + \omega_k  \ge  \Delta J(d^k,\rho_k;x^k, \eta^k)+ \omega_k  \ge   \beta_l (\Delta l(d^k,0;x^k)+ \omega_k).\]
Therefore, \eqref{psst} is not triggered in any iteration $k\in \Scal\cap\Tcal$ with $k \geq \bar k$.  By the definition of $\Tcal$, this means that  
$\Scal\cap\Tcal$ is finite.  Overall, we have shown in this case that $\Scal \cap \Dcal$ and $\Scal \cap \Tcal$ are finite, meaning $\Scal$ is finite.  However, this contradicts the statement of the lemma, which defines $\Scal$ to be infinite.
 \eitemize
Overall, since Case $(c)$ leads to a contradiction, it follows that either Case $(a)$ or $(b)$ must occur, which proves part $(iv)$. 
\end{proof}

We are now prepared to prove a theorem about the behavior of the algorithm when the penalty parameter is driven to zero.  The theorem involves a statement about points satisfying the well-known Mangasarian-Fromovitz constraint qualificatioin (MFCQ).  Defining $\Ecal = \{1,\ldots, \bar m\}$,  $\Ical = \{\mbar+1,\dots,m\}$,
\[
\begin{aligned}  \Acal(x) & = \{i\in\{\bar m+1, \ldots, m\} : c_i(x) = 0\},\\ \text{and}\ \ 
 \Ncal(x) & = \{i\in\{\bar m+1, \ldots, m\} : c_i(x) < 0\},
\end{aligned}
\] 
we now recall this qualification then state and prove our theorem.

\begin{definition} A point $x$ satisfies the MFCQ for problem \eqref{prob.nlp} if $v(x) = 0$, 
$\{\nabla c_i(x) : i \in \Ecal\}$ are linearly independent, and there exists $d\in\mathbb{R}^n$ such that 
\[ \baligned c_i (x) + \langle \nabla c_i(x), d\rangle &= 0\ \ \text{for all}\ \ i\in\Ecal \\ \text{and} \quad c_i(x) + \langle \nabla c_i(x), d\rangle &< 0\ \ \text{for all}\ \ i\in\Ical, \ealigned \]
or, equivalently, 
\[ \langle \nabla c_i(x), d \rangle =0\ \ \text{for all}\ \ i\in\Ecal\quad \text{and}\quad \langle \nabla c_i(x), d\rangle < 0\ \ \text{for all}\ \ i\in\Acal(x).\]
\end{definition}
The dual form \cite{Solodov10} of MFCQ states that $\zeta^i=0, i\in\Ecal\cup\Acal(x)$ is the unique solution of the linear system 
\[   \sum_{i\in\Ecal\cup\Acal(x)} \zeta^i \nabla c_i(x) = 0, \  \zeta^i \ge 0, i\in \Acal(x).\]

\begin{theorem}\label{thm.global bd}
  Suppose Assumption~\ref{global} holds and $\rho_k\to 0$.  Then,
every limit point of $\{x^k\}_{k\in\Dcal\cup\Tcal}$ is either an infeasible stationary point or a feasible point where the MFCQ does not hold.
\end{theorem}
\begin{proof} 
  By Theorem~\ref{thm.global1}$(i)$, any limit point of $\{x^k\}_{k\in\Dcal\cup\Tcal}$ is either feasible or an infeasible stationary point.  If any such point is infeasible, then there is nothing left to prove.  We may thus proceed by letting $x^*$ represent a feasible limit point of $\{x^k\}_{k\in\Dcal\cup\Tcal}$.  Our goal is to show that the MFCQ fails to hold at $x^*$.

  Let $\Scal\subseteq\Dcal\cup\Tcal$ be an infinite set such that $\{ x^k\}_{k\in\Scal} \to x^*$.  By Theorem~\ref{lem.dual}$(iv)$, it follows that there exists a nonzero limit point $\umbf^*$ of $\{\umbf^{(j_k)}\}_{k \in \Scal\cap\Dcal} \cup \{\umbf^k\}_{k \in \Scal\cap\Tcal}$.  In addition, from Lemma~\ref{lem.d_to_zero}, it follows that $(d,\umbf) = (0, \umbf^*)$ is stationary for the feasibility subproblem 
at $x^*$. Therefore, it follows from \eqref{eq.zeta} and the fact under Assumption~\ref{global} that $d=0$ lies in the interior of $X$ that $u_{m+1}^* = 0$ and
\[
  \begin{aligned}
    u_i^* &= \begin{cases}  \zeta_*^i \frac{\nabla c_i(x^*)}{\|\nabla c_i(x^*)\|_2}\ \ \text{with}\ \  \zeta^i_*\in[-1,1] & \text{for all } i\in\Ecal\\
    \zeta_*^i \frac{\nabla c_i(x^*)}{\|\nabla c_i(x^*)\|_2}\ \ \text{with}\ \ \zeta^i_*\in[0,1]  &  \text{for all } i\in\Ical \end{cases} \\ \text{meaning that}\ \ 
 \delta^*(u_i^*| C_i^*) &=  - \zeta^i_*  \tfrac{c_i(x^*)}{\|\nabla c_i(x^*)\|_2}\ \  \text{for all } i\in \Ecal\cup\Ical.
\end{aligned} 
\]
 It follows that 
 \[\begin{aligned}
0 & = v(x^*) = J(0, 0; x^*,\cdot )= D(\umbf^*, 0; x^*, \cdot ) \\
& = -\tfrac{1}{2}\langle u_0^*,  H(0; x^*, \cdot )^{-1}u^*_0 \rangle -\sum_{i\in\Ecal\cup\Ical} \|\nabla c_i(x^*)\|_2\delta^*(u_i^*|C_i) -\delta^*(u^*_{m+1}|X)\\
& = -\tfrac{1}{2}\langle u_0^*,  H(0; x^*, \cdot )^{-1}u^*_0 \rangle + \sum_{i\in\Ecal\cup\Ical} \zeta^i_*c_i(x^*)\\
& = -\tfrac{1}{2}\langle u_0^*,  H(0; x^*, \cdot )^{-1}u^*_0 \rangle + \sum_{i\in\Ncal(x^*)} \zeta^i_*c_i(x^*). 
\end{aligned}.\]
Since $H(0; x^*, \cdot )$ is positive definite and $\sum_{i\in\Ncal(x^*)} \zeta^i_*c_i(x^*)\le 0$, it follows that
\[ \tfrac{1}{2}\langle u_0^*,  H(0; x^*, \cdot )^{-1}u^*_0 \rangle = 0\quad\text{and}\quad \sum_{i\in\Ncal(x^*)} \zeta^i_*c_i(x^*)=0,\] 
 yielding $u_0^*=0$ and $\zeta_*^i=0$ for all $i\in\Ncal(x^*)$.
  Overall, we have shown that the constraints of \eqref{prob.D} imply that
\begin{equation}\label{sumis0lam} 
\sum_{i\in\Ecal\cup\Acal(x^*)} \zeta_*^i \nabla c_i(x^*) = 0.
\end{equation}
Therefore,  $x^*$ violates the dual form of the MFCQ  because   $\zeta^i_*, i\in\Ecal\cup\Acal(x^*)$ are not all zero. 
%\eqref{sumis0lam} implies that $x^*$ violates the dual form of the MFCQ \cite{Solodov10}. 
Since we have reached a contradiction, it follows that the MFCQ cannot hold at $x^*$, as desired.
\end{proof}

%\blue{I believe the proof is complete here since this violates the dual form of the MFCQ.
%That is, \eqref{sumis0lam} implies that $x^*$ violates the dual form of the MFCQ.}
%\medskip
%\noindent
% In order to reach a contradiction, now assume that the MFCQ holds at $x^*$.  Then, there exists a vector $z\in\mathbb{R}^n$ such that 
%\begin{equation}\label{mfcq0}
% \langle z, \nabla c_i(x^*)\rangle =0\ \text{for}\ i\in\Ecal\ \text{and}\ \langle z, \nabla c_i(x^*)\rangle < 0\ \text{for}\ i\in \Acal(x^*).
% \end{equation}
%From \eqref{sumis0lam}, this means that
%\begin{equation}
%\label{0AE}
% 0 = \sum_{i\in\Ecal}\zeta_*^i \langle z, \nabla c_i(x^*)\rangle  + \sum_{i\in\Acal(x^*)} \zeta_*^i \langle z, \nabla c_i(x^*)\rangle
% = \sum_{i\in\Acal(x^*)}\zeta_*^i \langle z, \nabla c_i(x^*)\rangle.\end{equation}
%  Combining \eqref{mfcq0} and \eqref{0AE}, one finds that $\zeta^i_* = 0$ for $i\in\Acal(x^*)$,   so \eqref{sumis0lam} implies
%\[ \sum_{i\in\Ecal} \zeta_*^i \nabla c_i(x^*) = 0.\]
% Since the MFCQ holds at $x^*$, meaning that the vectors $\nabla c_i(x^*)$ for $i\in\Ecal$ are linearly independent, one can conclude that $\zeta_*^i = 0$ for $i \in \Ecal$.  However, since we have shown that $\zeta_*^i$ for $i \in \Ecal \cup \Ical$ are all zero, we have reached a contradiction---due to \eqref{eq.zeta}---to the fact that $\umbf^*$ is nonzero.  Since we have reached a contradiction, it follows that the MFCQ cannot hold at $x^*$, as desired.

We summarize the results of all of our theorems in the following corollary.

\begin{corollary}  Suppose Assumption~\ref{global} holds. Then, one of the following occurs.
\begin{enumerate}
\item[(i)]  $\rho_k \to \rho_*$ for some constant $\rho_*>0$ and each limit point of $\{x^k\}$ either corresponds 
to a KKT point or an infeasible stationary point for problem~\eqref{prob.nlp}.
\item[(ii)] $\rho_k \to 0$ and all limit points of $\{x^k\}$ are infeasible stationary points for \eqref{prob.nlp}.
\item[(iii)]  $\rho_k \to 0$, all limit points of $\{x^k\}$ are feasible for \eqref{prob.nlp}, and the MFCQ fails to hold at all limit points 
of $\{x^k\}_{k\in\Dcal \cup\Tcal}$. 
 
\end{enumerate}
 
\end{corollary}

%*********
% Section
%*********
\section{Implementation}\label{sec.implementation}

In this section, we discuss techniques that can be used for implementing our method.  In \S\ref{sec.low_rank}, we describe details about how L-BFGS Hessian approximations could be updated. In \S\ref{sec.solvers}, we introduce a coordinate descent method as an example subproblem solver that could be used with our method.

%************
% Subsection
%************
\subsection{Discussion on L-BFGS Hessian approximation}\label{sec.low_rank}

In large-scale settings, it is often intractable to compute and store exact Hessians. Instead, limited-memory approximations of the Hessian could be used, e.g., based on L-BFGS \cite{byrd1994representations, Noce80}. In this section, we describe how to update the    Hessian approximation and its inverse when $\rho$ is updated by \eqref{dust}.

Assume the Hessian approximations have the form
\bequationn
  H_\rho  = \sig  I+\Psi \Sigma^{-1}\Psi^T\ \ \text{and}\ \ H_0  = \gam I+\Phi \Gamma^{-1}\Phi^T,
\eequationn
where $\Psi \in \R^{n\times r}$ with $r \ll n$ and $\Phi \in \R^{n\times l }$ with $l \ll n$ are low rank matrices, and $\Sigma\in\R^{r \times r}$ and $\Gamma\in\R^{l\times l}$ are invertible. We investigate the inverse of $H_\rho$ by using the following generalized matrix inversion formula. For any given invertible $A\in\bR^{n\times n}$, invertible $S\in\bR^{l\times l}$, and $U, V\in \bR^{n\times l}$, the Sherman-Morrison formula yields
%\begin{equation}\label{smwk}
%  (A+USV^T)^{-k} = A^{-k}-A^{-k}U(S^{-1}+V^TA^{-k}U)^{-1}WV^TA^{-k}\, ,
%\end{equation}
%where
%$$
%  W=\sum^{k-1}_{j=0} (S^{-1}(S^{-1}+V^TA^{-k}U)^{-1})^{(j)}\, .
%$$
%In particular, 
\begin{equation}\label{smwk1}
  (A+USV^T)^{-1}= A^{-1}-A^{-1}U(S^{-1}+V^TA^{-1}U)^{-1}V^TA^{-1}.
\end{equation}
%Let 
%\[
%  H_{\rho,\mu} = H_\rho + \frac{m}{\mu} I= (\sigma + \frac{m}{\mu}) I + \Psi \Sigma^{-1}\Psi^T.
%\]
Using \eqref{smwk1},   the inverses of $H_0$ and $H_\rho$   are given by  
\begin{align*}
  H_0^{-1} = \frac{1}{\gamma}\left[ I - \Phi(\gamma\Gamma+\Phi^T\Phi)^{-1}\Phi^T \right]\ \ \text{and}\ \ H_\rho^{-1} = \frac{1}{\sigma}\left[ I - \Psi(\sigma \Sigma + \Psi^T\Psi)^{-1} \Psi^T  \right].
 % H_{\rho,\mu}^{-1} & =  \frac{1}{\sigma+\frac{m}{\mu}}\left[ I - \Psi[(\sigma +\frac{m}{\mu})\Sigma + \Psi^T\Psi]^{-1} \Psi^T  \right].
\end{align*}
These can be rewritten in a compact form.  Defining
\begin{align*}
  \Theta^T_1 = (\gamma \Gamma + \Phi^T \Phi)^{-1} \Phi^T,\ \ \Theta^T_2 &= (\sigma \Sigma + \Psi^T\Psi)^{-1} \Psi^T,
 % \text{and}\ \ \Theta^T_2 & = \left[\(\sigma +\frac{m}{\mu}\)\Sigma + \Psi^T\Psi\right]^{-1} \Psi^T,
\end{align*}
it follows that
\begin{equation}\label{theta.inverse}
  H_0^{-1}  =  \frac{1}{\gamma} (I - \Phi \Theta_1^T)\quad\text{and}\quad
  H_\rho^{-1}       = \frac{1}{\sigma}\left[ I - \Psi \Theta_2^T \right].
 % H_{\rho,\mu}^{-1} & =  \frac{1}{\sigma+\frac{m}{\mu}}\left[ I - \Psi \Theta_2^T  \right].
\end{equation}

After reducing $\rho$ to a smaller value $\bar \rho < \rho$, one finds that
\begin{align*}
  H_{\bar \rho} =  \bar\rho H_f + H_0 &= \frac{\bar \rho}{\rho} ( H_0+\rho H_f) + (1-\frac{\bar\rho}{\rho}) H_0 \\
  & =\tau H_\rho + (1-\tau) H_0 \\
  & = \bar\sigma I + \tau \Psi\Sigma^{-1}\Psi^T + (1-\tau) \Phi \Gamma^{-1} \Phi^T\\
  & = H_\tau + (1-\tau)\Phi\Gamma^{-1}\Phi^T,
\end{align*}
with 
\[
  \tau = \frac{\bar \rho}{\rho},\quad \bar\sigma =\tau\sigma+(1-\tau)\gamma,\quad\text{and}\quad   H_{\tau} =\bar\sigma I + \tau \Psi\Sigma^{-1}\Psi^T  .
\] 
Therefore, we have
\[
\begin{aligned}
H_{\tau}^{-1} & =  \frac{1}{\bar \sigma} [ I - \Psi \Theta_3^T ]&& \text{with}\quad   \Theta_3^T = ( \frac{\bar\sigma}{\tau}\Sigma + \Psi^T\Psi)^{-1}\Psi^T,\\ \text{and}\ \ 
H_{\bar\rho}^{-1} & = H_{\tau}^{-1} - H_{\tau}^{-1} \Phi \Theta_4^T H_{\tau}^{-1}&&\text{with}\quad   \Theta_4^T  = \left[\frac{1}{1-\tau}\Gamma + \Phi^T H_{\tau}^{-1} \Phi  \right]^{-1} \Phi^T.
\end{aligned}
\]
%Substituting $\sigma$ with $\sigma+\frac{m}{\mu}$,  and repeating \eqref{inverse.newrho}, we have the inverse of $H_{\bar\rho,\mu}$:
%\begin{align}\label{inverse.newrhomu}
%  \tilde\sigma & =\tau\sigma+\frac{m}{\mu}+(1-\tau)\gamma, \   H_{\tilde\tau} =\tilde\sigma I + \tau \Psi\Sigma^{-1}\Psi^T\\
%  \Theta_5^T & = ( \frac{\bar\sigma}{\tau}\Sigma + \Psi^T\Psi)^{-1}\Psi^T\\
%  H_{\tilde\tau}^{-1} & =  \frac{1}{\bar \sigma} [ I - \Psi \Theta_3^T ] \\ 
%  \Theta_6^T & = \left[\frac{1}{1-\tau}\Gamma + \Phi^T H_{\tau}^{-1} \Phi  \right]^{-1} \Phi^T\\
%  H_{\bar\rho,\mu} & = H_{\tilde\tau}^{-1} - H_{\tilde\tau}^{-1} \Phi \Theta_4^T H_{\tilde\tau}^{-1}.
%\end{align}

%************
% Subsection
%************

\subsection{Subproblem Solver}\label{sec.solvers}

As an example of a subproblem solver that can be used within our approach, we present a coordinate descent algorithm to solve \eqref{prob.J}.  For simplicity, let us assume that $X=\bR^n$. 
%with properties satisfying Assumption~\ref{ass.algorithm} so that \ref{dust} and \ref{psst} can be incorporated in the proposed algorithms.   
We have the following two subproblems:
\begin{align}
  \label{opt} \min_{x \in \bR^n} \ & J(x; \rho) := \half x^T H_\rho x +\rho g^T x + \sum_{i = 1}^{\mbar} |a_i^T x + b_i| + \sum_{i = \mbar+1}^m (a_i^T x + b_i)_+ \\
  \label{fea} \text{and}\ \  \min_{z \in \bR^n} \ & J(z; 0) :=   \half z^T H_0 z +  \sum_{i = 1}^{\mbar} |a_i^T z + b_i| + \sum_{i = \mbar+1}^m (a_i^T z + b_i)_+ .
\end{align}
Lagrangian duals of \eqref{opt} and \eqref{fea} are, respectively,
\begin{align}
  \label{dual_opt} \max_{l \le \eta \le c} \ & D(\eta; \rho) : = - \frac{1}{2} (A^T\eta - \rho g)^T H_\rho^{-1} (A^T\eta - \rho g) + \eta^Tb\\
  \label{dual_fea} \text{and}\ \ \max_{l \le \lambda \le c} \ & D(\lambda; 0) : = -\frac{1}{2} \lambda^T A H_0^{-1} A^T \lambda + \lambda^T b
\end{align}
where $l = [\pmb{-1}_{\mbar}^T, \pmb{0}_{m - \mbar}^T]^T$, $A=[a_1, \ldots, a_m]^T$ and $c = \pmb{1}_m$. The solutions of \eqref{opt} and \eqref{fea} can be recovered by those of \eqref{dual_opt} and \eqref{dual_fea} as $x = - H_\rho^{-1}(\rho g + A^T \eta)$ and $z = - H_0^{-1}A^T\lambda$, respectively. If we solve the feasibility dual problem \eqref{dual_fea}, this will give us a better estimate of $r_v$ at the extra cost of solving for $\lambda$. If this cost becomes prohibitive, we can use $\eta$ instead of $\lambda$ in the calculation of $r_v$. This might lead to more iterations for the subproblem solver. Algorithm \ref{cord} shows one iteration update of a coordinate descent algorithm. Note that subproblems \eqref{sub_opt} and \eqref{sub_fea} minimize one dimensional quadratics over a box constraint; hence, these have closed form solutions. \begin{algorithm}
  \caption{Coordinate Descent Algorithm} 
  \label{cord}
  \begin{algorithmic}[1]
    %\State  (Update $\eta^{k}$ and $\lambda^{k}$)
    \For{$i = 1, \dots, m$} 
      \State Set
      \begin{align}
        \label{sub_opt} \eta_i^{k} & := \argmax_{l_i \le \eta_i \le c_i} D(\eta_1^{k}, \dots, \eta_{i-1}^{k}, \eta_i, \eta_{i+1}^{k-1}, \dots, \eta_m^{k-1}; \rho^{k-1}) \\
        \label{sub_fea} \text{and}\ \ \lambda_i^{k} & := \argmax_{l_i \le \lambda_i \le c_i} D(\lambda_1^{k}, \dots, \lambda_{i-1}^{k}, \lambda_i, \lambda_{i+1}^{k-1}, \dots, \lambda_m^{k-1}; 0).
      \end{align}
      \State Update $x^k := -H_\rho^{-1}(\rho g + A^T \eta^k)$.
      \State Set $\rho_{k}$ by applying \eqref{dust}.
    \EndFor
  \end{algorithmic}
\end{algorithm}

We will now discuss how to make one sweep over all coordinates in an efficient manner when we use   Hessian approximations. Since \eqref{sub_opt} and \eqref{sub_fea} have similar structure, we will use  \eqref{sub_fea} to demonstrate the implementation details.

Using \eqref{theta.inverse}, subproblem \eqref{dual_fea} can be written as 
\begin{align}\label{cord_dual}
  \max_{l \le \lambda \le c} \ & D(\lambda; 0) : = -\frac{1}{2\gamma} \lambda^T A A^T \lambda + \frac{1}{2\gamma} \lambda^T A\Phi\Theta_1^T A^T \lambda + \lambda^T b.
\end{align}
In large scale settings, it is not practical to calculate and store $AA^T$. Usually, $A$ will have a nice sparse structure, while $AA^T$ does not. Defining $Q := A\Phi$ and $\tilde Q := A\Theta_1$, subproblem \eqref{cord_dual} becomes
\begin{align}
  \max_{l \le \lambda \le c} \ & D(\lambda; 0) : = - \frac{1}{2\gamma} \lambda^T A A^T \lambda +  \frac{1}{2\gamma} \lambda^T Q \tilde Q^T \lambda + \lambda^T b.
\end{align}
The partial derivative of $D(\lambda; 0)$ with respect to $\lambda_i$ is given by
\begin{align}\label{dd}
  \frac{\partial D(\lambda; 0)}{\partial \lambda_i} := \frac{1}{\gamma} \sum_{j = 1}^m (- a_i^T a_j + q_i \tilde q_j^T )\lambda_j + b_i, 
\end{align}
where $q_i$ and $\tilde q_i$ are the $i$-th row of $Q$ and $\tilde Q$ respectively. 
Then, 
%the solution of 
\eqref{sub_fea} becomes
%\begin{align}
%  \lambda_i^k  & = \begin{cases}
%  & l_i \quad \mbox{if }  a_i^T a_i - q_i \tilde q_j^T = 0 \mbox{ and } \frac{\partial D(\lambda; 0)}{\partial \lambda_i} < 0 
%  \\
%  & [l_i, c_i] \quad \mbox{if } a_i^T a_i - q_i \tilde q_j^T = 0 \mbox{ and } \frac{\partial D(\lambda; 0)}{\partial \lambda_i} =  0 
%  \\ 
%  & c_i \quad \mbox{if } a_i^T a_i - q_i \tilde q_j^T = 0 \mbox{ and } \frac{\partial D(\lambda; 0)}{\partial \lambda_i} >  0 
%  \\
%  &\mbox{mid}\{\frac{\gamma b_i - \sum\limits_{j = 1}^{i-1}(a_i^T a_j - q_i \tilde q_j^T) \lambda_j^{k} -  \sum\limits_{j = i+1}^{n}(a_i^T a_j - q_i \tilde q_j^T) \lambda_j^{k-1}  }{a_i^T a_i - q_i \tilde q_j^T}, l_i, c_i\}  \ \mbox{if } a_i^T a_i - q_i \tilde q_j^T \neq 0. 
%  \end{cases}
%\end{align}
\begin{align}
  \lambda_i^k  & = \begin{cases}
   l_i & \mbox{if }  a_i^T a_i - q_i \tilde q_j^T = 0 \mbox{ and } \frac{\partial D(\lambda; 0)}{\partial \lambda_i} < 0 
  \\
   [l_i, c_i] & \mbox{if } a_i^T a_i - q_i \tilde q_j^T = 0 \mbox{ and } \frac{\partial D(\lambda; 0)}{\partial \lambda_i} =  0 
  \\ 
   c_i & \mbox{if } a_i^T a_i - q_i \tilde q_j^T = 0 \mbox{ and } \frac{\partial D(\lambda; 0)}{\partial \lambda_i} >  0 
  \\
  \mu_i 
  & \mbox{if } a_i^T a_i - q_i \tilde q_j^T \neq 0, 
  \end{cases}
\end{align}
where 
\[
\mu_i:=\mbox{mid}\left\{\frac{\gamma b_i - \sum\limits_{j = 1}^{i-1}(a_i^T a_j - q_i \tilde q_j^T) \lambda_j^{k} -  \sum\limits_{j = i+1}^{n}(a_i^T a_j - q_i \tilde q_j^T) \lambda_j^{k-1}  }{a_i^T a_i - q_i \tilde q_j^T},\ l_i,\ c_i\right\}.
\] 
%Now we see that 
Hence, the main calculation for the solution of \eqref{sub_fea} is the partial derivative \eqref{dd}. Direct computation of \eqref{dd} takes $O(mn + ml)$ operations. 
In \cite{rdcord}, it is shown that coordinate descent will become competitive if there is an efficient way to compute the partial derivative. Here, if we keep track of the vectors $v:= \sum_{j = 1}^m \lambda_j a_j$ and $p:= \sum_{j = 1}^m \lambda_j \tilde q_j$, then the complexity of the update of the derivative becomes $O(n + l)$ which is much better than $O(mn + ml)$. First notice that if we have $v$ and $p$ for 
the most recent $\lambda$, then 
\[
  \frac{\partial D(\lambda; 0)}{\partial \lambda_i}  = \frac{1}{\gamma} (-a_i^T v + q_i p^T) + b_i,
\]
i.e. given $v$ and $p$, calculating \eqref{dd} takes only $O(n + l)$ operations. Next let us see how to update $v$ and $p$. Assume we update $\lambda_i^{k-1}$ to $\lambda_i^k$, then
\begin{equation*}
  v  \gets v + (\lambda_i^k - \lambda_i^{k-1})a_i\quad\text{and}\quad
  p \gets p + (\lambda_i^k - \lambda_i^{k-1})\tilde q_i. 
\end{equation*}
This shows the update of $v$ and $p$ is $O(n + l)$. In summary, the total complexity for each coordinate update is $O(n + l)$. Moreover, if $A$ is a sparse matrix with an average of $n_s$ nonzeros per row, then the complexity becomes $O(n_s+ l)$.

%*********
% Section
%*********
\section{Numerical Experiments}\label{sec.numerical}
In this section, we present our experimental results on both feasible and infeasible test sets. For each iteration,  
Algorithm~\ref{cord} described in \S\ref{sec.solvers} is used to solve for dual variables, which in turn is used
to obtain the corresponding primal variables. Our code is implemented using Python and tested on a 2014 MacBook Air with 4 GB memory and 1.4 GHz Intel Core i5.   %, and our proposed penalty parameter updating strategy is  incorporated.

\subsection{Feasible test} \label{hs.feasible.test}We tested on 126 CUTEr Hock-Schittkowski (\texttt{hs}) problems \cite{Hock1980} which are all feasible.  We set the parameters stated in Algorithm \ref{alg.sqo.full} as $\gamma = 0.5$, $\rho_{(-1)} = 1$, $\beta_\phi = 0.7$, $\beta_v = 0.1$, $\beta_{l} = 0.6 \beta_\phi(1-\beta_v)$, $\omega_{0} = 10^{-2}$, $\theta_\rho = 0.9$, $\theta_\omega = 0.7$, $\theta_\alpha = 10^{-4}$, and $\eta^{0} = \pmb{0}_m$ with $x^{0}$ set as defined for each CUTEr problem. The maximum iteration limit for the subproblem solver was set as $10^6$, while a maximum iteration limit for Algorithm \ref{alg.sqo.full} was set to be $200$. We defined the maximum constraint violation $v_\infty(x)$ and the optimality KKT error $\epsilon_{opt}(x)$ as 
\begin{align*}
v_\infty(x) & :=  \max\{|c_i(x)| \ i = 1, \cdots, \bar{m},\  (c_i(x))_{+} \ i = \bar{m}+1, \cdots, m\}, \\
\epsilon_{opt}(x) & := \max{\left\{\norm{\grad f(x) + \sum_{i = 1}^m \eta_i \grad c_i(x) }_\infty, \norm{\eta \circ c(x)}_\infty \right\}},
\end{align*}
where $\circ$ denotes element-wise product. We terminate the algorithm if $ v_\infty(x) \le 10^{-5}$ and  $\epsilon_{opt}(x) \le 10^{-4}$, or the maximum iteration number $200$ is reached. These 126 problems are of small size; hence we use the exact Hessian in our implementation.  If the Hessian, call it $H$, is not positive definite, then we apply the following modification to adjust its negative eigenvalues.  Let $H = U \Lambda U^T$ be the eigen-decomposition of $H$, where $\Lambda = \diag \{\lambda_1, \cdots, \lambda_n\}$.  For a prescribed constant $\tau > 0$ (e.g., we use $\tau = 10^{-4}$ in these experiments), we reset $\lambda_i \gets \max\{\lambda_i,\tau\}$ and replace $H$ with $U\tilde\Lambda U^T$ where $\tilde\Lambda$ is the corresponding modification of $\Lambda$.  We also perform the following modification to control the condition number of the Hessian (approximation) employed in the algorithm: If $\text{cond}(H) > t_{c} > 0$ (e.g., we use $t_c = 10^{6}$ in these experiments), then we replace $H$ by $\alpha H + (1-\alpha)I$ where $\alpha$ is the largest value in $[0,1]$ such that the resulting matrix has condition number less than or equal to $t_c$.
%\textcolor{blue}{If the Hessian, say $H$, is not positive definite, then we apply the following modification(s) to adjust its negative eigenvalues and condition number.  Let $H = U \Lambda U^T$ be the eigen-decomposition of $H$, where $\Lambda = \diag \{\lambda_1, \cdots, \lambda_n\}$.  We replace all $\lambda_i < 0$ with a prescribed small constant $\tau = 10^{-4}$. Let $\tilde H$ be the modified positive definite Hessian.  If the condition number $\text{cond}(\tilde H)$ of the modified Hessian is greater than $t_{c} = 10^{6}$, then we use in the algorithm the modified matrix $H_{mod} = (t_c/\text{cond}(\tilde H)) \tilde H + \tau I.$ Otherwise, we use $H_{mod} = \tilde H.$}

\begin{table}[ht]
\centering
\caption{Performance comparison  of SQuID and the proposed algorithm on feasible problems.}
\label{tab.squid.fea}
\begin{tabular}{|c|c|cccc|}\hline
Problem type & Algorithm & Succeed & Fail & Infeasible & Total \\ \hline
\multirow{2}{5em}{Feasible \texttt{hs} problems} & SQuID  & 110 (90.16\%) & 11 (9.02\%) & 1 (0.82\%)  &  122\\ \cline{2-6}
&Proposed  & 115 (91.20\%) & 11 (8.80\%)  & 0 & 126 \\ \hline
\end{tabular}
\end{table}

\begin{figure}[t]
 \begin{center}
  \includegraphics[scale = 0.45]{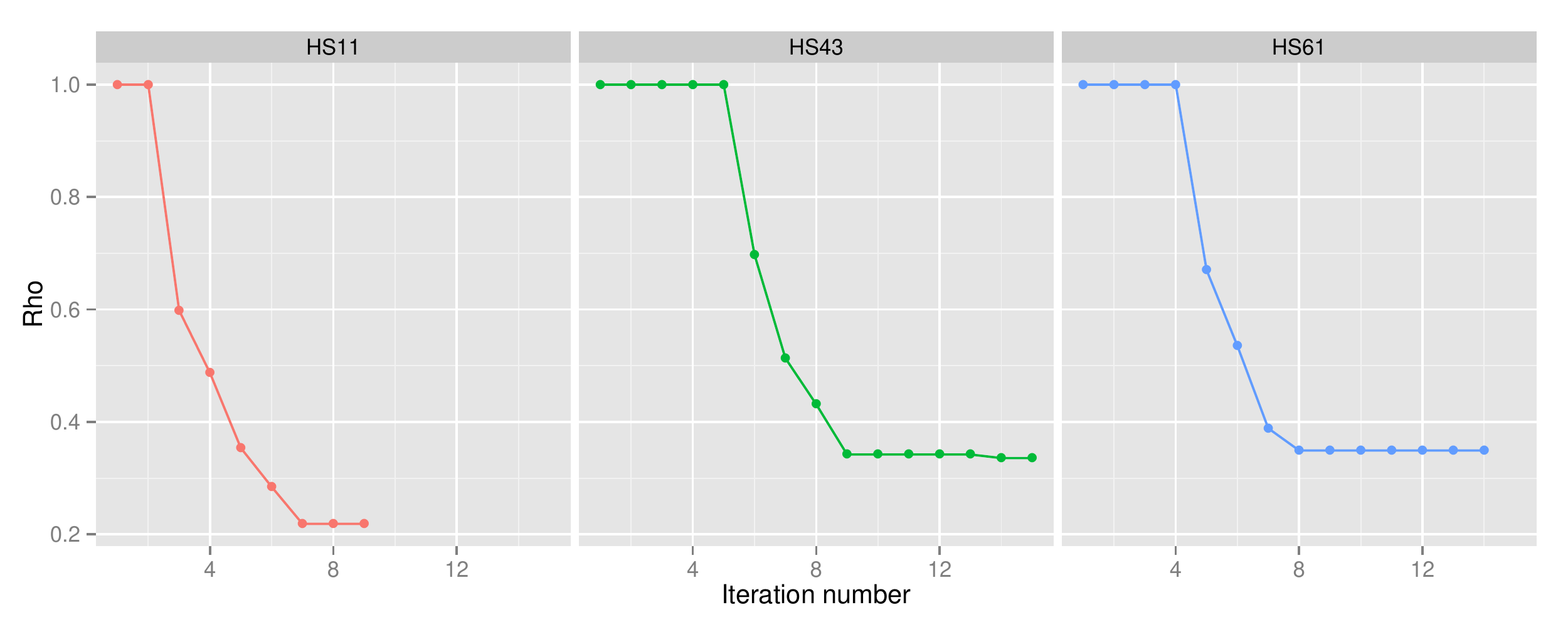}
 \end{center}
 \caption{$\rho$ values for problems \texttt{hs11}, \texttt{hs43} and \texttt{hs61}.}
 \label{fig:rho-update}
\end{figure}

For these experiments, we have the following observations.
\begin{itemize}
\item  
Out of 126 CUTEr \texttt{hs} problems, our algorithm successfully solved $115$, which is a success rate of about $91 \%  \approx 115/126.$   
\footnote{The termination criterion of SQuID in \cite{BurkCurtWang14} is based on the relative KKT residual scaled by $\rho$.}
Our proposed method outperforms the SQuID algorithm proposed in \cite{BurkCurtWang14}, which is also a penalty-SQP method with automatic infeasibility detection, although it requires two exact QP solves per iteration. The comparison statistics\footnote{The performance statistics for SQuID is obtained from \cite{BurkCurtWang14}, where the overall number of \texttt{hs} problems is 122 due to compiling errors.} are shown in Table~\ref{tab.squid.fea}.
%For the complete set of numerical test results see \cite{burke2018dynamic}.
Table \ref{tab:CUTEr-test-0}  summarizes the detailed output for these $115$ successful cases, where ``\# iter'' means the number of iterations and ``\# $f$'' denotes the number of function evaluations. 
\item   Our \eqref{dust} updating strategy works very well in these experiments, and does not cause $\rho$ to become excessively small for most cases.  To illustrate the behavior of the penalty parameter updates, we plot $\rho$ values for three sample problems---\texttt{hs11}, \texttt{hs43} and \texttt{hs61}---in Figure \ref{fig:rho-update}.
\item  The parameter $\omega$ did not require much tuning. We used $\omega_0 = 10^{-2}$ across all problems and achieved our $91\%$ success rate. We also ran the experiment with $\omega_0 = 10^{-1}$ and saw the same set of 115 problems solved successfully.
\item We test the sensitivity of our algorithm with respect to the parameter $\beta_\phi$. We ran the same experiments with $\beta_\phi = 0.5$ and $\beta_\phi = 0.99$. We have 113 successful cases for $\beta_\phi = 0.5$,
(see Table~\ref{tab:CUTEr-test-beta-phi-p5}), 
and 111 successful cases for $\beta_\phi = 0.99$
(see Table~\ref{tab:CUTEr-test-beta-phi-p99}). 
The additional failure cases in $\beta_\phi = 0.5$ and $\beta_\phi = 0.99$ compared to $\beta_\phi = 0.7$ are all due to subproblem exceeding the maximum iteration number.
\item   
Coordinate descent performs poorly on ill-conditioned subproblems. We observed that some subproblems require more than $5\times10^5$ steps to reach the specified accuracy. Since the focus of this paper is on the $\rho$ update strategy, we did not 
explore other subproblem solvers that might have performed better. 
Instead, we used a large iteration limit for the subproblem solver. 
%In the future work, it is worthwhile to explore other robust subproblem solvers.
\item   
In a few cases, the Hessian modification strategy described above did not work well. 
For example, for problems \texttt{hs72} and \texttt{hs75}, 
we had to reduce the modification constant to $10^{-8}$ to achieve convergence, since the scale of the Hessian for both problems is around $10^{-4}$. For problem \texttt{hs93}, convergence is observed with modification constant $10^{-2}$.
%We notice that the size of the positive eigenvalues of these problems are around $10^{-4}$, so our modification significantly changed the structure of the Hessian preventing convergence within the iteration limit.
\end{itemize}

\subsection{Infeasible test} As in \cite{BurkCurtWang14}, we modified the 126 CUTEr Hock-Schittkowski (\texttt{hs}) problems by adding bound constraints $x_1 \le 0$ and $x_1 \ge 1$ to make all \texttt{hs} problems infeasible; we refer to these problems as \texttt{hs\_inf}. All the parameters used for this infeasible test set are the same as mentioned for the feasible test set, except we increase the maximum iteration limit for the subproblem solver to 20000. 
 Defining the feasibility KKT error $\epsilon_{fea}(x)$ as 
\begin{align*}
\epsilon_{fea}(x) := \max \left\{  \norm{ \sum_{i = 1}^m \eta_i \grad c_i(x) }_\infty,\right. & \norm{(e-\eta^\Ecal)\circ [c_\Ecal(x)]^+}_\infty,  \norm{(e+\eta^\Ecal)\circ [c_\Ecal(x)]^-}_\infty,  \\
 &\left. \norm{(e-\eta^\Ical)\circ [c_\Ical(x)]^+}_\infty,  \norm{\eta^\Ical \circ [c_\Ical(x)]^-}_\infty \right\},
\end{align*}
we use the same stopping criteria as in \cite{BurkCurtWang14}, except that we do not necessarily need to drive $\rho$ to 0; hence we drop ``$\rho \le 10^{-8}$'' from the stopping criteria used in \cite{BurkCurtWang14}.

\begin{table}[ht]
\centering
\caption{Performance comparison  of SQuID and the proposed algorithm on infeasible problems.}
\label{tab.squid.inf}
\begin{tabular}{|c|c|ccc|}\hline
Problem type & Algorithm & Succeed  & Fail  & Total \\ \hline
\multirow{2}{9em}{Infeasible \texttt{hs} problems (\texttt{hs\_inf})} & SQuID  & 111 (90.24\%) & 12 (9.76\%)  &  123  \\ \cline{2-5}
&Proposed  & 116 (92.10\%) & 10 (7.90\%)  & 126 \\ \hline
\end{tabular}
\end{table}

For these experiments, we have the following observations.
\begin{itemize}
\item  Out of 126   \texttt{hs\_inf} problems, our algorithm successfully solved $116$, which is a success rate of about $92 \%  \approx 116/126.$ Our proposed method also outperforms SQuID on infeasible problems.  The comparison statistics\footnote{The performance statistics for SQuID is obtained from \cite{BurkCurtWang14}, where the overall number of \texttt{hs} problems is 123 due to compiling errors.} are shown  in Table~\ref{tab.squid.inf}.  
Table~\ref{tab:CUTEr-inf-test} summarizes the detailed output for these $116$ successful cases. 
\item In a few cases, the Hessian modification strategy described above did not work well. For example, for problems  \texttt{hs104\_inf},  \texttt{hs114\_inf}, \texttt{hs8\_inf}, \texttt{hs23\_inf} and \texttt{hs93\_inf}, convergence is observed when we increase the modification constant from $10^{-4}$ to $10^{-2}$.
\end{itemize}

\subsection{Large scale test} We also applied our implementation to solve some large scale problems from the CUTEr test set; see Table~\ref{tab.large.scale.statitics}. The parameter settings used were the same as used in Section~\ref{hs.feasible.test}, except that we set the iteration limit for the subproblem solver to be 2000. We used L-BFGS for the Hessian approximations which 
pairs well with the coordinate descent algorithm giving a $O(n + \ell)$ total complexity for each coordinate update. 
%where $\ell$ is the memory length \cite{burke2018dynamic}. 
Table~\ref{tab.large.scale.test} presents the results for successful runs.  For the remaining problems not shown, 
the coordinate descent QP algorithm could not reach the desired accuracy within the maximum number of subproblem iterations.  
We leave further investigation into the most effective iterative QP solver for
these problems to future work since
this is beyond the scope of the paper. 
%Nonetheless,
%the experiments serve well to illustrate the effectiveness of the dynamic penalty parameter updating strategy.
%where beyond the scope of this paper a more powerful iterative QP solver could be employed.

To recognize the benefits of our proposed algorithm compared to an alternative approach, let us consider the CPU times required to run the experiments whose results are shown in Table~\ref{tab.large.scale.test} compared to the CPU times that would be required by SQuID from \cite{BurkCurtWang14}.  The aforementioned implementation of SQuID was not able to terminate successfully on any of the problems in Table~\ref{tab.large.scale.statitics} within 10 minutes.  The primary expense is solving the QP subproblems to high accuracy in each iteration.  By contrast, the result shown in Table~\ref{tab.large.scale.test} that required the most CPU time was the run for problem \texttt{LUKVLE10}, where the entire run terminated in 64 seconds.  The benefits of our proposed algorithm are clear when solving large-scale problems.  (On a contemporary laptop to computer, the state-of-the-art code Ipopt~\cite{WaecBieg06} solves problem \texttt{LUKVLE10} in only a couple of seconds, but that code benefits from two decades of software development.)

\begin{table}[h] \label{tab.large.scale.statitics}
\footnotesize
\centering
\caption{CUTEr 13 large scale problems.} 
\begin{tabular}{l|rrr}
\toprule
  Problem &    \# constraints  &   \# variables &    \# equalities    \\
\midrule
  \texttt{DTOC1NA} &   3996 &   5998 &  3996   \\
  \texttt{DTOC1NB} &   3996 &   5998 &  3996   \\
  \texttt{DTOC1ND} &   3996 &   5998 &  3996   \\
      \texttt{EG3} &  20000 &  10001 &     1        \\
% \texttt{EIGENACO} &   1275 &   2550 &  1275   \\
  \texttt{GILBERT} &      1 &   5000 &     1   \\
 \texttt{JANNSON4} &      2 &  10000 &     0   \\
  \texttt{LUKVLE1} &   9998 &  10000 &  9998   \\
 \texttt{LUKVLE10} &   9998 &  10000 &  9998   \\
  \texttt{LUKVLE3} &      2 &  10000 &     2   \\
  \texttt{LUKVLE6} &   4999 &   9999 &  4999   \\
 \texttt{LUKVLI13} &   6664 &   9998 &     0   \\
  \texttt{LUKVLI3} &      2 &  10000 &     0   \\
  \texttt{LUKVLI6} &   4999 &   9999 &     0   \\
\bottomrule
\end{tabular}
\end{table}

\begin{table}[h] \label{tab.large.scale.test}
\footnotesize
\caption{Test results on CUTEr 13 large scale problems.} 
\begin{tabular}{l|cccccc}
\toprule
  Problem &     \# iter &  \# $f$ &    $f(x^*)$ &     $v(x^*)$ &  KKT & Final $\rho$ \\
\midrule
  \texttt{DTOC1NA}    & 13  &        13 &  4.138866e+00 &  2.215482e$-$06 &  1.878763e$-$05 &      0.751447 \\
  \texttt{DTOC1NB}    & 13  &        13 &  7.138849e+00 &  4.835061e$-$07 &  3.550798e$-$05 &      0.849347 \\
  \texttt{DTOC1ND}    & 14  &        19 &  4.760303e+01 &  1.799990e$-$07 &  4.807062e$-$05 &      0.815373 \\
      \texttt{EG3}          & 10  &        10 &  8.048306e$-$06 &  0.000000e+00 &  7.317141e$-$05 &      0.479603 \\
% \texttt{EIGENACO}   &  28 &        75 &  1.467067e$-$15 &  0.000000e+00 &  2.907696e$-$08 &      0.375413 \\
  \texttt{GILBERT}      &  74 &        74 &  2.459468e+03 &  2.170219e$-$08 &  4.047202e$-$06 &      0.024360 \\
 \texttt{JANNSON4}   &  79 &        80 &  9.801970e+03 &  6.956902e$-$08 &  1.830188e$-$05 &      0.009923 \\
  \texttt{LUKVLE1}      &  13 &        25 &  4.821043e$-$14 &  3.087659e$-$08 &  5.364314e$-$05 &      0.960000 \\
 \texttt{LUKVLE10}     &191 &       191 &  3.534934e+03 &  2.224607e$-$09 &  9.783636e$-$05 &      0.282103 \\
  \texttt{LUKVLE3}      &  41 &        49 &  2.758658e+01 &  9.747758e$-$14 &  4.949755e$-$05 &      0.318856 \\
  \texttt{LUKVLE6}      &  39 &        68 &  6.286441e+05 &  1.436051e$-$12 &  6.916637e$-$05 &      0.360397 \\
 \texttt{LUKVLI13}       &  65 &        76 &  1.321855e+02 &  3.212148e$-$09 &  7.052314e$-$05 &      0.293858 \\
  \texttt{LUKVLI3}        &  70 &        78 &  1.157754e+01 &  9.010570e$-$13 &  6.644757e$-$05 &      0.442002 \\
  \texttt{LUKVLI6}        &  43 &        63 &  6.286441e+05 &  1.390753e$-$11 &  6.766861e$-$05 &      0.195366 \\
\bottomrule
\end{tabular}
\end{table}

\section{Conclusion}\label{sec.conclusion}

In this paper, we have proposed a penalty-SQP framework for solving nonlinear optimization problems.  The novelty of this work 
is a dynamic penalty parameter updating strategy that is carried out within the QP subproblem solver, so that at the end of the QP solve, 
a search direction and a new penalty parameter are both obtained.  
The key idea is to force improvement toward feasibility  
whenever optimality and complementarity are sufficiently improved.  
This enables the SQP algorithm to finish penalty parameter updating and infeasibility detection via  \emph{inexact} solves for only \emph{one} subproblem in each iteration, a feature which is not shared with most contemporary solvers which require
two subproblem solves per iteration.

The convergence properties that we have proved for our algorithm guarantees the effectiveness of our updating strategy under reasonable assumptions. %The SQP algorithm finds optimal solution when the MFCQ and converges to an infeasible stationary point when penalty parameter is driven to 0. 
The empirical effects of our strategy are demonstrated in numerical results on small CUTEr examples.  We remark, however, that the performance 
could be further enhanced 
%by developing 
with the development of a more efficient QP subproblem solver 
and a more robust approach to addressing ill-conditioning of the Hessian approximation.
%a carefully convexifying the Hessian to handle ill-condition or ill-possness. 
%Future work includes an examination of the success of the coordinate descent
%subproblem solver combined with the with low-rank approximation of the Hessian
%in the large-scale setting.
%For large-scale problems, it would be interesting to see the performance of our proposed algorithm and updating strategy combined with low-rank approximation of the Hessian, which we leave  a subject of future research. 

\begin{footnotesize}
\begin {longtable}{l||r|r||r|r|r|r}
 \caption {CUTEr \texttt{hs} test results, 115 successful cases out of 126   problems.} \label{tab:CUTEr-test-0} \\
 
\hline

\multicolumn{1}{c}{Problem}    & \multicolumn{1}{c}{\# iter}         &  \multicolumn{1}{c}{  \# $f$}  & 
 \multicolumn{1}{c}{   $f(x^*)$} &  \multicolumn{1}{c}{$v(x^*)$}  &  \multicolumn{1}{c}{ KKT  } &  \multicolumn{1}{c}{  Final $\rho$}   \\  \hline
\hline
\endfirsthead
\endhead

\multicolumn{7}{c}%
{{  \tablename\ \thetable{} -- continued from previous page}} \\
\hline
\multicolumn{1}{c}{Problem}  & \multicolumn{1}{c}{  \# iter}     &  \multicolumn{1}{c}{  \# $f$}  & 
 \multicolumn{1}{c}{   $f(x^*)$} &  \multicolumn{1}{c}{     $v(x^*)$}  &  \multicolumn{1}{c}{     KKT  } &    \multicolumn{1}{c}{      Final $\rho$} 
  \\  \hline
\hline
\endhead

\hline \multicolumn{7}{c}{{Continued on next page}} \\ \hline
\endfoot

\hline\hline
\endlastfoot

         \texttt{hs1} &      24 &        34 &  4.215353e$-$17 &  0.000000e+00 &  2.983621e$-$09 &  1.000000 \\
     \texttt{hs10} &       8 &         9 & $-$1.000001e+00 &  1.551523e$-$06 &  3.360370e$-$06 &  1.000000 \\
    \texttt{hs100} &      11 &        22 &  6.806301e+02 &  2.403271e$-$06 &  8.021402e$-$06 &  0.540664 \\
 \texttt{hs100lnp} &      10 &        31 &  6.806301e+02 &  9.444675e$-$06 &  1.615661e$-$06 &  0.150095 \\
 \texttt{hs100mod} &       7 &        23 &  6.786796e+02 &  3.188191e$-$09 &  4.598100e$-$09 &  0.572194 \\
    \texttt{hs101} &      51 &       145 &  1.809765e+03 &  3.225274e$-$08 &  4.579758e$-$07 &  0.000132 \\
    \texttt{hs102} &      43 &       107 &  9.118803e+02 &  1.106847e$-$07 &  1.286046e$-$06 &  0.000274 \\
    \texttt{hs103} &      45 &       135 &  5.436642e+02 &  3.902120e$-$06 &  7.284432e$-$06 &  0.000567 \\
    \texttt{hs104} &      13 &        47 &  4.200002e+00 &  1.768636e$-$06 &  4.698347e$-$05 &  0.900000 \\
    \texttt{hs105} &      20 &       990 &  1.044612e+03 &  1.067283e$-$07 &  8.477262e$-$05 &  0.004175 \\
    \texttt{hs107} &      22 &        26 &  5.054972e+03 &  8.384681e$-$06 &  1.727293e$-$05 &  0.000688 \\
    \texttt{hs108} &      18 &        32 & $-$6.749664e$-$01 &  1.208205e$-$07 &  5.516680e$-$05 &  0.900000 \\
    \texttt{hs109} &      78 &       220 &  5.362069e+03 &  7.846857e$-$07 &  8.667739e$-$05 &  0.128672 \\
     \texttt{hs11} &       8 &         9 & $-$8.498465e+00 &  1.225127e$-$07 &  1.349376e$-$06 &  0.218726 \\
    \texttt{hs110} &       3 &         5 & $-$4.577848e+01 &  0.000000e+00 &  6.067373e$-$09 &  1.000000 \\
    \texttt{hs111} &      21 &        33 & $-$4.776110e+01 &  4.440919e$-$06 &  2.258500e$-$05 &  0.042014 \\
 \texttt{hs111lnp} &      20 &        30 & $-$4.776119e+01 &  9.859016e$-$06 &  2.462468e$-$05 &  0.041440 \\
    \texttt{hs112} &      21 &        23 & $-$4.776117e+01 &  7.792329e$-$06 &  2.725259e$-$06 &  0.044630 \\
    \texttt{hs113} &      16 &        17 &  2.430625e+01 &  1.363357e$-$06 &  9.915806e$-$06 &  0.387420 \\
    \texttt{hs117} &      14 &        21 &  3.235087e+01 &  7.313662e$-$06 &  2.220430e$-$05 &  0.011973 \\
    \texttt{hs118} &      19 &        20 &  9.329922e+02 &  9.973521e$-$06 &  1.820177e$-$06 &  0.064666 \\
    \texttt{hs119} &      20 &        21 &  2.448993e+02 &  8.555172e$-$06 &  1.909610e$-$06 &  0.146587 \\
     \texttt{hs12} &       5 &         9 & $-$3.000000e+01 &  3.991066e$-$07 &  2.934243e$-$07 &  1.000000 \\
     \texttt{hs14} &      20 &        89 &  1.393453e+00 &  7.451143e$-$06 &  4.922505e$-$06 &  0.414328 \\
     \texttt{hs15} &      10 &        11 &  3.603797e+02 &  4.207535e$-$08 &  1.368772e$-$05 &  0.002465 \\
     \texttt{hs16} &      20 &        21 &  2.314403e+01 &  9.780175e$-$06 &  9.391250e$-$05 &  0.014781 \\
     \texttt{hs17} &       9 &        11 &  1.000000e+00 &  0.000000e+00 &  9.445528e$-$07 &  0.018248 \\
     \texttt{hs18} &       7 &         8 &  5.000000e+00 &  2.771074e$-$09 &  4.843401e$-$08 &  1.000000 \\
     \texttt{hs19} &      33 &        47 & $-$6.961824e+03 &  8.499568e$-$06 &  8.737095e$-$06 &  0.000508 \\
      \texttt{hs2} &       7 &         9 &  4.941229e+00 &  0.000000e+00 &  4.601831e$-$07 &  1.000000 \\
     \texttt{hs20} &      25 &        26 &  4.019818e+01 &  2.882271e$-$06 &  9.114286e$-$05 &  0.008031 \\
     \texttt{hs21} &       2 &         3 & $-$9.996000e+01 &  4.440892e$-$16 &  4.999500e$-$09 &  1.000000 \\
  \texttt{hs21mod} &       8 &         9 & $-$9.596000e+01 &  2.220446e$-$16 &  4.371328e$-$18 &  0.162885 \\
     \texttt{hs22} &      11 &       127 &  1.000002e+00 &  0.000000e+00 &  2.386921e$-$06 &  1.000000 \\
     \texttt{hs23} &      20 &        21 &  1.999981e+00 &  9.518632e$-$06 &  5.376723e$-$06 &  0.282430 \\
     \texttt{hs24} &      13 &        14 & $-$9.998651e$-$01 &  2.585931e$-$12 &  7.171626e$-$05 &  0.531441 \\
     \texttt{hs25} &       1 &         2 &  3.283500e+01 &  0.000000e+00 &  2.005805e$-$08 &  1.000000 \\
     \texttt{hs26} &      13 &        28 &  2.172765e$-$10 &  4.361497e$-$06 &  2.031413e$-$07 &  1.000000 \\
    \texttt{hs268} &       3 &         4 &  8.608487e$-$06 &  0.000000e+00 &  2.526260e$-$05 &  1.000000 \\
     \texttt{hs27} &       7 &        11 &  4.000000e$-$02 &  1.781389e$-$19 &  5.580537e$-$05 &  1.000000 \\
     \texttt{hs28} &       2 &         3 &  1.117108e$-$13 &  0.000000e+00 &  2.034355e$-$07 &  1.000000 \\
     \texttt{hs29} &       8 &         9 & $-$2.262742e+01 &  5.095551e$-$10 &  1.211382e$-$05 &  0.680483 \\
      \texttt{hs3} &       7 &         8 &  2.338799e$-$04 &  0.000000e+00 &  9.672226e$-$05 &  1.000000 \\
     \texttt{hs30} &      18 &        75 &  1.000103e+00 &  0.000000e+00 &  6.717072e$-$05 &  0.590490 \\
     \texttt{hs31} &       9 &        11 &  6.000000e+00 &  1.192496e$-$08 &  7.648903e$-$05 &  0.104597 \\
     \texttt{hs32} &      12 &        13 &  1.000880e+00 &  1.018075e$-$13 &  9.627754e$-$05 &  0.109419 \\
     \texttt{hs33} &       9 &        93 & $-$4.000000e+00 &  0.000000e+00 &  3.969005e$-$10 &  0.088629 \\
     \texttt{hs34} &      19 &        28 & $-$8.340328e$-$01 &  8.685515e$-$06 &  3.003686e$-$07 &  0.810000 \\
     \texttt{hs35} &       1 &         2 &  1.111111e$-$01 &  0.000000e+00 &  8.332954e$-$05 &  1.000000 \\
    \texttt{hs35i} &       1 &         2 &  1.111111e$-$01 &  0.000000e+00 &  8.332954e$-$05 &  1.000000 \\
  \texttt{hs35mod} &       2 &         3 &  2.500000e$-$01 &  1.110223e$-$16 &  3.332426e$-$05 &  1.000000 \\
     \texttt{hs36} &      16 &        17 & $-$3.299993e+03 &  0.000000e+00 &  9.244868e$-$05 &  0.003757 \\
     \texttt{hs37} &       8 &        11 & $-$3.456000e+03 &  2.426503e$-$12 &  3.681286e$-$05 &  0.006363 \\
     \texttt{hs38} &      37 &        57 &  4.607342e$-$11 &  0.000000e+00 &  4.104272e$-$06 &  1.000000 \\
     \texttt{hs39} &      22 &        23 & $-$1.000009e+00 &  9.375743e$-$06 &  6.443443e$-$06 &  0.687275 \\
   \texttt{hs3mod} &       2 &         3 &  8.026142e$-$14 &  0.000000e+00 &  1.162326e$-$07 &  1.000000 \\
      \texttt{hs4} &       2 &         3 &  2.666667e+00 &  0.000000e+00 &  3.149394e$-$15 &  0.228768 \\
     \texttt{hs40} &      19 &        20 & $-$2.500008e$-$01 &  9.306140e$-$06 &  2.419174e$-$06 &  1.000000 \\
     \texttt{hs41} &      12 &        85 &  1.925925e+00 &  8.483081e$-$06 &  1.021823e$-$06 &  0.440257 \\
     \texttt{hs42} &       4 &        18 &  1.385786e+01 &  2.177929e$-$10 &  1.588578e$-$05 &  0.313419 \\
     \texttt{hs43} &      15 &        17 & $-$4.399990e+01 &  1.315519e$-$08 &  3.235023e$-$05 &  0.324783 \\
     \texttt{hs44} &      18 &        19 & $-$1.499991e+01 &  5.129125e$-$06 &  1.160088e$-$05 &  0.079766 \\
  \texttt{hs44new} &      20 &        21 & $-$1.500002e+01 &  6.312695e$-$06 &  3.358322e$-$06 &  0.047101 \\
     \texttt{hs45} &      17 &        18 &  1.000018e+00 &  6.508614e$-$06 &  1.250044e$-$05 &  0.590490 \\
     \texttt{hs46} &      17 &        18 &  4.352054e$-$09 &  9.197854e$-$06 &  3.362897e$-$07 &  1.000000 \\
     \texttt{hs47} &      16 &        20 &  1.134167e$-$09 &  9.982054e$-$06 &  1.208686e$-$06 &  0.135085 \\
     \texttt{hs48} &       8 &         9 &  2.516051e$-$19 &  5.230553e$-$06 &  1.162349e$-$09 &  1.000000 \\
     \texttt{hs49} &      13 &        14 &  2.791394e$-$07 &  7.660539e$-$12 &  3.325774e$-$05 &  1.000000 \\
      \texttt{hs5} &       4 &         7 & $-$1.913223e+00 &  0.000000e+00 &  3.354020e$-$05 &  0.656100 \\
     \texttt{hs50} &      10 &        11 &  3.510103e$-$17 &  6.655753e$-$06 &  1.845457e$-$08 &  0.088629 \\
     \texttt{hs51} &       2 &         3 &  6.496671e$-$17 &  1.204511e$-$08 &  9.999249e$-$09 &  1.000000 \\
     \texttt{hs52} &      22 &        23 &  5.326608e+00 &  6.900875e$-$06 &  3.165318e$-$06 &  0.101755 \\
     \texttt{hs53} &      21 &        22 &  4.092979e+00 &  9.854539e$-$06 &  4.525700e$-$06 &  0.129100 \\
     \texttt{hs54} &       6 &         7 & $-$1.561253e$-$01 &  4.934009e$-$10 &  9.544783e$-$05 &  1.000000 \\
     \texttt{hs55} &      20 &        21 &  6.666664e+00 &  7.546464e$-$06 &  1.148850e$-$06 &  0.919332 \\
     \texttt{hs56} &       9 &        11 & $-$3.456002e+00 &  2.128966e$-$06 &  1.199327e$-$05 &  0.479441 \\
     \texttt{hs57} &       1 &         2 &  3.064627e$-$02 &  0.000000e+00 &  2.696159e$-$06 &  1.000000 \\
     \texttt{hs59} &       8 &        12 & $-$7.802789e+00 &  0.000000e+00 &  3.334643e$-$06 &  0.900000 \\
      \texttt{hs6} &       9 &        24 &  8.091820e$-$10 &  1.605262e$-$07 &  2.842575e$-$05 &  1.000000 \\
     \texttt{hs60} &       5 &         6 &  3.256820e$-$02 &  9.958889e$-$08 &  1.894968e$-$07 &  1.000000 \\
     \texttt{hs61} &      13 &        47 & $-$1.436461e+02 &  1.588448e$-$06 &  4.775740e$-$07 &  0.338698 \\
     \texttt{hs62} &       5 &         7 & $-$2.627251e+04 &  1.526557e$-$16 &  4.647296e$-$07 &  0.001456 \\
     \texttt{hs63} &      14 &        18 &  9.617152e+02 &  9.119443e$-$06 &  5.246424e$-$06 &  0.470499 \\
     \texttt{hs64} &      43 &        44 &  6.299843e+03 &  1.254110e$-$08 &  7.139929e$-$05 &  0.041838 \\
     \texttt{hs65} &       6 &         7 &  9.535284e$-$01 &  5.315868e$-$06 &  4.367160e$-$07 &  1.000000 \\
     \texttt{hs66} &       6 &        10 &  5.181619e$-$01 &  7.088741e$-$06 &  1.275973e$-$06 &  0.900000 \\
     \texttt{hs67} &      16 &        17 & $-$1.162119e+03 &  0.000000e+00 &  6.495889e$-$05 &  1.000000 \\
      \texttt{hs7} &       7 &         8 & $-$1.732051e+00 &  9.245062e$-$07 &  1.247472e$-$06 &  1.000000 \\
     \texttt{hs70} &       5 &         6 &  1.877865e$-$01 &  0.000000e+00 &  5.938470e$-$05 &  0.656100 \\
     \texttt{hs71} &      20 &        29 &  1.701402e+01 &  6.474528e$-$06 &  6.121931e$-$07 &  0.585588 \\
     \texttt{hs72} &      55 &        56 &  7.276756e+02 &  8.987068e$-$08 &  1.880063e$-$05 &  0.000015 \\
     \texttt{hs73} &      18 &        19 &  2.989474e+01 &  1.903667e$-$09 &  8.426837e$-$06 &  0.032691 \\
     \texttt{hs74} &      21 &        22 &  5.126498e+03 &  7.848646e$-$06 &  2.976747e$-$06 &  0.122491 \\
     \texttt{hs75} &     181 &       531 &  5.174413e+03 &  5.201418e$-$06 &  7.131113e$-$07 &  0.000250 \\
     \texttt{hs76} &       8 &       476 & $-$4.681819e+00 &  5.280851e$-$07 &  9.252311e$-$05 &  0.387420 \\
    \texttt{hs76i} &       7 &        80 & $-$4.681711e+00 &  2.944135e$-$16 &  5.144605e$-$05 &  0.478297 \\
     \texttt{hs77} &      19 &        21 &  2.415058e$-$01 &  7.872245e$-$06 &  6.733883e$-$07 &  1.000000 \\
     \texttt{hs78} &      20 &        21 & $-$2.919704e+00 &  8.018328e$-$06 &  2.979643e$-$06 &  0.656100 \\
     \texttt{hs79} &      17 &        66 &  7.877677e$-$02 &  7.944258e$-$06 &  5.286831e$-$08 &  1.000000 \\
      \texttt{hs8} &      15 &        48 & $-$1.000000e+00 &  5.982302e$-$06 &  2.948125e$-$10 &  1.000000 \\
     \texttt{hs80} &      17 &        84 &  5.394955e$-$02 &  7.682287e$-$06 &  2.102159e$-$07 &  1.000000 \\
     \texttt{hs81} &      19 &        20 &  5.394951e$-$02 &  8.821814e$-$06 &  2.438733e$-$07 &  0.900000 \\
     \texttt{hs86} &      19 &       375 & $-$3.234871e+01 &  3.029890e$-$06 &  1.877376e$-$06 &  0.052335 \\
     \texttt{hs88} &      33 &        37 &  1.362657e+00 &  2.312826e$-$12 &  3.800417e$-$07 &  0.000645 \\
     \texttt{hs89} &      31 &        66 &  1.362657e+00 &  6.336861e$-$13 &  1.383187e$-$07 &  0.000693 \\
      \texttt{hs9} &       2 &         3 & $-$5.000000e$-$01 &  6.821210e$-$13 &  6.842740e$-$05 &  1.000000 \\
     \texttt{hs90} &      33 &        49 &  1.362657e+00 &  7.690347e$-$12 &  4.771747e$-$06 &  0.000652 \\
     \texttt{hs91} &      34 &        48 &  1.362657e+00 &  8.400233e$-$14 &  3.786938e$-$08 &  0.000660 \\
     \texttt{hs92} &      33 &        44 &  1.362657e+00 &  1.792880e$-$11 &  1.663324e$-$05 &  0.000677 \\
      \texttt{hs93} &      57 &        59 &  1.350759e+00 &  2.145353e$-$10 &  9.788245e$-$05 &  0.011972  \\
     \texttt{hs95} &      26 &        27 &  1.567252e$-$02 &  3.024953e$-$10 &  3.415001e$-$05 &  0.011169 \\
     \texttt{hs96} &      22 &        40 &  1.587067e$-$02 &  8.966170e$-$11 &  6.374995e$-$05 &  0.011070 \\
     \texttt{hs97} &      34 &        60 &  4.071230e+00 &  5.619911e$-$08 &  7.635879e$-$05 &  0.001126 \\
     \texttt{hs98} &      49 &       198 &  4.071243e+00 &  1.332748e$-$08 &  9.266559e$-$05 &  0.001105 \\
     \texttt{hs99} &      18 &        19 & $-$8.310799e+08 &  7.331983e$-$06 &  1.215967e$-$06 &  1.000000 \\
  \end{longtable}

\end{footnotesize}

\begin{footnotesize}
\begin {longtable}{l||r|r||r|r|r|r}
 \caption {CUTEr \texttt{hs} test results, 113 successful cases out of 126   problems with $\beta_\phi = 0.5$.} \label{tab:CUTEr-test-beta-phi-p5} \\
 
\hline

\multicolumn{1}{c}{Problem}    & \multicolumn{1}{c}{\# iter}         &  \multicolumn{1}{c}{  \# $f$}  & 
 \multicolumn{1}{c}{   $f(x^*)$} &  \multicolumn{1}{c}{$v(x^*)$}  &  \multicolumn{1}{c}{ KKT  } &  \multicolumn{1}{c}{  Final $\rho$}   \\  \hline
\hline
\endfirsthead
\endhead

\multicolumn{7}{c}%
{{  \tablename\ \thetable{} -- continued from previous page}} \\
\hline
\multicolumn{1}{c}{Problem}  & \multicolumn{1}{c}{  \# iter}     &  \multicolumn{1}{c}{  \# $f$}  & 
 \multicolumn{1}{c}{   $f(x^*)$} &  \multicolumn{1}{c}{     $v(x^*)$}  &  \multicolumn{1}{c}{     KKT  } &    \multicolumn{1}{c}{      Final $\rho$} 
  \\  \hline
\hline
\endhead

\hline \multicolumn{7}{c}{{Continued on next page}} \\ \hline
\endfoot

\hline\hline
\endlastfoot
      \texttt{hs1} &      24 &        34 &  4.215353e$-$17 &  0.000000e+00 &  2.983621e$-$09 &  1.000000 \\
     \texttt{hs10} &       8 &         9 & $-$1.000001e+00 &  1.551523e$-$06 &  3.360370e$-$06 &  1.000000 \\
    \texttt{hs100} &      12 &        82 &  6.806301e+02 &  1.266571e$-$10 &  2.188096e$-$05 &  0.568586 \\
 \texttt{hs100lnp} &      12 &        54 &  6.806301e+02 &  4.923703e$-$06 &  1.323439e$-$05 &  0.135085 \\
 \texttt{hs100mod} &       7 &        23 &  6.786796e+02 &  1.128043e$-$09 &  2.302428e$-$09 &  0.665922 \\
    \texttt{hs101} &      68 &       231 &  1.809744e+03 &  4.591308e$-$06 &  2.266775e$-$05 &  0.000153 \\
    \texttt{hs104} &      13 &        21 &  4.200000e+00 &  2.658767e$-$10 &  7.616564e$-$08 &  0.718165 \\
    \texttt{hs105} &      54 &      1060 &  1.044612e+03 &  2.560362e$-$08 &  3.225074e$-$05 &  0.003757 \\
    \texttt{hs107} &      27 &        31 &  5.054978e+03 &  7.030538e$-$06 &  3.275745e$-$05 &  0.000761 \\
    \texttt{hs108} &      18 &        38 & $-$6.746727e$-$01 &  0.000000e+00 &  8.059244e$-$05 &  0.656100 \\
    \texttt{hs109} &      91 &       208 &  5.362069e+03 &  5.071941e$-$06 &  7.901648e$-$05 &  0.143136 \\
     \texttt{hs11} &      10 &        11 & $-$8.498464e+00 &  1.123963e$-$08 &  1.512764e$-$07 &  0.251157 \\
    \texttt{hs110} &       3 &         5 & $-$4.577848e+01 &  0.000000e+00 &  6.067373e$-$09 &  1.000000 \\
    \texttt{hs111} &      26 &        42 & $-$4.776110e+01 &  2.671024e$-$06 &  1.824005e$-$06 &  0.044867 \\
 \texttt{hs111lnp} &      29 &        45 & $-$4.776109e+01 &  3.505706e$-$06 &  1.208402e$-$05 &  0.037393 \\
    \texttt{hs112} &      24 &        27 & $-$4.776118e+01 &  9.084188e$-$06 &  3.331079e$-$06 &  0.045371 \\
    \texttt{hs113} &      15 &        16 &  2.430637e+01 &  3.236107e$-$06 &  3.610382e$-$05 &  0.387420 \\
    \texttt{hs117} &      19 &        21 &  3.235008e+01 &  3.046132e$-$06 &  1.097142e$-$05 &  0.007070 \\
    \texttt{hs118} &      24 &        25 &  9.329922e+02 &  7.148723e$-$06 &  1.631322e$-$06 &  0.076232 \\
    \texttt{hs119} &      26 &        27 &  2.448993e+02 &  7.988931e$-$06 &  2.257583e$-$06 &  0.160258 \\
     \texttt{hs12} &       5 &         9 & $-$3.000000e+01 &  3.991066e$-$07 &  2.934243e$-$07 &  1.000000 \\
     \texttt{hs14} &      23 &       119 &  1.393453e+00 &  7.451723e$-$06 &  5.102185e$-$06 &  0.429421 \\
     \texttt{hs15} &      11 &        13 &  3.603798e+02 &  2.109424e$-$13 &  1.423044e$-$08 &  0.004071 \\
     \texttt{hs16} &      24 &        25 &  2.314415e+01 &  7.882727e$-$06 &  7.697258e$-$05 &  0.014781 \\
     \texttt{hs17} &       9 &        11 &  1.000000e+00 &  0.000000e+00 &  1.061723e$-$06 &  0.018248 \\
     \texttt{hs18} &       9 &        10 &  5.000000e+00 &  1.234419e$-$08 &  1.001887e$-$07 &  1.000000 \\
     \texttt{hs19} &      48 &        54 & $-$6.961818e+03 &  3.596275e$-$06 &  2.625654e$-$06 &  0.000594 \\
      \texttt{hs2} &       7 &         9 &  4.941229e+00 &  0.000000e+00 &  4.601831e$-$07 &  1.000000 \\
     \texttt{hs20} &      32 &        33 &  4.019833e+01 &  2.088097e$-$06 &  9.533308e$-$05 &  0.009253 \\
     \texttt{hs21} &       2 &         3 & $-$9.996000e+01 &  4.440892e$-$16 &  4.999500e$-$09 &  1.000000 \\
  \texttt{hs21mod} &      10 &        45 & $-$9.596000e+01 &  0.000000e+00 &  8.731893e$-$11 &  0.189722 \\
     \texttt{hs22} &      18 &       430 &  1.000005e+00 &  0.000000e+00 &  4.011241e$-$05 &  1.000000 \\
     \texttt{hs23} &      16 &        17 &  2.000266e+00 &  0.000000e+00 &  6.251621e$-$05 &  0.150095 \\
     \texttt{hs24} &      16 &        17 & $-$9.998796e$-$01 &  0.000000e+00 &  8.782431e$-$05 &  0.729000 \\
     \texttt{hs25} &       1 &         2 &  3.283500e+01 &  0.000000e+00 &  2.005805e$-$08 &  1.000000 \\
     \texttt{hs26} &      13 &        28 &  2.172765e$-$10 &  4.361497e$-$06 &  2.031413e$-$07 &  1.000000 \\
    \texttt{hs268} &       3 &         4 &  8.608487e$-$06 &  0.000000e+00 &  2.526260e$-$05 &  1.000000 \\
     \texttt{hs27} &       7 &        11 &  4.000000e$-$02 &  1.781389e$-$19 &  5.580537e$-$05 &  1.000000 \\
     \texttt{hs28} &       2 &         3 &  1.117108e$-$13 &  0.000000e+00 &  2.034355e$-$07 &  1.000000 \\
     \texttt{hs29} &       8 &         9 & $-$2.262742e+01 &  6.786021e$-$10 &  1.127105e$-$05 &  0.798638 \\
      \texttt{hs3} &       7 &         8 &  2.338799e$-$04 &  0.000000e+00 &  9.672226e$-$05 &  1.000000 \\
     \texttt{hs30} &      17 &        62 &  1.000015e+00 &  4.979195e$-$06 &  6.160904e$-$05 &  0.900000 \\
     \texttt{hs31} &       8 &        10 &  5.999995e+00 &  7.977267e$-$07 &  1.615447e$-$05 &  0.122758 \\
     \texttt{hs32} &      15 &        16 &  1.001279e+00 &  5.218048e$-$15 &  9.183154e$-$05 &  0.071790 \\
     \texttt{hs33} &       4 &         5 & $-$4.000003e+00 &  3.003592e$-$07 &  5.131525e$-$06 &  0.088629 \\
     \texttt{hs34} &      23 &        28 & $-$8.340328e$-$01 &  6.899262e$-$06 &  2.413293e$-$07 &  0.810000 \\
     \texttt{hs35} &       1 &         2 &  1.111111e$-$01 &  0.000000e+00 &  8.332954e$-$05 &  1.000000 \\
    \texttt{hs35i} &       1 &         2 &  1.111111e$-$01 &  0.000000e+00 &  8.332954e$-$05 &  1.000000 \\
  \texttt{hs35mod} &       3 &         7 &  2.500000e$-$01 &  0.000000e+00 &  1.082126e$-$07 &  0.810000 \\
     \texttt{hs36} &      19 &        20 & $-$3.299992e+03 &  0.000000e+00 &  9.398680e$-$05 &  0.003757 \\
     \texttt{hs37} &       8 &        18 & $-$3.456000e+03 &  0.000000e+00 &  3.881216e$-$05 &  0.006363 \\
     \texttt{hs38} &      38 &        58 &  7.867591e$-$14 &  0.000000e+00 &  1.866336e$-$07 &  1.000000 \\
     \texttt{hs39} &      27 &        28 & $-$1.000008e+00 &  7.670604e$-$06 &  5.813017e$-$06 &  0.757858 \\
   \texttt{hs3mod} &       2 &         3 &  8.026142e$-$14 &  0.000000e+00 &  1.162326e$-$07 &  1.000000 \\
      \texttt{hs4} &       2 &         3 &  2.666667e+00 &  0.000000e+00 &  3.149394e$-$15 &  0.228768 \\
     \texttt{hs40} &      23 &        24 & $-$2.500012e$-$01 &  7.931705e$-$06 &  2.342802e$-$06 &  1.000000 \\
     \texttt{hs41} &      17 &        86 &  1.925925e+00 &  6.326517e$-$06 &  1.242304e$-$06 &  0.516082 \\
     \texttt{hs42} &       4 &        18 &  1.385786e+01 &  2.198841e$-$10 &  3.639101e$-$06 &  0.367839 \\
     \texttt{hs43} &      15 &        17 & $-$4.399977e+01 &  5.769252e$-$10 &  7.534945e$-$05 &  0.325810 \\
     \texttt{hs44} &      21 &        22 & $-$1.500003e+01 &  9.486390e$-$06 &  1.097773e$-$06 &  0.007070 \\
  \texttt{hs44new} &      21 &        22 & $-$1.299990e+01 &  7.165140e$-$06 &  7.888032e$-$06 &  0.042391 \\
     \texttt{hs45} &      18 &        19 &  1.000077e+00 &  0.000000e+00 &  2.427232e$-$05 &  0.531441 \\
     \texttt{hs46} &      21 &        54 &  1.908729e$-$09 &  8.068528e$-$06 &  1.842347e$-$07 &  1.000000 \\
     \texttt{hs47} &      20 &        28 &  1.495011e$-$10 &  6.225418e$-$06 &  3.890209e$-$07 &  0.109419 \\
     \texttt{hs48} &       9 &        10 &  2.756164e$-$20 &  4.143635e$-$06 &  4.604039e$-$10 &  1.000000 \\
     \texttt{hs49} &      17 &       143 &  2.791394e$-$07 &  4.239720e$-$12 &  3.325774e$-$05 &  1.000000 \\
      \texttt{hs5} &       5 &         8 & $-$1.913223e+00 &  0.000000e+00 &  3.368922e$-$08 &  0.590490 \\
     \texttt{hs50} &      11 &        12 &  3.197115e$-$16 &  2.625468e$-$07 &  5.107095e$-$08 &  0.006363 \\
     \texttt{hs51} &       2 &         3 &  5.165478e$-$17 &  4.610553e$-$09 &  9.999748e$-$09 &  1.000000 \\
     \texttt{hs52} &      27 &        28 &  5.326601e+00 &  8.035451e$-$06 &  4.260170e$-$06 &  0.113117 \\
     \texttt{hs53} &      26 &        27 &  4.092982e+00 &  9.247958e$-$06 &  4.663561e$-$06 &  0.144093 \\
     \texttt{hs54} &       6 &         7 & $-$1.561253e$-$01 &  9.712104e$-$11 &  9.544783e$-$05 &  1.000000 \\
     \texttt{hs55} &      22 &        23 &  6.666665e+00 &  9.411365e$-$06 &  1.599698e$-$06 &  0.810000 \\
     \texttt{hs56} &       9 &        11 & $-$3.456007e+00 &  5.250197e$-$06 &  5.892236e$-$05 &  0.478297 \\
     \texttt{hs57} &       1 &         2 &  3.064627e$-$02 &  0.000000e+00 &  2.696159e$-$06 &  1.000000 \\
     \texttt{hs59} &      21 &        60 & $-$7.802789e+00 &  5.169181e$-$08 &  4.257062e$-$06 &  1.000000 \\
      \texttt{hs6} &       9 &        24 &  8.091820e$-$10 &  1.605262e$-$07 &  2.842575e$-$05 &  1.000000 \\
     \texttt{hs60} &       5 &         6 &  3.256820e$-$02 &  9.958889e$-$08 &  1.894968e$-$07 &  1.000000 \\
     \texttt{hs61} &      17 &       110 & $-$1.436461e+02 &  5.755363e$-$06 &  2.279362e$-$06 &  0.446153 \\
     \texttt{hs62} &       5 &         7 & $-$2.627251e+04 &  1.526557e$-$16 &  4.647296e$-$07 &  0.001456 \\
     \texttt{hs63} &      19 &        51 &  9.617152e+02 &  5.037713e$-$06 &  2.648938e$-$06 &  0.430467 \\
     \texttt{hs64} &      45 &        46 &  6.299843e+03 &  2.414132e$-$08 &  8.508737e$-$05 &  0.048452 \\
     \texttt{hs65} &       6 &         7 &  9.535284e$-$01 &  5.332807e$-$06 &  4.381076e$-$07 &  1.000000 \\
     \texttt{hs66} &       6 &        44 &  5.181617e$-$01 &  2.335037e$-$06 &  1.396659e$-$06 &  0.900000 \\
     \texttt{hs67} &      17 &        18 & $-$1.162119e+03 &  6.817694e$-$06 &  6.292961e$-$06 &  1.000000 \\
     \texttt{hs69} &     117 &       209 & $-$9.567129e+02 &  4.448109e$-$13 &  9.574185e$-$05 &  0.000002 \\
      \texttt{hs7} &       7 &         8 & $-$1.732051e+00 &  9.245062e$-$07 &  1.247472e$-$06 &  1.000000 \\
     \texttt{hs70} &      19 &        34 &  1.875462e$-$01 &  0.000000e+00 &  9.767005e$-$05 &  0.656100 \\
     \texttt{hs71} &      26 &        36 &  1.701402e+01 &  7.582508e$-$06 &  6.233309e$-$07 &  0.645973 \\
     \texttt{hs72} &      76 &        77 &  7.277010e+02 &  7.386719e$-$10 &  1.621663e$-$05 &  0.000018 \\
     \texttt{hs73} &      18 &        21 &  2.989574e+01 &  1.459998e$-$10 &  3.044686e$-$05 &  0.030903 \\
     \texttt{hs74} &      29 &        30 &  5.126498e+03 &  5.202308e$-$06 &  1.885142e$-$06 &  0.137573 \\
     \texttt{hs76} &      17 &       325 & $-$4.681787e+00 &  2.220446e$-$16 &  4.714559e$-$05 &  0.430467 \\
    \texttt{hs76i} &      22 &       437 & $-$4.681822e+00 &  2.293698e$-$06 &  8.936387e$-$05 &  0.387420 \\
     \texttt{hs77} &      23 &        26 &  2.415057e$-$01 &  6.573382e$-$06 &  5.622722e$-$07 &  1.000000 \\
     \texttt{hs78} &      24 &        25 & $-$2.919703e+00 &  7.440837e$-$06 &  3.113752e$-$06 &  0.729000 \\
     \texttt{hs79} &      19 &        31 &  7.877664e$-$02 &  6.433766e$-$06 &  1.264371e$-$07 &  1.000000 \\
      \texttt{hs8} &      17 &        19 & $-$1.000000e+00 &  8.877931e$-$06 &  4.375083e$-$10 &  1.000000 \\
     \texttt{hs80} &      19 &        20 &  5.394949e$-$02 &  9.291850e$-$06 &  3.046230e$-$07 &  1.000000 \\
     \texttt{hs81} &      23 &        24 &  5.394956e$-$02 &  7.427161e$-$06 &  2.011085e$-$07 &  0.900000 \\
     \texttt{hs86} &      29 &       302 & $-$3.234877e+01 &  7.952470e$-$06 &  3.592137e$-$06 &  0.038152 \\
     \texttt{hs88} &      43 &        47 &  1.362657e+00 &  2.171298e$-$13 &  6.694802e$-$08 &  0.000701 \\
     \texttt{hs89} &      37 &        72 &  1.362657e+00 &  5.467216e$-$13 &  4.170063e$-$07 &  0.000734 \\
      \texttt{hs9} &       2 &         3 & $-$5.000000e$-$01 &  6.821210e$-$13 &  6.842740e$-$05 &  1.000000 \\
     \texttt{hs90} &      44 &        51 &  1.362657e+00 &  7.816595e$-$13 &  8.366618e$-$05 &  0.000699 \\
     \texttt{hs91} &      47 &        59 &  1.362657e+00 &  1.651768e$-$14 &  1.537274e$-$08 &  0.000694 \\
     \texttt{hs92} &      46 &        57 &  1.362657e+00 &  6.403740e$-$14 &  3.249272e$-$08 &  0.000696 \\
     \texttt{hs93} &      59 &        61 &  1.350760e+02 &  4.464122e$-$10 &  8.911687e$-$05 &  0.011973 \\
     \texttt{hs95} &      28 &        33 &  1.572551e$-$02 &  5.708108e$-$13 &  7.454286e$-$05 &  0.010023 \\
     \texttt{hs96} &      28 &        36 &  1.570997e$-$02 &  8.056444e$-$12 &  2.606842e$-$05 &  0.011554 \\
     \texttt{hs97} &      49 &       114 &  3.135805e+00 &  2.494120e$-$08 &  9.749914e$-$05 &  0.001008 \\
     \texttt{hs98} &      78 &       359 &  3.135808e+00 &  7.729489e$-$09 &  2.828064e$-$05 &  0.000731 \\
     \texttt{hs99} &      23 &        26 & $-$8.310799e+08 &  8.955496e$-$06 &  1.485217e$-$06 &  1.000000 \\
      \end{longtable}

\end{footnotesize}

\begin{footnotesize}
\begin {longtable}{l||r|r||r|r|r|r}
 \caption {CUTEr \texttt{hs} test results, 111 successful cases out of 126   problems with $\beta_\phi = 0.99$.} \label{tab:CUTEr-test-beta-phi-p99} \\
 
\hline

\multicolumn{1}{c}{Problem}    & \multicolumn{1}{c}{\# iter}         &  \multicolumn{1}{c}{  \# $f$}  & 
 \multicolumn{1}{c}{   $f(x^*)$} &  \multicolumn{1}{c}{$v(x^*)$}  &  \multicolumn{1}{c}{ KKT  } &  \multicolumn{1}{c}{  Final $\rho$}   \\  \hline
\hline
\endfirsthead
\endhead

\multicolumn{7}{c}%
{{  \tablename\ \thetable{} -- continued from previous page}} \\
\hline
\multicolumn{1}{c}{Problem}  & \multicolumn{1}{c}{  \# iter}     &  \multicolumn{1}{c}{  \# $f$}  & 
 \multicolumn{1}{c}{   $f(x^*)$} &  \multicolumn{1}{c}{     $v(x^*)$}  &  \multicolumn{1}{c}{     KKT  } &    \multicolumn{1}{c}{      Final $\rho$} 
  \\  \hline
\hline
\endhead

\hline \multicolumn{7}{c}{{Continued on next page}} \\ \hline
\endfoot

\hline\hline
\endlastfoot
      \texttt{hs1} &      24 &        34 &  4.215353e$-$17 &  0.000000e+00 &  2.983621e$-$09 &  1.000000 \\
     \texttt{hs10} &       8 &         9 & $-$1.000001e+00 &  1.454144e$-$06 &  7.004923e$-$05 &  0.964609 \\
    \texttt{hs100} &       7 &        16 &  6.806301e+02 &  3.072384e$-$06 &  3.762049e$-$06 &  0.410794 \\
 \texttt{hs100lnp} &       8 &        28 &  6.806301e+02 &  1.243959e$-$06 &  1.726943e$-$06 &  0.150095 \\
 \texttt{hs100mod} &       7 &        23 &  6.786796e+02 &  1.918734e$-$07 &  1.932158e$-$07 &  0.405803 \\
    \texttt{hs104} &       9 &        14 &  4.200000e+00 &  7.857448e$-$11 &  1.338461e$-$09 &  0.852182 \\
    \texttt{hs105} &      21 &        24 &  1.044612e+03 &  6.440900e$-$09 &  9.183609e$-$05 &  0.005154 \\
    \texttt{hs107} &       8 &        11 &  5.054992e+03 &  4.009292e$-$06 &  2.372680e$-$05 &  0.000813 \\
    \texttt{hs108} &       9 &        19 & $-$8.660181e$-$01 &  8.033984e$-$06 &  3.876695e$-$05 &  1.000000 \\
    \texttt{hs109} &      57 &       145 &  5.362069e+03 &  5.287754e$-$07 &  9.963619e$-$05 &  0.106742 \\
     \texttt{hs11} &       7 &         8 & $-$8.498465e+00 &  2.140915e$-$07 &  2.327186e$-$06 &  0.160679 \\
    \texttt{hs110} &       3 &         5 & $-$4.577848e+01 &  0.000000e+00 &  6.067373e$-$09 &  1.000000 \\
    \texttt{hs111} &      13 &        23 & $-$4.776109e+01 &  1.828854e$-$06 &  8.498378e$-$05 &  0.037562 \\
 \texttt{hs111lnp} &      17 &        31 & $-$4.776110e+01 &  9.677666e$-$07 &  9.685628e$-$05 &  0.036232 \\
    \texttt{hs112} &       9 &        10 & $-$4.776111e+01 &  6.414382e$-$06 &  1.509000e$-$06 &  0.033584 \\
    \texttt{hs113} &       6 &         7 &  2.430621e+01 &  7.842192e$-$06 &  3.828328e$-$06 &  0.430467 \\
    \texttt{hs117} &       5 &        10 &  3.234984e+01 &  2.854480e$-$07 &  9.921332e$-$05 &  0.016423 \\
    \texttt{hs118} &       8 &         9 &  9.329922e+02 &  9.331469e$-$06 &  1.385033e$-$06 &  0.051875 \\
    \texttt{hs119} &      11 &        12 &  2.448995e+02 &  5.010450e$-$06 &  8.440035e$-$07 &  0.105371 \\
     \texttt{hs12} &       5 &         9 & $-$3.000000e+01 &  3.972022e$-$07 &  3.691353e$-$05 &  0.947984 \\
     \texttt{hs14} &       8 &        40 &  1.393453e+00 &  7.453297e$-$06 &  3.473740e$-$06 &  0.292307 \\
     \texttt{hs15} &      20 &        21 &  3.064994e+02 &  3.356833e$-$07 &  7.885639e$-$05 &  0.000901 \\
     \texttt{hs17} &       9 &        10 &  1.000000e+00 &  0.000000e+00 &  1.311337e$-$07 &  0.018248 \\
     \texttt{hs18} &       6 &         7 &  5.000000e+00 &  0.000000e+00 &  4.941352e$-$08 &  1.000000 \\
     \texttt{hs19} &      21 &        31 & $-$6.961814e+03 &  2.065775e$-$07 &  1.789287e$-$05 &  0.000485 \\
      \texttt{hs2} &       7 &         9 &  4.941229e+00 &  0.000000e+00 &  4.601831e$-$07 &  1.000000 \\
     \texttt{hs20} &      16 &        17 &  4.019853e+01 &  1.032912e$-$06 &  9.289010e$-$05 &  0.006655 \\
     \texttt{hs21} &       2 &         3 & $-$9.996000e+01 &  0.000000e+00 &  1.059788e$-$07 &  1.000000 \\
  \texttt{hs21mod} &       6 &         7 & $-$9.596000e+01 &  2.220446e$-$16 &  1.217378e$-$18 &  0.137065 \\
     \texttt{hs22} &       4 &         5 &  1.000000e+00 &  1.463081e$-$09 &  1.174890e$-$05 &  0.812471 \\
     \texttt{hs23} &       9 &        10 &  1.999983e+00 &  8.315618e$-$06 &  9.354214e$-$06 &  0.245995 \\
     \texttt{hs24} &       8 &         9 & $-$9.999917e$-$01 &  1.786971e$-$11 &  2.342727e$-$05 &  1.000000 \\
     \texttt{hs25} &       1 &         2 &  3.283500e+01 &  0.000000e+00 &  2.005805e$-$08 &  1.000000 \\
     \texttt{hs26} &      13 &        28 &  2.172765e$-$10 &  4.361497e$-$06 &  2.031413e$-$07 &  1.000000 \\
    \texttt{hs268} &       3 &         4 &  8.608487e$-$06 &  0.000000e+00 &  2.526260e$-$05 &  1.000000 \\
     \texttt{hs27} &       7 &        11 &  4.000000e$-$02 &  1.781389e$-$19 &  5.580537e$-$05 &  1.000000 \\
     \texttt{hs28} &       2 &         3 &  1.117108e$-$13 &  0.000000e+00 &  2.034355e$-$07 &  1.000000 \\
     \texttt{hs29} &       8 &         9 & $-$2.262742e+01 &  1.387226e$-$09 &  1.643885e$-$05 &  0.509159 \\
      \texttt{hs3} &       7 &         8 &  2.338799e$-$04 &  0.000000e+00 &  9.672226e$-$05 &  1.000000 \\
     \texttt{hs30} &       9 &        11 &  1.000002e+00 &  5.926657e$-$06 &  1.713723e$-$05 &  0.900000 \\
     \texttt{hs31} &       6 &         8 &  6.000002e+00 &  0.000000e+00 &  3.527147e$-$06 &  0.078262 \\
     \texttt{hs32} &       3 &         4 &  1.000217e+00 &  8.459899e$-$14 &  4.972783e$-$05 &  0.228768 \\
     \texttt{hs33} &       4 &         6 & $-$3.999996e+00 &  2.597922e$-$10 &  3.209629e$-$06 &  0.088629 \\
     \texttt{hs34} &       8 &        15 & $-$8.340328e$-$01 &  8.697522e$-$06 &  3.382240e$-$07 &  0.900000 \\
     \texttt{hs35} &       1 &         2 &  1.111111e$-$01 &  0.000000e+00 &  8.332954e$-$05 &  1.000000 \\
    \texttt{hs35i} &       1 &         2 &  1.111111e$-$01 &  0.000000e+00 &  8.332954e$-$05 &  1.000000 \\
  \texttt{hs35mod} &       2 &         3 &  2.500037e$-$01 &  1.110223e$-$16 &  6.582735e$-$06 &  0.900000 \\
     \texttt{hs36} &       9 &        10 & $-$3.300001e+03 &  9.869097e$-$06 &  9.110111e$-$06 &  0.005726 \\
     \texttt{hs37} &       7 &         8 & $-$3.456000e+03 &  1.595168e$-$11 &  9.119249e$-$05 &  0.006363 \\
     \texttt{hs38} &      37 &        58 &  6.716615e$-$09 &  0.000000e+00 &  9.209631e$-$05 &  1.000000 \\
     \texttt{hs39} &      10 &        11 & $-$1.000010e+00 &  9.589344e$-$06 &  5.495945e$-$06 &  0.573155 \\
   \texttt{hs3mod} &       2 &         3 &  8.026142e$-$14 &  0.000000e+00 &  1.162326e$-$07 &  1.000000 \\
      \texttt{hs4} &       2 &         3 &  2.666667e+00 &  0.000000e+00 &  3.149394e$-$15 &  0.228768 \\
     \texttt{hs40} &       9 &        10 & $-$2.500017e$-$01 &  7.450609e$-$06 &  2.239323e$-$06 &  1.000000 \\
     \texttt{hs41} &       6 &         7 &  1.925926e+00 &  1.833479e$-$06 &  6.510768e$-$05 &  0.327095 \\
     \texttt{hs42} &       4 &        18 &  1.385786e+01 &  2.149501e$-$10 &  4.274437e$-$05 &  0.234510 \\
     \texttt{hs43} &      10 &        12 & $-$4.399999e+01 &  1.159501e$-$10 &  2.722291e$-$06 &  0.313363 \\
     \texttt{hs44} &      12 &        13 & $-$1.500007e+01 &  9.977008e$-$06 &  9.552142e$-$06 &  0.109419 \\
  \texttt{hs44new} &      13 &        14 & $-$1.500005e+01 &  6.267152e$-$06 &  6.000245e$-$06 &  0.109419 \\
     \texttt{hs45} &      14 &        15 &  1.000002e+00 &  0.000000e+00 &  6.280274e$-$07 &  0.506331 \\
     \texttt{hs46} &      15 &        18 &  1.940719e$-$08 &  6.448602e$-$06 &  1.154214e$-$06 &  1.000000 \\
     \texttt{hs47} &      15 &        20 &  4.310114e$-$09 &  2.798472e$-$06 &  6.038699e$-$06 &  0.313811 \\
     \texttt{hs48} &       4 &         5 &  8.943420e$-$21 &  5.673768e$-$06 &  7.563741e$-$10 &  1.000000 \\
     \texttt{hs49} &      13 &        14 &  2.791394e$-$07 &  7.657430e$-$12 &  3.325774e$-$05 &  1.000000 \\
      \texttt{hs5} &       5 &         8 & $-$1.913223e+00 &  0.000000e+00 &  4.292412e$-$06 &  0.387420 \\
     \texttt{hs50} &       8 &         9 &  2.141999e$-$11 &  4.935936e$-$06 &  5.865571e$-$06 &  0.150095 \\
     \texttt{hs51} &       2 &         3 &  6.496671e$-$17 &  1.204511e$-$08 &  9.999249e$-$09 &  1.000000 \\
     \texttt{hs52} &      10 &        11 &  5.326603e+00 &  7.782006e$-$06 &  2.944367e$-$06 &  0.082081 \\
     \texttt{hs53} &       9 &        10 &  4.092983e+00 &  9.355034e$-$06 &  2.718295e$-$06 &  0.086098 \\
     \texttt{hs54} &       6 &         7 & $-$1.561253e$-$01 &  8.805858e$-$10 &  9.544783e$-$05 &  1.000000 \\
     \texttt{hs55} &       8 &         9 &  6.666669e+00 &  6.324330e$-$06 &  1.964817e$-$06 &  0.702330 \\
     \texttt{hs56} &       7 &         9 & $-$3.456003e+00 &  2.522962e$-$06 &  2.063936e$-$05 &  0.383819 \\
     \texttt{hs57} &       1 &         2 &  3.064627e$-$02 &  0.000000e+00 &  2.696159e$-$06 &  1.000000 \\
     \texttt{hs59} &      12 &        28 & $-$7.802789e+00 &  1.375042e$-$09 &  7.514656e$-$07 &  1.000000 \\
      \texttt{hs6} &       9 &        24 &  8.091820e$-$10 &  1.605262e$-$07 &  2.842575e$-$05 &  1.000000 \\
     \texttt{hs60} &       5 &         6 &  3.256820e$-$02 &  9.958889e$-$08 &  1.894968e$-$07 &  1.000000 \\
     \texttt{hs61} &      15 &       250 & $-$1.436461e+02 &  1.592513e$-$06 &  3.957012e$-$07 &  0.282333 \\
     \texttt{hs62} &       5 &         7 & $-$2.627251e+04 &  1.526557e$-$16 &  4.647296e$-$07 &  0.001456 \\
     \texttt{hs63} &       8 &        11 &  9.617152e+02 &  7.274287e$-$06 &  1.803124e$-$05 &  0.380778 \\
     \texttt{hs64} &      43 &        44 &  6.299843e+03 &  9.138578e$-$09 &  8.249539e$-$05 &  0.031639 \\
     \texttt{hs65} &       5 &         6 &  9.535288e$-$01 &  3.728086e$-$07 &  1.739860e$-$07 &  1.000000 \\
     \texttt{hs66} &       3 &         8 &  5.181609e$-$01 &  6.792149e$-$06 &  1.674698e$-$06 &  1.000000 \\
     \texttt{hs67} &      16 &        17 & $-$1.162119e+03 &  0.000000e+00 &  3.710131e$-$07 &  1.000000 \\
      \texttt{hs7} &       7 &         8 & $-$1.732051e+00 &  9.245062e$-$07 &  1.247472e$-$06 &  1.000000 \\
     \texttt{hs70} &      18 &        23 &  1.875514e$-$01 &  0.000000e+00 &  9.833571e$-$05 &  0.656100 \\
     \texttt{hs71} &      10 &        17 &  1.701402e+01 &  5.412577e$-$06 &  2.626199e$-$05 &  0.387655 \\
     \texttt{hs72} &      37 &        38 &  7.276793e+02 &  1.050624e$-$09 &  1.520144e$-$05 &  0.000015 \\
     \texttt{hs73} &       5 &         6 &  2.989515e+01 &  9.876709e$-$07 &  1.702878e$-$05 &  0.029923 \\
     \texttt{hs74} &      14 &        15 &  5.126498e+03 &  1.689201e$-$06 &  5.878016e$-$05 &  0.093021 \\
     \texttt{hs75} &     138 &       482 &  5.174413e+03 &  8.891761e$-$06 &  2.216519e$-$05 &  0.000283 \\
     \texttt{hs76} &       3 &         4 & $-$4.681787e+00 &  3.049008e$-$16 &  1.646295e$-$05 &  0.531441 \\
    \texttt{hs76i} &       2 &         3 & $-$4.681771e+00 &  5.551115e$-$16 &  2.528737e$-$05 &  0.531441 \\
     \texttt{hs77} &      10 &        12 &  2.415047e$-$01 &  5.346488e$-$06 &  2.731655e$-$06 &  1.000000 \\
     \texttt{hs78} &       8 &         9 & $-$2.919696e+00 &  6.902499e$-$06 &  1.418866e$-$06 &  0.329716 \\
     \texttt{hs79} &       8 &        10 &  7.877678e$-$02 &  4.464599e$-$06 &  3.741699e$-$08 &  1.000000 \\
      \texttt{hs8} &       7 &         8 & $-$1.000000e+00 &  6.707219e$-$06 &  3.305371e$-$10 &  1.000000 \\
     \texttt{hs80} &       8 &        11 &  5.394964e$-$02 &  5.273635e$-$06 &  1.839952e$-$07 &  1.000000 \\
     \texttt{hs81} &       8 &         9 &  5.394952e$-$02 &  8.630110e$-$06 &  2.148511e$-$07 &  0.900000 \\
     \texttt{hs86} &       3 &         4 & $-$3.234849e+01 &  4.202104e$-$06 &  1.382075e$-$05 &  0.071790 \\
     \texttt{hs88} &      25 &        39 &  1.362657e+00 &  2.457157e$-$14 &  2.541669e$-$08 &  0.000656 \\
     \texttt{hs89} &      25 &        60 &  1.362657e+00 &  1.128935e$-$13 &  5.610106e$-$08 &  0.000678 \\
      \texttt{hs9} &       2 &         3 & $-$5.000000e$-$01 &  6.821210e$-$13 &  6.842740e$-$05 &  1.000000 \\
     \texttt{hs90} &      25 &        43 &  1.362657e+00 &  2.226114e$-$14 &  2.424408e$-$08 &  0.000655 \\
     \texttt{hs91} &      26 &        51 &  1.362657e+00 &  2.715207e$-$11 &  6.579864e$-$06 &  0.000603 \\
     \texttt{hs92} &      25 &        36 &  1.362657e+00 &  2.460974e$-$14 &  2.543059e$-$08 &  0.000656 \\
     \texttt{hs93} &      64 &        66 &  1.350760e+02 &  4.470868e$-$11 &  9.377169e$-$05 &  0.010797 \\
     \texttt{hs95} &      30 &       155 &  1.561995e$-$02 &  2.167244e$-$11 &  1.255226e$-$06 &  0.006987 \\
     \texttt{hs96} &      15 &        22 &  1.571116e$-$02 &  5.942702e$-$09 &  1.456985e$-$05 &  0.007057 \\
     \texttt{hs97} &      23 &        45 &  4.071230e+00 &  5.430853e$-$08 &  7.079277e$-$05 &  0.000927 \\
     \texttt{hs98} &      23 &        41 &  4.071231e+00 &  5.204100e$-$08 &  7.454991e$-$05 &  0.000917 \\
     \texttt{hs99} &       8 &         9 & $-$8.310799e+08 &  5.530622e$-$06 &  9.172217e$-$07 &  1.000000 \\
         \end{longtable}

\end{footnotesize}

\begin{footnotesize}
\begin {longtable}{l||r|r||r|r|r|r}

\caption {Infeasible CUTEr \texttt{hs}  test results, 116 successful cases out of 126 infeasible  problems.}\label{tab:CUTEr-inf-test}\\
 
\hline
\multicolumn{1}{c}{Problem}    & \multicolumn{1}{c}{\# iter}         &  \multicolumn{1}{c}{  \# $f$}  & 
 \multicolumn{1}{c}{   $f(x^*)$} &  \multicolumn{1}{c}{$v(x^*)$}  &  \multicolumn{1}{c}{ KKT  } &  \multicolumn{1}{c}{  Final $\rho$}   \\  \hline
\hline
\endfirsthead
\endhead

\multicolumn{7}{c}%
{{  \tablename\ \thetable{} -- continued from previous page}} \\
\hline
\multicolumn{1}{c}{Problem}  & \multicolumn{1}{c}{  \# iter}     &  \multicolumn{1}{c}{  \# $f$}  & 
 \multicolumn{1}{c}{   $f(x^*)$} &  \multicolumn{1}{c}{ $v(x^*)$}  &  \multicolumn{1}{c}{  KKT  } &    \multicolumn{1}{c}{      Final $\rho$} 
  \\  \hline
\hline
\endhead

\hline \multicolumn{7}{c}{{Continued on next page}} \\ \hline
\endfoot

\hline\hline
\endlastfoot

  \texttt{hs100\_inf} &       1 &         2 &  7.050369e+02 &  1.000000e+00 &  2.737830e$-$06 &  1.368915e$-$07 \\ 
  \texttt{hs100lnp\_inf} &      12 &        27 &  6.962388e+02 &  1.000322e+00 &  1.038439e$-$08 &  9.697737e$-$03 \\
   \texttt{hs100mod\_inf} &       1 &         2 &  7.050369e+02 &  1.000000e+00 &  2.737830e$-$06 &  1.368915e$-$07 \\
  \texttt{hs101\_inf} &      41 &       117 &  2.962431e+03 &  1.137023e+00 &  1.519645e$-$02 &  1.213144e$-$06 \\
  \texttt{hs102\_inf} &      45 &       103 &  2.998736e+03 &  1.031044e+00 &  5.193474e$-$03 &  3.135696e$-$05 \\
  \texttt{hs103\_inf} &      46 &       131 &  2.786811e+03 &  1.000003e+00 &  3.274512e$-$06 &  1.955742e$-$05 \\
  \texttt{hs104\_inf} &      22 &        33 &  4.200000e+00 &  1.011048e+00 &  1.779853e$-$06 &  3.311568e$-$03 \\
  \texttt{hs105\_inf} &       2 &         5 &  1.170198e+03 &  1.000000e+00 &  0.000000e+00 &  3.866220e$-$08 \\
  \texttt{hs106\_inf} &      55 &        59 &  2.000000e+03 &  2.275122e+00 &  9.272761e$-$04 &  1.642320e$-$02 \\
  \texttt{hs107\_inf} &      29 &        55 &  5.055009e+03 &  1.000001e+00 &  4.884375e$-$07 &  1.559948e$-$04 \\
  \texttt{hs108\_inf} &      20 &        25 & $-$8.660197e$-$01 &  1.000002e+00 &  1.321342e$-$06 &  5.000000e$-$01 \\
  \texttt{hs109\_inf} &      64 &       140 &  4.777502e+03 &  3.656297e+02 &  3.723947e$-$02 &  6.828624e$-$03 \\
   \texttt{hs10\_inf} &      10 &        45 & $-$1.000000e+00 &  1.000000e+00 &  4.253509e$-$11 &  5.000000e$-$01 \\
  \texttt{hs111\_inf} &      31 &        42 & $-$4.529062e+01 &  1.000000e+00 &  1.557059e$-$07 &  4.458838e$-$02 \\
   \texttt{hs111lnp\_inf} &      23 &       101 & $-$4.529060e+01 &  1.000001e+00 &  1.103729e$-$07 &  1.079561e$-$02 \\
  \texttt{hs112\_inf} &      24 &        25 & $-$4.776110e+01 &  1.000002e+00 &  1.159142e$-$06 &  4.198856e$-$02 \\
  \texttt{hs113\_inf} &      22 &        48 &  4.240308e+01 &  1.000000e+00 &  1.803741e$-$06 &  3.380110e$-$02 \\
  \texttt{hs114\_inf} &     129 &       137 & $-$8.509095e+02 &  6.615003e+02 &  1.557715e+00 &  2.503156e$-$02 \\
  \texttt{hs117\_inf} &       1 &         2 &  2.398758e+03 &  1.000000e+00 &  8.926099e$-$06 &  1.368915e$-$07 \\
  \texttt{hs118\_inf} &      23 &        24 &  9.094002e+02 &  1.000002e+00 &  7.007491e$-$07 &  2.084545e$-$02 \\
  \texttt{hs119\_inf} &      16 &        17 &  2.450109e+02 &  1.000057e+00 &  1.161125e$-$08 &  2.431977e$-$08 \\
   \texttt{hs11\_inf} &       9 &        10 & $-$7.998667e+00 &  1.000000e+00 &  2.613585e$-$08 &  1.902388e$-$01 \\
   \texttt{hs12\_inf} &       1 &         2 & $-$7.451564e$-$02 &  1.000000e+00 &  0.000000e+00 &  1.368915e$-$07 \\
   \texttt{hs13\_inf} &       5 &         6 &  1.932833e+00 &  1.000000e+00 &  0.000000e+00 &  1.215767e$-$01 \\
   \texttt{hs14\_inf} &      24 &        55 &  1.393462e+00 &  1.000002e+00 &  1.634047e$-$06 &  3.695115e$-$01 \\
   \texttt{hs15\_inf} &       4 &         5 &  4.657848e$-$01 &  1.000000e+00 &  1.668519e$-$06 &  9.697737e$-$03 \\
   \texttt{hs16\_inf} &       6 &         7 &  6.417104e$-$01 &  1.000000e+00 &  0.000000e+00 &  2.906321e$-$06 \\
   \texttt{hs17\_inf} &      33 &       829 &  2.052191e+00 &  1.000000e+00 &  8.590665e$-$07 &  1.405123e$-$10 \\
   \texttt{hs18\_inf} &      55 &        85 &  1.568117e+02 &  2.000000e+00 &  6.661338e$-$16 &  9.740833e$-$04 \\
    \texttt{hs1\_inf} &       4 &         5 &  1.258025e+01 &  1.000000e+00 &  0.000000e+00 &  8.709974e$-$10 \\
   \texttt{hs20\_inf} &      11 &        12 &  1.178896e+02 &  1.000000e+00 &  0.000000e+00 &  1.083106e$-$09 \\
   \texttt{hs21\_inf} &       3 &        68 & $-$9.999000e+01 &  1.000000e+00 &  0.000000e+00 &  7.504732e$-$02 \\
  \texttt{hs21mod\_inf} &       7 &         8 & $-$9.598670e+01 &  1.000000e+00 &  1.351003e$-$14 &  1.781441e$-$01 \\
   \texttt{hs22\_inf} &       5 &        40 &  1.000009e+00 &  1.000000e+00 &  1.623803e$-$08 &  5.000000e$-$01 \\
   \texttt{hs23\_inf} &      18 &        20 &  1.999963e+00 &  1.000019e+00 &  1.664523e$-$05 &  3.109002e$-$01 \\
   \texttt{hs24\_inf} &       1 &         2 & $-$1.336948e$-$02 &  1.000000e+00 &  1.626221e$-$09 &  1.368915e$-$07 \\
   \texttt{hs25\_inf} &      10 &        13 &  3.283500e+01 &  1.000000e+00 &  7.749357e$-$14 &  1.000000e+00 \\
  \texttt{hs268\_inf} &       2 &         4 &  3.180734e+03 &  1.000000e+00 &  1.412204e$-$13 &  1.873928e$-$14 \\
   \texttt{hs26\_inf} &      13 &        36 &  6.949105e$-$02 &  1.000000e+00 &  4.654743e$-$10 &  1.213261e$-$08 \\
   \texttt{hs27\_inf} &       5 &         7 &  9.265421e$-$02 &  2.000000e+00 &  3.150539e$-$14 &  1.232023e$-$07 \\
   \texttt{hs28\_inf} &       2 &         3 &  5.001065e$-$08 &  1.000000e+00 &  8.255296e$-$11 &  1.368915e$-$07 \\
   \texttt{hs29\_inf} &       1 &         2 & $-$1.001498e+00 &  1.000000e+00 &  1.567583e$-$07 &  1.368915e$-$07 \\
    \texttt{hs2\_inf} &      11 &        47 &  2.500000e+01 &  1.000000e+00 &  1.756152e$-$10 &  1.901548e$-$12 \\
   \texttt{hs30\_inf} &       1 &         2 &  2.990881e+00 &  1.000000e+00 &  0.000000e+00 &  1.368915e$-$07 \\
   \texttt{hs31\_inf} &      28 &       527 &  1.873971e+01 &  1.000000e+00 &  8.464849e$-$07 &  6.844574e$-$08 \\
   \texttt{hs32\_inf} &       1 &         2 &  7.049285e+00 &  1.000000e+00 &  1.216813e$-$06 &  1.368915e$-$07 \\
   \texttt{hs33\_inf} &       1 &         2 & $-$3.000761e+00 &  1.000000e+00 &  1.673118e$-$06 &  1.368915e$-$07 \\
   \texttt{hs34\_inf} &       1 &         2 & $-$7.605082e$-$04 &  1.000000e+00 &  0.000000e+00 &  1.368915e$-$07 \\
   \texttt{hs35\_inf} &       1 &         2 &  2.228000e+00 &  1.000000e+00 &  0.000000e+00 &  1.368915e$-$07 \\
  \texttt{hs35i\_inf} &       1 &         2 &  2.228000e+00 &  1.000000e+00 &  0.000000e+00 &  1.368915e$-$07 \\
    \texttt{hs35mod\_inf} &       1 &         2 &  2.234820e+00 &  1.000000e+00 &  4.563049e$-$07 &  1.368915e$-$07 \\
   \texttt{hs36\_inf} &       6 &         7 & $-$2.634143e+02 &  1.000000e+00 &  1.519082e$-$08 &  1.815123e$-$10 \\
   \texttt{hs37\_inf} &       7 &        76 & $-$3.150257e+02 &  1.000000e+00 &  3.474836e$-$16 &  2.144845e$-$06 \\
   \texttt{hs38\_inf} &       5 &         6 &  2.801498e+01 &  1.000000e+00 &  0.000000e+00 &  3.807341e$-$09 \\
    \texttt{hs3\_inf} &      10 &        59 & $-$2.846398e$-$14 &  1.000000e+00 &  1.907600e$-$10 &  1.000000e+00 \\
 \texttt{hs3mod\_inf} &       2 &         3 &  9.983686e$-$01 &  1.000000e+00 &  0.000000e+00 &  1.368915e$-$07 \\
   \texttt{hs40\_inf} &      23 &        24 & $-$2.500005e$-$01 &  1.000003e+00 &  1.605618e$-$06 &  1.000000e+00 \\
   \texttt{hs41\_inf} &       4 &         5 &  1.950869e+00 &  1.000000e+00 &  9.205041e$-$09 &  6.026930e$-$08 \\
   \texttt{hs42\_inf} &       1 &        13 &  1.400000e+01 &  2.000000e+00 &  7.605085e$-$07 &  1.368915e$-$07 \\
   \texttt{hs43\_inf} &       1 &         2 & $-$4.101071e$-$01 &  1.000000e+00 &  0.000000e+00 &  1.368915e$-$07 \\
   \texttt{hs44\_inf} &       1 &         2 & $-$1.520438e$-$03 &  1.000000e+00 &  1.521016e$-$07 &  1.368915e$-$07 \\
     \texttt{hs44new\_inf} &       1 &         2 & $-$1.002280e+00 &  1.000000e+00 &  0.000000e+00 &  1.368915e$-$07 \\
   \texttt{hs45\_inf} &       3 &         4 &  1.445187e+00 &  1.000000e+00 &  6.341053e$-$11 &  5.143153e$-$10 \\
   \texttt{hs46\_inf} &       2 &         9 &  3.277487e+00 &  1.000000e+00 &  1.715330e$-$07 &  1.475042e$-$09 \\
   \texttt{hs47\_inf} &      12 &        24 &  2.172792e$-$01 &  1.000040e+00 &  4.105068e$-$09 &  1.638935e$-$09 \\
   \texttt{hs48\_inf} &      18 &       178 &  7.951999e$-$06 &  1.000002e+00 &  1.978232e$-$06 &  1.711143e$-$08 \\
   \texttt{hs49\_inf} &       7 &         8 &  7.144037e$-$03 &  1.000001e+00 &  8.551244e$-$09 &  1.368915e$-$07 \\
    \texttt{hs4\_inf} &       2 &         3 &  6.265554e$-$01 &  1.000000e+00 &  1.521016e$-$07 &  1.368915e$-$07 \\
   \texttt{hs50\_inf} &      10 &        11 &  4.392824e$-$05 &  1.011312e+00 &  6.405816e$-$07 &  9.000000e$-$01 \\
   \texttt{hs51\_inf} &      20 &        49 &  1.048914e$-$10 &  1.000006e+00 &  2.581146e$-$06 &  3.125000e$-$02 \\
   \texttt{hs52\_inf} &      22 &        23 &  5.499964e+00 &  1.000006e+00 &  2.523975e$-$06 &  2.257227e$-$02 \\
   \texttt{hs53\_inf} &      22 &        23 &  5.499973e+00 &  1.000005e+00 &  1.457471e$-$06 &  4.416682e$-$02 \\
   \texttt{hs54\_inf} &       6 &         8 & $-$2.282916e$-$90 &  1.000000e+00 &  4.000073e$-$04 &  1.000000e+00 \\
   \texttt{hs55\_inf} &      13 &        14 &  6.666659e+00 &  1.000053e+00 &  4.316596e$-$06 &  1.119802e$-$04 \\
   \texttt{hs56\_inf} &       7 &        13 & $-$2.379995e+00 &  1.000382e+00 &  3.469250e$-$08 &  3.138106e$-$01 \\
   \texttt{hs57\_inf} &       1 &         2 &  5.079755e$-$02 &  4.000000e$-$01 &  4.024329e$-$06 &  1.368915e$-$07 \\
   \texttt{hs59\_inf} &      50 &       125 &  5.534983e+00 &  1.076927e+01 &  4.356665e$-$05 &  7.221198e$-$02 \\
    \texttt{hs5\_inf} &       1 &         2 &  9.905029e$-$01 &  1.000000e+00 &  0.000000e+00 &  1.368915e$-$07 \\
   \texttt{hs60\_inf} &      11 &        43 &  7.914655e$-$02 &  1.000000e+00 &  2.996201e$-$07 &  1.368915e$-$07 \\
   \texttt{hs61\_inf} &      22 &        32 & $-$7.189256e+01 &  2.750853e+00 &  2.527758e$-$07 &  6.438145e$-$11 \\
   \texttt{hs62\_inf} &       2 &        38 & $-$2.569993e+04 &  1.000000e+00 &  1.712570e$-$19 &  6.844574e$-$08 \\
   \texttt{hs63\_inf} &      17 &        23 &  9.681069e+02 &  1.000113e+00 &  1.265109e$-$06 &  1.041756e$-$01 \\
   \texttt{hs65\_inf} &       3 &         4 &  6.400195e+00 &  1.000000e+00 &  3.431161e$-$06 &  4.295800e$-$08 \\
   \texttt{hs66\_inf} &       8 &       118 &  5.183511e$-$01 &  1.000000e+00 &  1.917211e$-$06 &  5.978711e$-$02 \\
   \texttt{hs67\_inf} &      49 &        55 & $-$5.179534e$-$01 &  1.000000e+00 &  0.000000e+00 &  1.930302e$-$08 \\
   \texttt{hs68\_inf} &       4 &         5 & $-$8.692912e$-$01 &  1.010548e+00 &  7.077349e$-$07 &  3.486784e$-$01 \\
   \texttt{hs69\_inf} &      10 &        50 & $-$9.433878e+02 &  1.000000e+00 &  3.979935e$-$08 &  1.933110e$-$08 \\
    \texttt{hs6\_inf} &       9 &        18 &  2.632899e$-$08 &  1.000000e+00 &  2.995589e$-$08 &  1.368915e$-$07 \\
   \texttt{hs70\_inf} &       3 &         4 &  2.133158e+00 &  1.000000e+00 &  5.329050e$-$08 &  5.892729e$-$08 \\
   \texttt{hs71\_inf} &      20 &        27 &  1.684920e+01 &  1.000022e+00 &  2.198781e$-$05 &  6.250000e$-$02 \\
   \texttt{hs73\_inf} &      12 &        13 &  2.996091e+01 &  1.000000e+00 &  2.680150e$-$15 &  1.230964e$-$02 \\
   \texttt{hs74\_inf} &      11 &        12 &  3.555000e+03 &  4.793332e+02 &  1.295432e$-$03 &  1.039186e$-$01 \\
   \texttt{hs75\_inf} &      10 &        11 &  3.555000e+03 &  4.793329e+02 &  9.790964e$-$04 &  1.105956e$-$01 \\
   \texttt{hs76\_inf} &       1 &         2 & $-$1.257980e+00 &  1.000000e+00 &  0.000000e+00 &  1.368915e$-$07 \\
  \texttt{hs76i\_inf} &       1 &         2 & $-$1.257980e+00 &  1.000000e+00 &  0.000000e+00 &  1.368915e$-$07 \\
   \texttt{hs77\_inf} &      53 &       180 &  6.956721e$-$01 &  1.000268e+00 &  2.285362e$-$05 &  1.868347e$-$03 \\
   \texttt{hs79\_inf} &      18 &        19 &  1.504972e$-$01 &  1.000012e+00 &  5.587084e$-$06 &  2.058911e$-$01 \\
    \texttt{hs7\_inf} &       8 &         9 & $-$1.732051e+00 &  1.000000e+00 &  5.038447e$-$09 &  1.000000e+00 \\
   \texttt{hs81\_inf} &      32 &        70 &  1.000008e+00 &  1.000009e+00 &  9.151413e$-$06 &  3.223196e$-$01 \\
   \texttt{hs83\_inf} &      17 &        18 & $-$2.986759e+04 &  1.415960e+01 &  7.181688e$-$07 &  1.049888e$-$03 \\
   \texttt{hs85\_inf} &      47 &        48 &  4.137154e+01 &  9.368586e+06 &  8.784126e$-$06 &  1.102048e$-$03 \\
   \texttt{hs86\_inf} &       1 &         2 &  1.224387e+01 &  1.000000e+00 &  5.627761e$-$06 &  1.368915e$-$07 \\
   \texttt{hs87\_inf} &      10 &        11 &  9.347058e+03 &  1.000208e+00 &  8.220663e$-$07 &  4.212087e$-$03 \\
   \texttt{hs88\_inf} &       5 &         6 &  1.169516e+00 &  1.000399e+00 &  1.011916e$-$07 &  1.232023e$-$07 \\
   \texttt{hs89\_inf} &      22 &       148 &  5.024977e+00 &  1.051266e+00 &  4.016240e$-$08 &  2.248197e$-$09 \\
    \texttt{hs8\_inf} &      13 &        15 & $-$1.000000e+00 &  1.955966e+00 &  3.366090e$-$06 &  1.000000e+00 \\
   \texttt{hs90\_inf} &      13 &        29 &  1.543613e+00 &  1.000000e+00 &  0.000000e+00 &  4.230379e$-$09 \\
   \texttt{hs91\_inf} &      12 &        39 &  1.842669e+00 &  1.000000e+00 &  0.000000e+00 &  3.426607e$-$09 \\
   \texttt{hs92\_inf} &      10 &        78 &  1.823176e+00 &  1.000000e+00 &  0.000000e+00 &  1.368915e$-$07 \\
   \texttt{hs93\_inf} &       5 &        14 &  1.692649e$-$14 &  3.070000e+00 &  5.126290e$-$21 &  2.996907e$-$04 \\
   \texttt{hs95\_inf} &       8 &        11 &  1.292305e+00 &  6.907751e$-$01 &  1.852601e$-$03 &  6.953210e$-$03 \\
   \texttt{hs96\_inf} &       6 &         8 & $-$1.137676e$-$01 &  1.002112e+00 &  2.102580e$-$03 &  1.126420e$-$02 \\
   \texttt{hs97\_inf} &      17 &        26 &  3.654643e+00 &  1.001483e+00 &  8.537328e$-$04 &  2.783380e$-$05 \\
   \texttt{hs98\_inf} &      14 &        23 &  3.779482e+00 &  1.003179e+00 &  2.049170e$-$03 &  3.672281e$-$05 \\
   \texttt{hs99\_inf} &      10 &        13 & $-$8.310797e+08 &  1.004413e+00 &  4.413353e$-$03 &  1.000000e+00 \\
    \texttt{hs9\_inf} &       1 &        19 &  1.439486e$-$10 &  1.000000e+00 &  1.710708e$-$07 &  5.983858e$-$07 \\
  \end{longtable}

\end{footnotesize}

%**************
% Bibliography
%**************
\bibliographystyle{siamplain}
\bibliography{biblio}

\end{document}